\documentclass[11pt,reqno]{amsart}
\usepackage{amsmath,amssymb,amscd,verbatim,color}
\usepackage[backref]{hyperref}

\usepackage{amsfonts}

\usepackage{color}

\usepackage{xcolor}

\newtheorem{lm}{Lemma}[section]
\newtheorem{thm}{Theorem}[section]
\newtheorem{cor}{Corollary}[section]
\newtheorem{re}{Remark}[section]
\newtheorem{prop}{Proposition}[section]
\newtheorem{definition}{Definition}[section]
\allowdisplaybreaks[3]

\renewcommand*{\backref}[1]{}
\renewcommand*{\backrefalt}[4]{\quad \tiny
  \ifcase #1 (\textbf{NOT CITED.})%
  \or    
  \else   
  \fi}

\makeatletter

\def\MRbibitem{\@ifnextchar[\my@lbibitem\my@bibitem}

\def\mybiblabel#1#2{\@biblabel{{\hyperref{http://www.ams.org/mathscinet-getitem?mr=#1}{}{}{#2}}}}

\def\myhyperanchor#1{\Hy@raisedlink{\hyper@anchorstart{cite.#1}\hyper@anchorend}}

\def\my@lbibitem[#1]#2#3#4\par{%
  \item[\mybiblabel{#2}{#1}\myhyperanchor{#3}\hfill]#4%
  \@ifundefined{ifbackrefparscan}{}{\BR@backref{#3}}%
  \if@filesw{\let\protect\noexpand\immediate
    \write\@auxout{\string\bibcite{#3}{#1}}}\fi\ignorespaces%
}

\def\my@bibitem#1#2#3\par{%
  \refstepcounter\@listctr
  \item[\mybiblabel{#1}{\the\value\@listctr}\myhyperanchor{#2}\hfill]#3%
  \@ifundefined{ifbackrefparscan}{}{\BR@backref{#2}}%
  \if@filesw\immediate\write\@auxout
    {\string\bibcite{#2}{\the\value\@listctr}}\fi\ignorespaces%
}

\makeatother

\title[]{}

\title[Smooth foliations without a bunching condition and Belitskii's $C^{1}$ linearization for RDS
]
{Smooth invariant foliations without a bunching condition and Belitskii's $C^{1}$ linearization for random
dynamical systems}

\keywords{
Invariant foliation; invariant distribution; bunching condition; random normal form; cohomological equation.
\\
This work was partially supported by from NSFC (11922105, 12090010,12090013,11971330,12171336,11831012) and National Key R\&D Program of China 2022YFA1005900.
}

\author{ Wenmeng Zhang}
\address[Wenmeng Zhang]
{School of Mathematical Sciences\\
    Chongqing Normal University\\
    Chongqing 401331, P. R. China}
\email[W.~M.~Zhang]{wmzhang@cqnu.edu.cn}

\author{Kening Lu}
\address[Kening Lu]
{School of Mathematics\\
    Sichuan University\\
    Chengdu, Sichuan 610064, P. R. China }
\email[K.~Lu]{keninglu@gmail.com}

\author{Weinian Zhang}
\address[Weinian Zhang]
{School of Mathematics\\
    Sichuan University\\
    Chengdu, Sichuan 610064, P. R. China }
\email[W.~N.~Zhang]{wnzhang@scu.edu.cn}

\begin{document}

\baselineskip 13pt
\begin{abstract}
Smooth linearization is one of the central themes in the study of dynamical systems.
The classical Belitskii's $C^1$ linearization theorem has been widely used in the investigation of dynamical behaviors such as bifurcations, mixing, and chaotic behaviors due to its minimal requirement of partial second order non-resonances and low regularity of systems.
In this article, we revisit Belitskii's $C^1$ linearization theorem by taking an approach based on smooth invariant foliations and study this problem for a larger class of dynamical systems ({\it random dynamical systems}). We assumed that the linearized system satisfies the condition of Multiplicative Ergodic Theorem and the associated Lyapunov exponents satisfy Belitskii's  partial second order non-resonant conditions. We first establish the existence of $C^{1,\beta}$ stable and unstable foliations without assuming the bunching condition for Lyapunov exponents, then prove a $C^{1,\beta}$ linearization theorem of Belitskii type for random dynamical systems. As a result, we show that the classical Belitskii's $C^1$ linearization theorem for a $C^{2}$ diffeomorphism $F$ indeed holds without assuming all eigenspaces of the linear system $DF(0)$ are invariant under the nonlinear system $F$, a requirement previously  imposed by Belitskii in his proof.

\end{abstract}

\maketitle

\parskip 5pt


\setcounter{page}{1} \tableofcontents


\section{Introduction}

Reducing a nonlinear system to the simplest form through a smooth conjugacy  is one of the fundamental themes in the study of dynamical systems. When the simplest form is a linear system, finding such a conjugacy is called the {\it smooth linearization problem}, which has widely been investigated under different frameworks such as local diffeormorphisms near fixed points (\cite{Hart64, KB-Book95, Poincare, Sternberg2}), circle diffeomorphisms (\cite{Arnold, FK-Ann09, Herman-IHES79,KT-Invent09, Moser2, Yoccoz1}), and toral diffeomorphisms (\cite{Gog-CMP17, HK-IHES90, LMM-CMP2, LMM-CMP1, SY-Adv19}).

The study of smooth linearization goes back to Poincar\'e (\cite{Poincare}).
He proved that analytic linearization can be realized for  an analytic diffeomorphism $F(x)=\Lambda x +O(|x|^2)$ in $\mathbb{C}^d$ when the eigenvalues $r_1,...,r_d$ of $\Lambda$ satisfy  the contraction (or expansion) condition:
\[
\max_{i=1,...,d}|r_i|   <1 \;(\text{or} \;\min_{i=1,...,d}|r_i|>1)
\]
and the {\it non-resonant conditions of all orders}:
\begin{align}
r_j\ne r_1^{m_1}\cdots r_d^{m_d}
\label{NR-cond}
\end{align}
for each $j=1,...,d$,
where $m_i\ge 0$ are all possible integers such that
$\sum_{i=1}^dm_i\ge 2$. Siegel (\cite{Siegel2}) replaced Poincar\'e's conditions by a small-divisor (Diophantine) condition:
\[
|r_j-r_1^{m_1} \cdots r_d^{m_d}|\geq
\frac{C}{|m|^\mu}, \quad 1\leq j\leq d, \;m\in \mathbb{Z}^d_+, \;
|m|=m_1+\cdots +m_d\geq 2,
\]
and proved that the analytic diffeomorphism
 is analytically conjugated to its linear part $\Lambda$. Siegel's condition allows $\Lambda$ to be hyperbolic having both contraction and expansion.
The simpler proof of Siegel's theorem was later given by
V. Arnold (\cite{AR}), Moser (\cite{Moser}) and Zehnder (\cite{Zehnder}) by using the KAM method.

$C^k$ smooth linearization in $\mathbb{R}^d$ was firstly obtained by Sternberg (\cite{Sternberg1, Sternberg2}) in 1950s. He proved that for each $k\in \mathbb{N}$ there is a sufficiently large  integer $N>0$
such that a $C^N$ local hyperbolic diffeomorphism $F(x)=\Lambda x +O(|x|^2)$ is $C^k$ conjugate to its linear part $\Lambda$
if the eigenvalues $r_1,\cdots, r_d$ of $\Lambda$ satisfy
 the {\it  non-resonant  condition {\rm (\ref{NR-cond})} up to order $N$}, i.e., (\ref{NR-cond}) holds  for all $m_i\ge 0$ with $2\le \sum_{i=1}^dm_i\le N$. Following this work, there is an extensive literature on this subject, see for example  Nagumo-Is\'e (\cite{Nagumo}), Chen (\cite{Chen}), Nelson (\cite{Nelson}),  Sell (\cite{Sell1}), and   Banyaga, de la Llave, and Wayne (\cite{BLW}).  A crucial condition in these works is the {\it high order non-resonant} condition. In the meantime, it also requires {\it high order regularity} of systems.

A fundamental question is whether a smooth linearization can be realized with the lowest order non-resonance, or even a part of the lowest order non-resonance (we call it the {\it partial lowest order non-resonance}), together with low order regularity of the systems. Without assuming any non-resonant conditions, Hartman (\cite{HartPAMS60}) and Grobman (\cite{Grobman}) independently showed in 1960's that a $C^1$ local hyperbolic diffeomorphism  can be $C^0$ linearized. However, the $C^0$ conjugacy usually does not preserve dynamical properties
such as the characteristic directions, the derivatives of the Poincar\'e mapping of global orbits,
the iteration rates of systems and the differentiable structure of smooth manifolds. In order to maintain such properties, the $C^1$ conjugacy is essential.

In 1960, Hartman (\cite{Hartman60}) proved that
a $C^{1,1}$ {\it contractive} diffeomorphism in $\mathbb{R}^d$ admits local $C^1$ linearization without assuming any non-resonances. In the same paper, he also proved that a $C^{1,1}$ typically hyperbolic diffeomorphism admits local $C^1$ linearization under the following condition:
\begin{align}
|r_1|/|r_{t}|<|r_{{t}+1}|^{-1}, \qquad
|r_{t+1}|/|r_d|<|r_{{t}}|,
\label{bunching}
\end{align} where the eigenvalues $r_1,\cdots, r_d$ are ordered  as
\begin{align}
|r_1|\ge\cdots\ge|r_{t}|>1>|r_{{t}+1}|\ge\cdots\ge |r_d|
\label{Hyb-cond}
\end{align}
for
$
{t}\in \{1,...,d-1\}
$.
Condition \eqref{bunching} means that the contractive spectrum and the expansive spectrum both lie in a ``narrow-band" area in $\mathbb{C}$, which is referred to as the {\it bunching condition}
(see for example \cite[pp.603-604]{KB-Book95}.
Since the bunching condition (\ref{bunching}) holds automatically in $\mathbb{R}^2$, as a result,  a $C^{1,1}$ typically hyperbolic diffeomorphism in $\mathbb{R}^2$ admits local $C^1$ linearization. Samovol (\cite{Sam}) in 1972 proved that a $C^N$ typically hyperbolic diffeomorphism with sufficiently large $N\in\mathbb{N}$
admits local $C^1$ linearization if the eigenvalues satisfy a partial second order non-resonant condition:
\begin{align}
r_j\ne r_ir_\kappa
\label{S-cond}
\end{align}
for $i=1,...,{t}$, $\kappa={t}+1,...,d$ and $j=1,...,d$. Note that condition \eqref{S-cond} is not only of the lowest order but also a part of the second order non-resonance. However, it retains the requirement of $C^N$ smoothness of the system for large $N$.
Samovol's strategy is to transform the $C^N$ system to its polynomial normal form by Sternberg's normal form theorem, which requires $N$ to be sufficiently large, and then to give an explicit form of the conjugacy, which can linearize the polynomial normal form and can be proved to be $C^1$ under (\ref{S-cond}).

In 1973, this high order regularity requirement was removed by Belitskii (\cite{Beliskii1}). 
He stated that
a $C^{1,1}$ typical hyperbolic diffeomorphism $F$ in $\mathbb{R}^d$ admits local $C^1$ linearization if
a condition
\begin{align}
|r_j|\ne |r_i|\,|r_\kappa|
\label{B-cond}
\end{align}
holds for $i=1,...,{t}$, $\kappa={t}+1,...,d$ and $j=1,...,d$. Clearly, this result gives the above Hartman's $C^1$ linearization results
since \eqref{bunching} implies \eqref{B-cond}. We notice that condition \eqref{B-cond} allows some second order resonances to appear. Condition (\ref{B-cond}), called the {\it strong non-resonance} of order 2,
is a little stronger than (\ref{S-cond}), but
counter examples given in \cite[pp.139-142]{Beliskii2} show that such a condition is sharp, i.e., it cannot be removed for $C^1$ linearization of a $C^2$ system in general.

Belitskii's theorem is a truly remarkable and optimal result for $C^1$ linearization, and has
many applications to homoclinic bifurcations (\cite{AAIS-book94,HomSan-10,HomKoNa-ARMA01}),
heteroclinic tangencies (\cite{LaRo-JDE15}),
mixing of hyperbolic systems
(\cite{FMT-Ann07}),
invariant measures of partially hyperbolic systems (\cite{Hom-PAMS17}), and solutions of
higher order semilinear equations (\cite{Budd-SIAM-JMA89, DFG-MathAnn10}).
However,
Belitskii's proof was commented ``In our opinion, the proofs of Theorems given by Belitskii should be recognized as insufficient and
incomprehensible" by Bronstein and Kopanskii (\cite[p.191]{Bronkopa-book94}).
Actually, in his proof, Belitskii needs to establish the surjectivity of nonlinear operators
deduced from the conjugacy equation, where
all eigenspaces of $\Lambda$ are in fact assumed to be invariant under $F(x)=\Lambda x +O(|x|^2)$.
 This is a very restrictive assumption, which will be elaborated just after the following Corollary.



In this article, we revisit Belitskii's $C^1$ linearization theorem, taking a different approach based on the smooth invariant foliations,
which was not involved in Samovol and Belitskii's proofs,
and study this problem for {\it random dynamical systems}, a larger class of dynamical systems. We first establish the existence of $C^{1,\beta}$ stable and unstable foliations without assuming the bunching condition \eqref{bunching} but only the partial second order non-resonant condition \eqref{B-cond}.
Then we prove a Belitskii type of $C^{1,\beta}$ linearization theorem for random dynamical systems under a partial second order non-resonant condition in terms of  Lyapunov exponents.


Let $(\Omega, \mathcal{F}, \mathbb{P})$ be a probability space and
$(\theta^n)_{n\in \mathbb{Z}}$ be a measurable
$\mathbb{P}$-measure preserving dynamical system on $\Omega$. A
random dynamical system (or a cocycle) on the space $\mathbb{R}^d$
over the dynamical system ${\theta}^n$ is a measurable mapping
$$
F: \mathbb{Z} \times \Omega \times \mathbb{R}^d \to
\mathbb{R}^d, \quad (n, \omega, x) \mapsto F(n, \omega, x),
$$
such that the mapping $F(n,\omega,\cdot):
\mathbb{R}^d \to \mathbb{R}^d$ forms a cocycle over $\theta^n$:
\[
F(0, \omega,\cdot)={\rm id}~~\mbox{(the identity mapping)}, \quad \hbox{ for all }\; \omega \in \Omega,
\]
\[
F(n+m,\omega,\cdot)=F(n,\theta^{m}\omega,\cdot)\circ F(m,\omega,\cdot),
\quad \hbox{ for all }\; m, n \in \mathbb{Z}, \quad\omega \in
\Omega.
\]
Replacing $\mathbb{Z}$ by $\mathbb{R}$ gives a continuous time
random dynamical system. A typical example is the solution operator for a stochastic
differential equation (\cite{Arn98}).

We consider a $C^{2, \alpha}$ ($\alpha\in[0.1]$) random dynamical system $F(n, \omega, x)$ (that is, $F$ is
$C^{2,\alpha}$ in $x$ for each $n$ and $\omega$ and its derivatives are measurable with respect to $\omega$),  and write the time-one
mapping $F(1, \omega, x)$ as $F(\omega, x):=F(1, \omega,
x)$. Then $F(\omega, \cdot)$ is a random
diffeomorphism, whose $i$-th order derivative is denoted by $D^iF(\omega,\cdot)$ for every $i\in\mathbb{N}$. This random diffeomorphism generates the random dynamical system
\[
F(n, \omega, \cdot)=
\begin{cases} F(\theta^{n-1}\omega, \cdot)\circ \cdots \circ
F(\omega, \cdot), & n >0, \\
I, & n=0, \\
F^{-1}(\theta^n\omega, \cdot )\circ \cdots \circ
F^{-1}(\theta^{-1}\omega, \cdot), & n<0.
\end{cases}
\]

We assume that
$F(\omega, \cdot)$ has a fixed point $x=0$ for all $\omega\in \Omega$,
to which a large class of random dynamical systems can be converted (see \cite[p.310]{Arn98}),
and is a locally tempered $C^{2,\alpha}$ random diffeomorphism, that is, there is a
tempered ball $V(\omega)=B_{\rho(\omega)}(0)=\{x\in\mathbb{R}^d: \|x\|<
\rho(\omega)\}$, where $\rho(\omega)$ is a random variable
tempered from below (i.e., $\lim_{n\to\pm\infty}\frac1n\min\{0,\log\rho(\theta^n\omega)\}=0,\;
\mathbb{P}-a.s.$), such that
\begin{align*}
\|D^iF(\omega,x)\|\le M(\omega),
\quad\|D^2F(\omega,x)-D^2F(\omega,y)\|\le
L(\omega)\|x-y\|^\alpha, \quad \forall x,y\in V(\omega),
\end{align*}
where $M(\omega),L(\omega)>0$ are random variable tempered from above (i.e., $\lim_{n\to\pm\infty}\frac1n\max\{0,\log M(\theta^n\omega)\}=0,\;
\mathbb{P}-a.s.$). The size of a tempered ball may decrease as
$\omega$ varies, but these changes along each orbit
$\theta^n\omega$ are at a subexponential rate. The upper bounds of $M(\omega)$ and $L(\omega)$ may grow to infinity  as $\omega$ varies, but along
each orbit $\theta^n\omega$, it may increase only at a
subexponential rate. This nonuniform behavior is one of the
intrinsic features of random dynamical systems.

Two local tempered random diffeomporhisms $F$ and
$G$ with fixed point $x=0$ are said to be $C^{1, \beta}$ locally conjugate for $\beta\in [0,1]$
if there exists a $C^{1,\beta}$ random diffeomorphism
$\Phi(\omega,x)$ defined on a tempered ball $V(\omega)$ with
$\Phi(\omega,0)=0$ such that
$$
\Phi(\theta\omega,F(\omega,x))=G(\omega, \Phi(\omega,x))\quad \hbox{for}
\quad x\in V(\omega), \quad \mathbb{P}-a.s.
$$
This conjugacy relation implies that $h$ carries orbits of
$\varphi$ to orbits of $\psi$ when the orbits stay in the corresponding tempered ball.

We write $F(\omega, x)$ as
\[
F(\omega, x)= \Lambda(\omega)x + f(\omega, x)=\Lambda(\omega)x + O(\|x\|^2),
\]
where $\Lambda(\omega)=DF(\omega, 0) \in Gl(d, \mathbb{R})$. Assume that $\Lambda(\omega)$ satisfies the
conditions of the Multiplicative Ergodic Theorem, that is,
\[
\max\{0,\log\|\Lambda(\cdot)\|\} \in L^1(\Omega, \mathcal{F}, \mathbb{P})\quad
\hbox{and} \quad \max\{0,\log\|\Lambda^{-1}(\cdot)\|\}\in L^1(\Omega,
\mathcal{F}, \mathbb{P}),
\]
and that $\mathbb{P}$ is an ergodic invariant measure for $\theta$.
By the Multiplicative Ergodic Theorem
(see Section 2 for details), there exists a $\theta$-invariant set
$\tilde \Omega\subset \Omega$ of full measure such that for each
$\omega \in \tilde \Omega$, $\Lambda(n, \omega)$ (the linear random
dynamical system generated by $\Lambda(\omega)$) has $p$
Lyapunov exponents $\lambda_j$, $j=1, \cdots, p$. We assume that $\Lambda(n, \omega)$ is hyperbolic, namely, there is no zero Lyapunov exponent. We order them as follows:
\[
\lambda_1 >\cdots \lambda_\tau >0> \lambda_{\tau+1}\cdots > \lambda_{p}.
\]
 For the remainder of this article, we use $\Omega$ to denote $\tilde \Omega$ and assume that these statements are true for all $\omega\in \Omega$.

Our main result is the following theorem on $C^{1,\beta}$ foliations and linearization for random dynamical systems in the typically hyperbolic case, in which there are both positive and negative Lyapunov exponents.

\noindent {\bf Theorem.} {\it
Let $F:\Omega \times V(\omega)\to \mathbb{R}^d$ be a tempered $C^{2,\alpha}$ $(\alpha\in[0,1])$ random diffeomorphism with an $\epsilon$-slowly continuous second order derivative when $\alpha=0$.
Assume its linearization $\Lambda(\omega)$ at hyperbolic fixed point $x=0$ satisfies the condition of the Multiplicative Ergodic Theorem. If the Lyapunov exponents satisfy a partial second order non-resonant condition:
\begin{align}\label{Lya-non}
\lambda_j\ne \lambda_i+\lambda_\kappa
\end{align}
for all $i=1,...,{\tau}$, $\kappa={\tau}+1,...,{{p}}$ and $j=1,...,{{p}}$,
then the following assertions are true:
\begin{itemize}
\item[(i)] Foliations:  $F$ has $C^{1,\beta_\alpha}$ random stable and unstable foliations with $C^{2,\alpha}$ leaves in a tempered ball $U(\omega)\subset V(\omega)$, where $\beta_\alpha=0$ when $\alpha=0$ and $\beta_\alpha\in (0,\alpha]$ when $\alpha\in (0,1]$.
\item[(ii)] Conjugacy:
$F(\omega,\cdot)$ is conjugate
to its linear part $\Lambda(\omega)$ in a tempered ball $U(\omega)\subset V(\omega)$ by a $C^{1, \beta_\alpha}$  random diffeomorphism $\Phi:\Omega\times U(\omega)\to \mathbb{R}^d$ satisfying
$
\Phi(\omega,0)=0
$
and
$
D\Phi(\omega,0)={\rm id}.
$
\end{itemize}
}

\noindent{\bf Remark.} {\it  {\rm (1)}
The assumption of $\epsilon$-slowly continuity means that $e^{-\epsilon |n|}D^2 F(\theta^n \omega, \cdot)$ is equicontinuous with respect $n\in \mathbb{Z}$ for each fixed $\omega\in\Omega$. It
is only necessary for random dynamical systems in the case of $\alpha=0$ since, when $\alpha\in (0,1]$ or $F$ is deterministic, the second order derivative of $F(\omega,\cdot)$ is automatically
$\epsilon$-slowly continuous. {\rm (2)} The non-resonant condition \eqref{Lya-non} obviously does not imply the bunching condition:
\begin{equation}\label{R-bunching}
\lambda_1-\lambda_\tau < -\lambda_{\tau+1} \quad\text{and}\quad \lambda_{\tau+1}-\lambda_p < \lambda_{\tau}.
\end{equation}
{\rm (3)} When $F$ is deterministic,
the relationship between  Lyapunov exponent $\lambda$ and its corresponding eigenvalue
$r$ is $|r|=e^\lambda$. Thus condition
\eqref{Lya-non} becomes \eqref{B-cond}.}

As a result, we have the following $C^{1,\beta_\alpha}$ linearization theorem.

\noindent {\bf Corollary.} {\it
Let $F:V\to \mathbb{R}^d$ be a $C^{2,\alpha}$ {\rm ($\alpha \in [0,1]$)}  diffeomorphism over a neighborhood $V\subset \mathbb{R}^d$ of $x=0$. Suppose the eigenvalues of its linearization $\Lambda:=DF(0)$ at the hyperbolic fixed point $x=0$ satisfy  {\rm (\ref{B-cond})}.
Then $F$ is conjugate
to its linear part $\Lambda$ in a small neighborhood $U\subset V$ by a $C^{1,\beta_\alpha}$ diffeomorphism $\Phi: U\to \mathbb{R}^d$ satisfying
$
\Phi(0)=0
$
and
$
D\Phi(0)={\rm id},
$
where $\beta_\alpha=0$ when $\alpha=0$ and $\beta_\alpha\in (0,\alpha]$ when $\alpha\in (0,1]$.
}

When $\alpha=0$, this corollary is exactly the same as the above-mentioned Belitskii's theorem \cite{Beliskii1} except the $C^2$
smoothness, which is stronger than his $C^{1,1}$ smoothness. In \cite{Beliskii1}, Belitskii proved the above result under
an assumption that all eigenspaces of $\Lambda$ are invariant under the nonlinear system $F$.  To be more precise, in \cite[Lemma 3.1]{Beliskii1}, he assumed that
the subspaces $L_I:=\bigoplus_{i\in I\subset Q}L_i$ of $\mathbb{R}^d$ are invariant under $F$ for all subsets $I$ of $Q=\{1,\cdots,q\}$ with $q\in \mathbb{N}$, where each $L_i$ is the eigenspace of $\Lambda$ corresponding to the eigenvalue lying on the circle $\{z\in \mathbb{C}:|z|=\rho_i\}$ and those $\rho_i$, being the moduli of the eigenvalues of $\Lambda$, are ordered as
$$
\rho_1> \rho_2>\cdots >\rho_k >1> \rho_{k+1}>\cdots>\rho_q.
$$
In other words, Belitskii assumed that the invariant manifolds of $F$ associated with any spectral splitting according to their moduli  are all flat, i.e., linear spaces. This can be achieved only if one can construct a $C^{1,1}$ invariant manifold of $F$ for any spectral splitting:
\[
\{\rho_1, \cdots, \rho_q\}=\{\rho_i : i\in I\}\cup \{\rho_i : i\in Q\setminus I\}.
\]
In this case, the invariant manifold is given by the graph of a  $C^{1,1}$ mapping
\[
\phi: L_I=\bigoplus_{i\in I}L_i \mapsto L_{I'}=\bigoplus_{i\in Q\setminus I}L_i
\] Using these manifolds to change the coordinate system, one obtain  flat invariant manifolds for the new system.

By the theory of invariant manifolds, $F$ has $C^{1,1}$ stable manifold and unstable manifold, i.e., according to the spectral slitting:
\[
Q=\{\rho_i : \rho_i<1\} \cup \{\rho_i : \rho_i>1\}.
\] However, in order to have $C^{1,1}$ pseudo-stable and pseudo-unstable manifolds for $F$, a spectral gap condition is needed. de la Llave and Wayne \cite[Example 5.1]{LW-MZ95} gave an example of a 2-dimensional $C^\infty$ diffeomorphism with a small spectral gap, whose pseudo-stable manifold is not $C^{(\ln 3/\ln 2)+\varepsilon}$ for any $\varepsilon>0$, so the pseudo-stable manifold is not $C^{1,1}$ either. To better understand the issue here, we borrow their Example 5.1 and add one trivial stable component to get the diffeomorphism
\begin{equation} \label{example}
F(x)=\begin{pmatrix}
    1/2 & 0 & 0 \\
    0 & 2 & 0 \\
    0 & 0 & 3
  \end{pmatrix}
  \begin{pmatrix}
  x_1 \\
  x_2 \\
  x_3
  \end{pmatrix}
  +
  \begin{pmatrix}
  0 \\
  0 \\
  \varphi(x_2)
  \end{pmatrix},
\end{equation}
where $\varphi$ is a $C^\infty$ function with a compact support. Note that the eigenvalues satisfy the non-resonant condition \eqref{B-cond}, but the pseudo-stable manifold associated with the splitting $\{1/2, 2\}\cup\{3\}$ is not $C^{1,1}$. Consequently,
one can not transform equation \eqref{example} by a conjugacy to a $C^{1,1}$ system with invariant manifolds all being linear spaces. Furthermore, a spectral splitting such as
\[
Q=\{\rho_1, \rho_q\} \cup \{\rho_i : i\in Q, i\not= 1, q\}
\]
may not yield a $C^{1,1}$ invariant manifold (as a graph over $L_1\bigoplus L_q$) either. The existence of such a $C^{1,1}$ invariant manifold may need much stronger non-resonant conditions (see .e.g. \cite{L-JSP97}). In short, it is quite restricted to assume that all the invariant manifolds of $F$ based on spectral decompositions are linear spaces. We emphasize that the non-resonant condition \eqref{B-cond} implies neither the spectral gap condition needed for $C^{1,1}$ invariant manifolds as we see for example \eqref{example} nor the bunching condition \eqref{bunching} required in Hartman's theorem.

On the other hand,
Bronstein and Kopanskii stated in \cite[p.191]{Bronkopa-book94} that ``Theorem 9.4 coincides, in fact, with Theorem 5 (i.e., Belitskii's theorem), but the proof is new." However, we find that there is no ``Theorem 9.4" in their monograph \cite{Bronkopa-book94}. In fact, there are only Theorems 9.1 and 9.6 in \cite[Chapter II, Section 9]{Bronkopa-book94}
for the smooth linearization of mappings. The two theorems are the same type as Sternberg's theorem since $F$ was required to be sufficiently smooth, and the derivatives of $F$ up to sufficiently high order were required to be equal to $0$ at the origin, which is equivalent to a corresponding non-resonant condition. These conditions on  high order smoothness and vanishing derivatives make it possible to get a contractive operator, whose fixed point is used to construct the conjugacy, by the homotopy method (for Sternberg's theorem) or by the conjugacy relation (for Bronstein and Kopanskii's results).

In this article, we do not make any assumption on the invariant manifolds of $F$,  except for the non-resonant condition \eqref{Lya-non} and the $C^{2, \alpha}$ regularity of $F$. We take a different approach from Belitskii's. We construct $C^{1,\beta_\alpha}$ random stable and unstable foliations with  $C^{2,\alpha}$ leaves and utilize them to decouple the system into {\it a contraction and an expansion.}
The $C^{2, \alpha}$  smoothness of the leaves can be obtained by using the Hadamard graph transformation or the Lyapunov-Perron method. However, without the bunching condition on the Lyapunov exponents, it is  more complicated to prove the $C^{1,\beta_\alpha}$ smoothness of the foliations with respect to the base point only under the non-resonant condition \eqref{Lya-non}.

To obtain the $C^{1,\beta_\alpha}$ smoothness, it is sufficient to show that the tangent spaces of leaves of the stable and unstable foliations are $C^{1,\beta_\alpha}$ distributions thanks to Frobenius' theorem. We construct  continuous bases for these stable and unstable distributions which are uniformly $C^{1,\beta_\alpha}$ along
each local leaf of the stable and unstable foliations. This will be done by solving a random cohomological equation along each leaf of all intermediate (pseudo-center) foliations, which is the main step to overcome the difficulty of no bunching condition (see the end of section 7). The non-resonant condition \eqref{Lya-non} plays a key role in establishing the smoothness, which will be used to prove the result on random normal form and to show the invertibility of  cohomological operators (see (\ref{nobunch}) and (\ref{nobunch-2}) below). To deal with the non-uniformity of hyperbolicity, the Lyapunov norm is used in solving the random cohomology equation. Our proof also involves the roughness of non-uniformly exponential dichotomy and  Journ\'{e}'s lemma (\cite{Journe-88}).

When the Lyaponov exponents satisfy the bunching condition (\ref{R-bunching}) instead of \eqref{Lya-non},
the $C^{1,\beta_\alpha}$ smoothness of the random stable and unstable foliations with respect to the base point can be obtained by using the Lyapunov-Perron approach, see for example \cite{SLZ}.  There is an extensive study of regularity of invariant foliations. The optimal regularity of invariant foliations for hyperbolic dynamical systems under bunching conditions was obtained by Hasselblatt (\cite{Hass-ETDS94, Hass2, Hass3}), and Hasselblatt and Wilkinson (\cite{Hass-W}). Moreover,
Push, Shub and Wilkinson (\cite{PSW}) established the regularity of invariant foliations  of partially hyperbolic diffeomorphisms under bunching conditions.

Finally, we mention some existing results which are closely relevant to the results derived in this article. The theorems of random version were established for Poincar\'e's analytic linearization,
Siegel's analytic linearization, Sternberg's smooth linearization
and Hartman-Grobman's $C^0$ linearization (\cite{BaVa,LiLu1,LiLu3,LiLu2,Wanner}),
but no result for random dynamical systems on $C^1$ linearization
until very recently. In \cite{LZZ}, the authors of the current paper  proved  a random version of Hartman's $C^{1,\beta_\alpha}$ linearization
for a $C^{1,\alpha}$ ($\alpha\in (\alpha_0,1]$) random contraction (or expansion) in $\mathbb{R}^d$,
where $\alpha_0\in (0,1)$ is a constant depending on the Lyapunov exponents.


\section{Main result}
\setcounter{equation}{0}

In this section, we first review some of the basic concepts and facts for random dynamical systems taken from Arnold's book \cite{Arn98}. We then state our main results.


Let $(\Omega, \mathcal{F}, \mathbb{P})$ be a probability space with a sample set $\Omega$, a $\sigma$-algebra $\mathcal{F}\subset 2^{\Omega}$ and a probability measure $\mathbb{P}$ on $\mathcal{F}$, and let $(\Omega, \mathcal{F}, \mathbb{P}, (\theta^n)_{n\in \mathbb{Z}})$
(denoted by $\theta$ for short) be a metric dynamical system (see \cite[p.536]{Arn98}).
A  mapping
$$
{F}: \mathbb{Z} \times \Omega \times \mathbb{R}^d \to
\mathbb{R}^d, \quad (t, \omega, x) \mapsto {F}(t, \omega, x),
$$
where $\mathbb{Z}$ is endowed with its Borel $\sigma$-algebra
${\mathcal B}(\mathbb{Z})$,
is called a $C^N$ ($N\ge 1$) {\it random dynamical system} (abbreviated as RDS) on the measurable space
{\rm ($\mathbb{R}^d$, ${\mathcal B}$)} over $\theta$ if
\begin{itemize}
\item[(1)] $F(n,\omega, \cdot)$ is
$\mathcal{B}(\mathbb{Z})\otimes\mathcal{F}\otimes\mathcal{B}(\mathbb{R}^d)$-measurable
and $F(n,\omega, \cdot)$ is $C^N$ with measurable derivatives, and

\item[(2)] the mappings ${F}(n,\omega,\cdot):
\mathbb{R}^d \to \mathbb{R}^d$, $n\in \mathbb{Z}$, form a cocycle over $\theta$, i.e.,
\begin{align*}
&{F}(0, \omega,\cdot)={\rm id}~ \mbox{{\rm(}the identity mapping{\rm )}},   \quad\forall\omega \in \Omega,
\\
&{F}(n+m,\omega,\cdot)={F}(n,\theta^{m}\omega,\cdot)\circ {F}(m,\omega,\cdot),
\quad \forall n, m \in \mathbb{Z}, \quad\forall \omega \in \Omega.
\end{align*}
\end{itemize}
Clearly, the time-one mapping ${F}(\omega, \cdot):={F}(1, \omega, \cdot): \mathbb{R}^d\to \mathbb{R}^d$ is a $C^N$ random diffeomorphism.
Conversely, a $C^N$ random diffeomorphism ${F}(\omega, \cdot): \mathbb{R}^d \to \mathbb{R}^d$
generates a $C^N$ RDS by
\begin{align}
{F}(n, \omega, \cdot):=
\begin{cases} {F}(\theta^{n-1}\omega, \cdot)\circ \cdots \circ
{F}(\omega, \cdot), & n >0, \\
{\rm id}, & n=0, \\
{F}^{-1}(\theta^n\omega, \cdot )\circ \cdots \circ
{F}^{-1}(\theta^{-1}\omega, \cdot), & n<0.
\end{cases}
\label{DD-DS}
\end{align}
When $\theta$ is the identity mapping, the RDS becomes a usual deterministic system.

A random variable $K:\Omega \to (0, \infty)$ is said to be  {\it tempered} if it is {\it tempered from above}, i.e.,
$$
\lim_{n\to\pm\infty}\frac1n\max\{0,\log K(\theta^n\omega)\}=0,\qquad
\mathbb{P}-a.s.,
$$
and is {\it tempered from below}, i.e.,
$$
\lim_{n\to\pm\infty}\frac1n\min\{0,\log K(\theta^n\omega)\}=0,\qquad
\mathbb{P}-a.s.
$$
Note that $K(\omega)$ is tempered from above (or below) if and only if $1/K(\omega)$ is tempered from below (or above).
When $K(\omega)$ is tempered above, it follows from \cite{Arn98} that for any constant $\epsilon>0$ there is a tempered random variable $K_\epsilon(\omega)$
such that $K_\epsilon(\omega)\ge K(\omega)$ and
\begin{align}
e^{-\epsilon |n|} K_\epsilon(\omega)
\le K_\epsilon(\theta^n\omega)\le e^{\epsilon |n|} K_\epsilon(\omega), \quad
\forall n\in\mathbb{Z}.
\label{star}
\end{align} Such $K_\epsilon(\omega)$ is called an $\epsilon$-slowly varying random variable.
We call $V(\omega):=\{x\in \mathbb{R}^d:\|x\|\le \rho(\omega)\}$ a {\it tempered ball} if
$\rho(\omega)>0$ is a random variable tempered from below.
For an integer $N\in \mathbb{N}$ and a real number $\alpha\in [0,1]$, we call ${F}(\omega,x)$ a {\it tempered $C^{N,\alpha}$ random diffeomorphism} (see \cite{LiLu2})
if
${F}(\omega,\cdot)$ is a $C^N$ diffeomorphism and the $i$-th derivative $D^i{F}(\omega,x)$ are measurable with respect to $\omega$ such that
\begin{align}
\begin{split}
&\|D^i{F}(\omega,x)\|\le M(\omega),\quad \forall i=0,1,...,N,\quad\forall
\omega\in \Omega,
\\
&\|D^N{F}(\omega,x)-D^N{F}(\omega,y)\|\le L(\omega)\|x-y\|^\alpha,
\quad \forall x,y\in V(\omega),
\end{split}
\label{MMLL}
\end{align}
where $M(\omega),L(\omega)>0$ are random variables tempered from above. In particular, if $\alpha=0$ then it means that $F(\omega, x)$ is $C^N$. Moreover, in the case of $\alpha=0$, we can further assume that $D^NF(\omega, \cdot)$ is $\epsilon$-slowly continuous with respect to $x$ for each fixed $\omega\in \Omega$, where $\epsilon>0$ is a small constant,  i.e., $e^{-\epsilon |n|}D^N F(\theta^n \omega, \cdot)$ is equicontinuous with respect $n\in \mathbb{Z}$ for each fixed $\omega\in\Omega$. When $D^N F(\omega, x)$ is H\"older continuous with a tempered H\"older constant, then $D^N F(\omega, \cdot)$ is $\epsilon$-slowly continuous; when the system is deterministic, the uniform continuity of $D^N F$ implies that it is $\epsilon$-slowly continuous.

We  assume that ${F}(n, \omega, 0)=0$ for any $n\in \mathbb{Z}$ and
any $\omega \in\Omega$, and we use the notation $DF(\omega,\cdot)$ (instead of $D^1F(\omega,\cdot)$) for the first order derivative of $F$ with respect to $x$.

Consider a linear random diffeomorphism $\Lambda(\omega):=\Lambda(\omega,\cdot):\mathbb{R}^d \to \mathbb{R}^d$. Then it generates a linear random dynamical system
\begin{align}
{\Lambda}(n, \omega):=
\begin{cases} \Lambda(\theta^{n-1}\omega) \cdots \Lambda(\omega), & n >0, \\
{\rm id}, & n=0, \\
\Lambda^{-1}(\theta^n\omega) \cdots \Lambda^{-1}(\theta^{-1}\omega),& n<0.
\end{cases}
\label{RS-CYC}
\end{align}
Assume that $\Lambda(\omega)$ satisfies the
conditions of the Multiplicative Ergodic Theorem, that is,
\begin{equation}\label{MET}
\max\{0,\log\|\Lambda(\cdot)\|\} \in L^1(\Omega, \mathcal{F}, \mathbb{P})\quad
\hbox{and} \quad \max\{0,\log\|\Lambda^{-1}(\cdot)\|\}\in L^1(\Omega,
\mathcal{F}, \mathbb{P}).
\end{equation} and $\mathbb{P}$ is an ergodic invariant measure for $\theta$.
Then, by the Multiplicative Ergodic Theorem (\cite{Arn98}), there exists a $\theta$-invariant set
$\tilde \Omega\subset \Omega$ of full measure such that for each
$\omega \in \tilde \Omega$,  the linear random
dynamical system $\Lambda(n, \omega)$ has $p$
Lyapunov exponents: $\lambda_j$ with the multiplicity $d_j$, $j=1,...,p$, and the corresponding Oseledets spaces $E_j(\omega)$ with dimension $d_j$ form an invariant splitting of the phase space $\mathbb{R}^d$:
\[
\mathbb{R}^d= E_1(\omega)\oplus\cdots\oplus E_{{p}}(\omega).
\]
Since $\mathbb{P}$ is ergodic, all $\lambda_i$, $d_i$, and $p$ are constant.
For the remainder of this article, we use $\Omega$ to denote $\tilde \Omega$ and assume that these statements are true for all $\omega\in \Omega$.

Assume that   the linear random
dynamical system $\Lambda(n, \omega)$ is hyperbolic, i.e., there is no zero Lyapunov exponent.
We arrange the Lyapunov exponents as
\begin{align}
\lambda_1>\cdots> \lambda_{\tau}>0>\lambda_{\tau+1}>\cdots>\lambda_{{p}}
\label{LE-H}
\end{align}
for a certain
$\tau\in \{1,...,{p}-1\}$.
It follows from \cite[Lemma 2.8]{LiLu2} that  $\Lambda(n,\omega)$ is conjugate to a block diagonal linear random dynamical system by a random isomorphism and the corresponding Oseledets spaces are given by
\[
E_j(\omega)=X_j:=\{0\} \times \cdots \times \{0\} \times
\mathbb{R}^{d_j} \times \{0\} \times \cdots \times\{0\}
\subset \mathbb{R}^d.
\]
Moreover, for any small $\epsilon>0$ one has
\begin{align}
\begin{split}
&\|\Lambda(n,\omega)|_{X_j}\| \leq K(\omega)
e^{(\lambda_j-\epsilon)n},
\qquad \forall n< 0,
\\
&\|\Lambda(n,\omega)|_{X_j}\| \leq  K(\omega)
e^{(\lambda_j+\epsilon)n},
\qquad \forall n\geq 0,
\end{split}
\label{RS-ED}
\end{align}
for each $j=1,...,p$,
where $K(\omega)>0$ is a random variable tempered from above.
In particular, it satisfies the exponential dichotomy
\begin{align}
\begin{split}
&\|\Lambda(n,\omega)|_{X_u}\| \leq K(\omega)
e^{(\lambda_\tau-\epsilon)n},
\qquad \forall n< 0, \\
&\|\Lambda(n,\omega)|_{X_s}\| \leq  K(\omega)
e^{(\lambda_{\tau+1}+\epsilon)n},
\qquad \forall n\geq 0,
\end{split}
\label{RS-ED-2}
\end{align}
where $\lambda_{\tau}>0>\lambda_{\tau+1}$ and
$
X_u:=X_1\oplus \cdots \oplus X_{\tau},
$
$
X_s:=X_{{\tau}+1}\oplus \cdots \oplus X_{{p}}.
$
Notice that this exponential dichotomy is nonuniform since $ K(\omega)$ may be arbitrarily large as $\omega$ varies. However, along each orbit $(\theta^n \omega)_{n\in\mathbb{Z}}$,  $K(\omega)$  can increase only at a subexponential rate  which is one of the intrinsic features of RDS. The following definition on random conjugacy can be found in \cite{Arn98,LiLu2,LZZ}.
\begin{definition}
{\rm
Two tempered random diffeomporhisms $F(\omega,\cdot)$ and $G(\omega,\cdot)$
are said to be {\it $C^{k,\beta}$ conjugate},
where $k\ge 0$ is an integer and $\beta\in [0,1]$ is a real,
if there exists a $C^{k,\beta}$ random diffeomorphism
$\Phi(\omega,\cdot):V(\omega)\to \mathbb{R}^d$ such that
\begin{align}
\Phi(\theta\omega,F(\omega,x))=G(\omega, \Phi(\omega,x)),\qquad \forall x\in V(\omega),
\label{conjugacy}
\end{align}
where $V(\omega)$ is a tempered ball. In particular, as mentioned in the introduction,
$F$ is said to be {\it $C^{k,\beta}$ linearized}
if $G$ is a linear random diffeomporhism.
}
\end{definition}

Now, we restate our main theorem of this paper for the sake of completeness of this section. 

\begin{thm}\label{Main-ThmC1-RS}
Let $F:\Omega \times V(\omega)\to \mathbb{R}^d$ be a tempered $C^{2,\alpha}$ $(\alpha\in[0,1])$ random diffeomorphism with an $\epsilon$-slowly continuous second order derivative when $\alpha=0$. Assume that its linearization $\Lambda(\omega)=DF(\omega, 0)$ at the hyperbolic fixed point $x=0$ satisfies the condition \eqref{MET}. If  the Lyapunov exponents satisfy
a partial second order non-resonant condition
\begin{align}
\lambda_j\ne \lambda_i+\lambda_\kappa
\label{EC-RS}
\end{align}
for $i=1,...,{\tau}$, $\kappa={\tau}+1,...,{{p}}$ and $j=1,...,{{p}}$,
then the following assertions are true:
\begin{itemize}
\item[(i)] Foliations:  $F$ has $C^{1,\beta_\alpha}$ random stable and unstable foliations with $C^{2,\alpha}$ leaves in a tempered ball $U(\omega)\subset V(\omega)$, where $\beta_\alpha=0$ when $\alpha=0$ and $\beta_\alpha\in (0,\alpha]$ when $\alpha\in (0,1]$.
\item[(ii)] Conjugacy:
$F(\omega,\cdot)$ is conjugate
to its linear part $\Lambda(\omega)$ in a tempered ball $U(\omega)\subset V(\omega)$ by a $C^{1,\beta_\alpha}$  random diffeomorphism $\Phi:\Omega\times U(\omega)\to \mathbb{R}^d$ satisfying
$
\Phi(\omega,0)=0
$
and
$
D\Phi(\omega,0)={\rm id}.
$
\end{itemize}
\end{thm}

Remark that the assumption of $\epsilon$-slowly continuity is only necessary for random dynamical systems in the case of $\alpha=0$ since, when $\alpha\in (0,1]$ or $F$ is deterministic, the second order derivative of $F(\omega,\cdot)$ is automatically
$\epsilon$-slowly continuous.
Moreover, when $F$ is deterministic,
the relationship between  Lyapunov exponent $\lambda$ and its corresponding eigenvalue
$r$ is $|r|=e^\lambda$. Thus condition
\eqref{EC-RS} becomes (\ref{B-cond}). As a result, we have the following $C^{1}$ linearization theorem.

\begin{cor}\label{Main-CorC1}
Let $F:V\to \mathbb{R}^d$ be a $C^{2,\alpha}$ {\rm (}$\alpha\in [0,1]${\rm )} diffeomorphism over a neighborhood $V\subset \mathbb{R}^d$ of $x=0$. Suppose the eigenvalues of its linearization $\Lambda:=DF(0)$ at hyperbolic fixed point $x=0$ satisfy  {\rm (\ref{B-cond})}.
Then $F$ is conjugate
to its linear part $\Lambda$ in a small neighborhood $U\subset V$ by a $C^{1,\beta_\alpha}$ diffeomorphism $\Phi: U\to \mathbb{R}^d$ satisfying
$
\Phi(0)=0
$
and
$
D\Phi(0)={\rm id},
$
where $\beta_\alpha=0$ when $\alpha=0$ and $\beta_\alpha\in (0,\alpha]$ when $\alpha\in (0,1]$.
\end{cor}

This corollary extends Belitskii's $C^1$ linearization theorem to $C^{1,\beta_\alpha}$ linearization, where $\beta_\alpha>0$ whenever $\alpha>0$. Remark that $C^{1,\beta_\alpha}$ linearization with $\beta_\alpha>0$ has important applications in problems such as Lorenz attractors (\cite{HollMelb-07}) and topological entropy (\cite{HomWe-CMP99}).

\section{Random invariant distribution}
\setcounter{equation}{0}

In this section, we study the H\"older continuity of random invariant distributions. We first extend the local random diffeomorphism $F(\omega,x)=\Lambda(\omega)x+f(\omega,x)$ defined on a tempered ball $V(\omega)=\{x\in \mathbb{R}^d:\|x\|\le \rho(\omega)\}$ to a global one. In what follows, for a random variable (for example $K:\Omega \to (0, \infty)$), we use $K(\omega)^n$ to denote $K(\omega)$ to the power of $n\in\mathbb{Z}$ and let $U(\omega):=\{x\in \mathbb{R}^d:\|x\|\le \rho(\omega)/2\}$,
where
$$
\|x\|:=\max\{\|x_1\|,...,\|x_p\|\}\quad {\rm for} ~~ x=x_1+\cdots+x_p\in X_1\oplus\cdots \oplus X_p.
$$

Consider a random smooth cut-off function
 $u(\omega,\cdot):\mathbb{R}^d\to (0,1)$ such that
\begin{align}
u(\omega,x)=
\left\{
\begin{array}{lll}
1,& \forall x\in U(\omega),
\vspace{1ex}\\
0,& \forall x\in \mathbb{R}^d\backslash V(\omega),
\end{array}
\right.
\quad \rho(\omega)^r\|D^ru(\omega,x)\|\le C_u,~\forall r=1,2,3,
\label{cutoff-w}
\end{align}
where $C_u$ is a positive constant.
The construction of $u$ is given in  Appendix.

Let $\tilde F(\omega,x)$ be an extension of $F(\omega, x)$ given by
\begin{align*}
\tilde F(\omega,x)=\begin{cases}\Lambda(\omega)x+u(\omega,x)f(\omega,x),&\quad \text{for all}\; x\in V(\omega)\\
\Lambda(\omega)x, &\quad \text{for all}\; x\in \mathbb{R}^d\backslash V(\omega).
\end{cases}
\end{align*}
Then, by elementary estimations, we have the following result for the extension $\tilde F(\omega,x)$.

\begin{lm}\label{Extension}
The following properties hold for $\tilde F(\omega,x)$:
\begin{itemize}
  \item [(i)] $\tilde F(\omega,x)=F(\omega, x)$ for all $x\in U(\omega)$.
  \item [(ii)] The derivatives of $\tilde F$ satisfy the following estimates.
  \begin{align*}
\begin{split}
&\|D\tilde F(\omega,x)-\Lambda(\omega)\|\le \delta(\omega),\qquad \|D^2\tilde F(\omega,x)\|\le M_{\epsilon}(\omega),
\\
&\|D^2\tilde F(\omega,x)-D^2\tilde F(\omega,y)\|\le L_{\epsilon,\rho}(\omega)\|x-y\|^\alpha, \quad \forall x,y\in \mathbb{R}^d,
\end{split}
\end{align*}
\end{itemize}
where $\delta(\omega):=M_\epsilon(\omega)\rho(\omega)$, $M_\epsilon(\omega)$ and $L_{\epsilon,\rho}(\omega)$ are $\epsilon$-slowly varying random variables for any $\epsilon>0$ such that $M_\epsilon(\omega)\ge (3C_u+1)M(\omega)$ and  $L_{\epsilon,\rho}(\omega)\ge 7C_uM(\omega)/\rho(\omega)+L(\omega)$. Here  $M(\omega)$ and $L(\omega)$ are given in \eqref{MMLL} with $N=2$.
\end{lm}
In what follows, we still let $F(\omega,\cdot)$ denote the global mapping $\tilde F(\omega,\cdot)$, and therefore
$F(\omega,\cdot):\mathbb{R}^d\to \mathbb{R}^d$ satisfies that
\begin{align}
\begin{split}
&\|DF(\omega,x)-\Lambda(\omega)\|\le \delta(\omega),\qquad \|D^2F(\omega,x)\|\le M_{\epsilon}(\omega),
\\
&\|D^2F(\omega,x)-D^2F(\omega,y)\|\le L_{\epsilon,\rho}(\omega)\|x-y\|^\alpha, \quad \forall x,y\in \mathbb{R}^d,
\label{c-off-RS}
\end{split}
\end{align}
and $D^2F(\omega,\cdot)$ is $\epsilon$-slowly continuous when $\alpha=0$.
Note that the third inequality of (\ref{c-off-RS}) holds for every H\"older exponents $\tilde \alpha\in [0,\alpha]$ because $F(\omega,x)-\Lambda(\omega)x$ is equal to $0$ outside the small neighborhood $V(\omega)$.

Throughout this paper, we choose  $\epsilon>0$ being sufficiently small such that
\begin{align}
\epsilon
<\frac{1}{100}
\min\Big\{
&1,~
\min_{i=1,...,p-1}\frac{\lambda_{j}-\lambda_{j+1}}
{2\lambda_{\max}}\min\{\lambda_{\tau},
-\lambda_{\tau+1}\},
\min_{\mbox{\tiny
$\begin{array}{c}
i\in\{1,...,\tau\},\kappa\in\{\tau+1,...,{p}\},
\\
j\in\{1,...,{p}\}
\end{array}$}}
|\lambda_i+\lambda_\kappa-\lambda_j|\Big\},
\label{def-ep}
\end{align}
where $\lambda_{\rm max}:=\max\{2\lambda_1,-2\lambda_{{p}}\}>0$. Note that the choice of $\epsilon$ depends on the Lyapunov exponents only.
Let the diameter of $V(\omega)$ be defined by
$$
\rho(\omega):= 1/(M_\epsilon(\omega)
\mathfrak{M}_\epsilon(\omega))>0
$$
for a given $\epsilon$-slowly varying tempered random variable $\mathfrak{M}_\epsilon(\omega)>0$ satisfying
\begin{align}
\begin{split}
\mathfrak{M}_\epsilon(\omega)\ge
\max\bigg\{
&4,~\delta_\lambda^{-1}K_\epsilon(\theta\omega),~
8 C_\lambda K_\epsilon(\omega)^3,~(2 M_\epsilon(\omega) K_\epsilon(\omega)^5)^{2},~(2 K_\epsilon(\omega)^2)^{4},
\\
&\big\{12K_\epsilon(\omega)^2 e^{\lambda_{j}-\lambda_{j+1}-3\epsilon}\big\}
^{4(2\lambda_{\max}-\lambda_{j+1}-\epsilon)
/(\lambda_{j}-\lambda_{j+1}-3\epsilon)},
\\
&2C_\lambda K_\epsilon(\omega)C_\epsilon(\omega), 2 e^{\lambda_{\max}}C^3_\lambda
K_\epsilon(\omega) K_\epsilon (\theta\omega)
\bigg\}
\end{split}
\label{MM-p}
\end{align}
where $K_\epsilon(\omega)\ge \max\{1, K_\lambda K(\omega)\}>0$ and $C_\epsilon(\omega):=4(K_\epsilon(\omega))^2
M_\epsilon(\omega)>0$ are $\epsilon$-slowly varying random variables, $\delta_\lambda>0$ is a small constant and $C_\lambda, K_\lambda\ge 1$ are constants. Note that the constants $\delta_\lambda, C_\lambda, K_\lambda$ are determined by the Lyapunov exponents only, which will be given in Lemma \ref{center-growth} and \eqref{project-difference}. Then, as seen in Lemma \ref{Extension},
\begin{align}\label{def-d}
\delta(\omega):=M_\epsilon(\omega)\rho(\omega)
=1/\mathfrak{M}_\epsilon(\omega)>0
\end{align}
is small and tempered.

Next, we obtain from (\ref{DD-DS}) that
\begin{align}
D{F}(n, \omega, x)=
\begin{cases}
D{F}(\theta^{n-1}\omega,{F}(n-1, \omega, x))\circ \cdots \circ
D{F}(\omega, {F}(0, \omega, x)), & n >0, \\
{\rm id}, & n=0, \\
\big\{D{F}(\theta^n\omega, D{F}(n, \omega, x))\big\}^{-1}\circ \cdots \circ
\big\{D{F}(\theta^{-1}\omega, D{F}(-1, \omega, x))\big\}^{-1}, & n<0.
\end{cases}
\label{DDFFnn}
\end{align}
For every $\omega\in \Omega$ and $x\in \mathbb{R}^d$, we set
\[
\varpi:=(\omega, x)\in {\bf \Omega}:=\Omega\times \mathbb{R}^d,\qquad \vartheta\varpi:=(\theta\omega,F(\omega,x)),\qquad {\bf \Lambda}(\varpi):=DF(\omega,  x).
\]
Then, $D{F}(n, \omega, x)$ is a cocycle generated by  ${\bf \Lambda}(\varpi)$ driven by   $\vartheta^n \varpi:=(\theta^n\omega,F(n, \omega,x))$. The first inequality of (\ref{c-off-RS}) implies that ${\bf \Lambda}(\varpi)$ is a small perturbation of $\Lambda(\omega)$ with the difference
\[
\|{\bf \Lambda}(\varpi)-\Lambda(\omega)\|\le \delta(\omega).
\] Note that $\Lambda(\omega)=\Lambda(\varpi)$ since $\Lambda(\omega)$ does not depend on $x$.
Using (\ref{RS-ED}) and  applying the roughness of the tempered exponential dichotomy (\cite[Theorem 1]{ZLZ}) to the system $e^{-(\lambda_j+\lambda_{j+1}) n/2} \Lambda(n,\omega)$ for each $j=1,\cdots, p-1$ with the splitting
\[
\mathbb{R}^d=\big(X_1\oplus \cdots \oplus X_{j}\big)\oplus \big(X_{j+1}\oplus \cdots \oplus X_{p}\big),
\]
we have the following lemma.
\begin{lm}\label{center-growth} There are positive constants $\delta_\lambda$ and $K_\lambda$ depending on the Lyapunov exponents only such that if $\delta(\omega)\leq  \delta_\lambda/K_\epsilon(\theta\omega)$, then there is an invariant splitting for  the cocycle ${\bf \Lambda}(n, \varpi)$ generated by  ${\bf \Lambda}(\varpi)$
\[
\mathbb{R}^d=E_1(\omega,x)\oplus \cdots\oplus E_{{p}}(\omega,x)
\] such that for each $j=1,\cdots, p$,
\begin{align*}
\begin{split}
&\|{\bf \Lambda}(n, \varpi)|_{E_j}\| \leq K_\lambda K(\omega)
e^{(\lambda_j-2\epsilon)n},
\qquad \forall n< 0,
\\
&\|{\bf \Lambda}(n, \varpi)|_{E_j}\| \leq  K_\lambda K(\omega)
e^{(\lambda_j+2\epsilon)n},
\qquad \forall n\geq 0.
\end{split}
\end{align*}
\end{lm}
The invariance means that ${\bf \Lambda}(n, \varpi) E_j( \varpi)=E_j(\vartheta^n \varpi)$.
Without loss of generality, we still let $\epsilon$ denote $2\epsilon$. Note that
\[D{F}(n, \omega, x)={\bf \Lambda}(n, \varpi)\; \text{ and } F(\omega, x)=F(1, \omega, x).
\]
Thus, we have
 \begin{align}
DF(\omega,x)\{E_j (\omega,x)\}
=E_j(\theta\omega,F(\omega,x)),\qquad \forall j=1,...,{p},
\label{Int-m-Dis}
\end{align}
and
\begin{align}
\begin{split}
&\|DF(n,\omega,x)|_{E_j(\omega,x)}\| \leq K_\epsilon(\omega)
e^{(\lambda_j-\epsilon)n},
\qquad \forall n< 0,
\\
&\|DF(n,\omega,x)|_{E_j(\omega,x)}\| \leq K_\epsilon(\omega)
e^{(\lambda_j+\epsilon)n},
\qquad \forall n\geq 0,
\end{split}
\label{int-m-d}
\end{align}
where $K_\epsilon(\omega)\ge K_\lambda K(\omega)$ is an $\epsilon$-slowly varying random variable.
 We point out that the measurability of $E_j(\omega,x)$ with respect to $\omega$ follows from that $E_j(\omega,x)$ is the tangent spaces of an intermediate foliation whose partial derivative is measurable with respect to $\omega$ (see (\ref{EWjj}) below in the next section).  As in the deterministic case (see \cite[Sections 2 and 6]{Pesin-book}), for each $j=1,...,p$,
 \[
E_j(\omega):=\{E_j(\omega, x)\subset \mathbb{R}^d:x\in \mathbb{R}^d\}
\] is called
a {\it random intermediate distribution} and  $x$ is  the {\it base point} of the distribution. In particular,
$$
E_u(\omega):=\{E_u(\omega, x)\subset \mathbb{R}^d:x\in \mathbb{R}^d\}\mbox{ and }
E_s(\omega):=\{E_s(\omega, x)\subset \mathbb{R}^d:x\in \mathbb{R}^d\}
$$
are called
the {\it random unstable distribution and random stable distribution} respectively, with fibers
\begin{align*}
&E_u(\omega,x):=E_1(\omega,x)\oplus \cdots\oplus E_{\tau}(\omega,x),
\\
&E_s(\omega,x):=E_{\tau+1}(\omega,x)\oplus \cdots\oplus E_{{p}}(\omega,x)
\end{align*}
satisfying the following:
\begin{itemize}
\item[(i)]
$\mathbb{R}^d=E_u(\omega, x) \oplus E_s(\omega, x);$
\item[(ii)] The splitting is invariant, i.e.,
 \begin{align}
\begin{split}
&DF(\omega,x)E_u (\omega,x)
=E_u(\theta\omega,F(\omega,x)),
\\
&DF(\omega,x)E_s (\omega,x)
=E_s(\theta\omega,F(\omega,x));
\end{split}
\label{Inv-Dis}
\end{align}
\item[(ii)] For each $\varrho_0\in (0,\,\min\{\lambda_{\tau}, -\lambda_{\tau+1}\})$, there exists a random variable $ K_\epsilon(\omega)>0$ tempered from above  such that
\begin{align}
\begin{split}
\|DF(n,\omega,x)|_{E_u (\omega,x)}\|& \leq K_\epsilon(\omega)
e^{\varrho_0  n},
\qquad \forall n< 0, \\
\|DF(n,\omega,x)|_{E_s (\omega,x)}\|& \leq  K_\epsilon(\omega)
e^{-\varrho_0  n},
\qquad \forall n\geq 0.
\end{split}
\label{ED}
\end{align}
\end{itemize}

Finally, we show that each $E_j(\omega,x)$ is H\"older continuous with respect to $x$. Recall that the distance between two subspaces $E$ and $\tilde E$ of $\mathbb{R}^d$ is given by
\begin{align}\label{dis-subsp}
{\rm dist}(E,\tilde E):=\max\Big(\max_{v_1\in E,\,\|v_1\|=1}{\rm dist}(v_1, \tilde E),\max_{v_2\in \tilde E,\,\|v_2\|=1}{\rm dist}(v_2, E)\Big),
\end{align}
where ${\rm dist}(v,E):=\min_{w\in E}\|v-w\|$.

\begin{thm} \label{thm-HD}
For all $x,y\in \mathbb{R}^d$, $\omega\in \Omega$ and all $j=1,...,{p}$,
\begin{align}
\begin{split}
&{\rm dist}(E_j(\omega, x),E_j(\omega, y))\le \|x-y\|^{\beta_E},
\\
&{\rm dist}(E_j(\omega, x),E_j(\omega, y))\le \delta_E(\omega),
\end{split}
\label{DB-us}
\end{align}
where
$
\beta_{E}:=\min_{\ell=1,...,p-1}(\lambda_{\ell}-\lambda_{\ell+1}-3\epsilon)
/(6\lambda_{\max})\in (0,1)$ and
$\delta_E(\omega):=\delta(\omega)^{\beta_{E}}.
$
\end{thm}

{\bf Proof}. First of all, we define the subspaces
\begin{align*}
&E_+(\omega, x):=E_1(\omega, x) \oplus \cdots \oplus E_j(\omega, x),
\\
&E_-(\omega, x):=E_{j+1}(\omega, x)\oplus \cdots \oplus E_{{p}}(\omega, x)
\end{align*}
and claim that for all integers $m\ge 0$
\begin{align}
\begin{split}
\|DF(m,\omega,x)\zeta_-\|& \leq 2K_\epsilon(\omega)
e^{(\lambda_{j+1}+\epsilon)m}\|\zeta_-\|,
\\
\|DF(m,\omega,x)\zeta_-^\bot\|& \ge (2K_\epsilon(\omega))^{-1} e^{(\lambda_{j}-2\epsilon)m}\|\zeta_-^\bot\|
\end{split}
\label{ED-md}
\end{align}
for $\zeta_-\in E_-(\omega,x)$ and $\zeta_-^\bot\in \{E_-(\omega,x)\}^\bot$, the orthogonal complement of $E_-(\omega,x)$. In fact, we have from (\ref{int-m-d}) that
\begin{align}
\begin{split}
&\|DF(n,\omega,x)\zeta_+\| \leq K_\epsilon(\omega)
e^{(\lambda_j-\epsilon)n}\|\zeta_+\|,
\qquad \forall n< 0, \\
&\|DF(n,\omega,x)\zeta_-\| \leq  K_\epsilon(\omega)
e^{(\lambda_{j+1}+\epsilon)n}\|\zeta_-\|,
\qquad \forall n\geq 0,
\end{split}
\label{int-m-+-}
\end{align}
for $\zeta_+\in E_+(\omega, x)$ and $\zeta_-\in E_-(\omega, x)$. Then the first inequality of (\ref{ED-md}) follows from
the second one of (\ref{int-m-+-}) immediately.

Now, we show that the second inequality of (\ref{ED-md}) holds. We first notice that  for $m\ge 0$
\begin{align*}
\big(DF(-m,\omega,x)\big)^{-1}=DF(m,\theta^{-m}\omega,F(-m,\omega,x))
\end{align*}
and therefore the first inequality of (\ref{int-m-+-}) yields
\begin{align*}
\|DF(m,\theta^{-m}\omega,F(-m,\omega,x))\tilde \zeta_+\|
\ge (1/K_\epsilon(\omega))e^{(\lambda_j-\epsilon) m}\|\tilde \zeta_+\|
\end{align*}
for $\tilde\zeta_+ \in  E_+(\theta^{-m}\omega, F(-m, \omega, x))$.
Replacing $F(-m,\omega,x)$ and $\theta^{-m}\omega$ with $x$ and  $\omega$ respectively,
we obtain
\begin{align}
\|DF(m,\omega,x)\zeta_+\|
\ge (1/K_\epsilon(\theta^{m}\omega))e^{(\lambda_j-\epsilon) m}\|\zeta_+\|
\ge (1/K_\epsilon(\omega))e^{(\lambda_j-2\epsilon) m}\|\zeta_+\|.
\label{ED-ge}
\end{align}
Let $\Pi_\pm(\omega, x):\mathbb{R}^d\to E_\pm(\omega, x)$ be the projections associated with $\mathbb{R}^d=E_+(\omega, x)\oplus E_-(\omega, x)$ and let
$\pi_\pm:\mathbb{R}^d\to X_\pm$ be the projections associated with $\mathbb{R}^d=X_+\oplus X_-$, where
\begin{align*}
X_+:=X_1 \oplus \cdots \oplus X_{j},\quad
X_-:=X_{j+1}\oplus \cdots \oplus X_{{p}}.
\end{align*}
Then, using the arguments as in \cite[(3.41)-(3.44) and (3.16)]{ZLZ},
we have
\begin{align} \label{project-difference}
\|\Pi_+(\omega, x)-\pi_+\|\le C_\lambda\delta(\omega)K_\epsilon(\omega),
\quad
\|\Pi_-(\omega, x)-\pi_-\|\le C_\lambda\delta(\omega)K_\epsilon(\omega),
\end{align}
where $C_\lambda>0$ is a constant depending only on the Lyapunov exponents.

Using the second inequality of \eqref{project-difference}, when $\delta(\omega)$ is small enough, there are
$p-j$ vectors
\begin{align*}
&v_1^-:=(v_{1,1},\cdots,v_{1,j}, 1 , \cdots , 0)^T,~...,~ v_{p-j}^-:=(v_{p-j,1},\cdots,v_{p-j,j}, 0 , \cdots , 1)^T\in E_-(\omega,x)
\end{align*}
with $v_{i,j}\in\mathbb{R}$ ($T$ means the transpose of a vector) such that
\begin{align*}
\max\{|v_{i,1}|,\cdots,|v_{i,j}|\}=\|v_{i}^--\pi_-v_{i}^-\|
=\|(\Pi_-(\omega, x)-\pi_-)v_{i}^-\|
\le C_\lambda\delta(\omega)K_\epsilon(\omega)\|v_{i}^-\|,
\end{align*}
which implies that  $\|v_{i}^-\|=\max\{|v_{i,1}|,\cdots,|v_{i,j}|,1\}=1$
since $\delta(\omega)>0$ is sufficiently small.
Therefore
\begin{align}
\max\{|v_{i,1}|,\cdots,|v_{i,j}|\}=\|v_{i}^--\pi_-v_{i}^-\|
\le C_\lambda\delta(\omega)K_\epsilon(\omega).
\label{small-v}
\end{align}
For any  vector
$
\zeta_-^\bot:=(\zeta_1,...\zeta_j, \zeta_{j+1},..., \zeta_p)\in (E_-(\omega, x))^{\bot},
$
since $\zeta_-^\bot$ is perpendicular to $v_1^-, \cdots, v_{p-j}^-$, we have
$
B\zeta_-^\bot=0,
$
where
\begin{align*}
B:=\begin{pmatrix}
    v_{1,1} &\cdots &v_{1,j}& 1  & \cdots & 0 \\
    \vdots  & \ddots & \vdots & \vdots & \ddots & \vdots \\
    v_{p-j,1} &\cdots &v_{p-j,j}& 0 &\cdots &1
  \end{pmatrix}.
\end{align*}
Thus, $B\zeta_-^\bot=0$ is equivalent to
\begin{align*}
\begin{pmatrix}
    \zeta_{j+1}\\
    \vdots\\
    \zeta_p
  \end{pmatrix}
  =
  \begin{pmatrix}
    -v_{1,1} &\cdots &-v_{1,j} \\
    \vdots & \ddots & \vdots \\
    -v_{p-j,1} &\cdots &-v_{p-j,j}
  \end{pmatrix}
  \begin{pmatrix}
    \zeta_{1}\\
    \vdots\\
    \zeta_j
  \end{pmatrix},
\end{align*}
which together with (\ref{small-v}) indicates that
\begin{align*}
\|\pi_- \zeta_-^\bot\|
\le j C_\lambda\delta(\omega)K_\epsilon(\omega)\|\zeta_-^\bot\|
\le C_\lambda\delta(\omega)K_\epsilon(\omega)\|\zeta_-^\bot\|.
\end{align*}
Here, we omit $j\le p$ (the number of the Lyapunov exponents) since $C_\lambda>0$ denotes a sufficiently large constant depending on the Lyapunov exponents. It
implies that
\[
\|\pi_+ \zeta_-^\bot\|\ge \|\zeta_-^\bot\|-\|\pi_- \zeta_-^\bot\|\ge (1- C_\lambda \delta(\omega)K_\epsilon(\omega))\|\zeta_-^\bot\|.
\]
Therefore,
\begin{align}
\|\Pi_+(\omega, x)\zeta_-^\bot\|
&\ge \|\pi_+\zeta_-^\bot\|-\|\Pi_+(\omega, x)\zeta_-^\bot-\pi_+\zeta_-^\bot\|
\nonumber\\
&\ge \{1-2C_\lambda\delta(\omega)K_\epsilon(\omega)\}\|\zeta_-^\bot\|
\label{eu-1}
\end{align}
and
\begin{align}
\|\Pi_-(\omega, x)\zeta_-^\bot\|
&\le
\|\pi_-\zeta_-^\bot\|+\|\Pi_-(\omega, x)\zeta_-^\bot-\pi_-\zeta_-^\bot\|
\nonumber\\
&\le 2C_\lambda\delta(\omega)K_\epsilon(\omega)\|\zeta_-^\bot\|.
\label{es-2}
\end{align}
It follows from (\ref{int-m-+-})-(\ref{es-2}) that
\begin{align}
&\|DF(m,\omega,x)\zeta_-^\bot\|
\nonumber\\
&\ge
\|DF(m,\omega,x)\Pi_+(\omega, x)\zeta_-^\bot\|-\|DF(m,\omega,x)\Pi_-(\omega, x)\zeta_-^\bot\|
\nonumber\\
&\ge (1/K_\epsilon(\omega))e^{(\lambda_j-2\epsilon) m}\|\Pi_+(\omega, x)\zeta_-^\bot\|-K_\epsilon(\omega)
e^{(\lambda_{j+1}+\epsilon) m}\|\Pi_-(\omega, x)\zeta_-^\bot\|
\nonumber\\
&\ge (1/K_\epsilon(\omega))
\{1-2C_\lambda\delta(\omega)K_\epsilon(\omega)\}e^{(\lambda_j-2\epsilon) m}\|\zeta_-^\bot\|
\nonumber\\
&\quad-
2C_\lambda\delta(\omega)K_\epsilon(\omega)^2e^{(\lambda_{j+1}+\epsilon) m}\|\zeta_-^\bot\|
\nonumber\\
&\ge (1/K_\epsilon(\omega))e^{(\lambda_{j}-2\epsilon) m}\|\zeta_-^\bot\|-4C_\lambda
\delta(\omega)K_\epsilon(\omega)^2e^{(\lambda_{j}-2\epsilon) m}\|\zeta_-^\bot\|
\nonumber\\
&\ge (2K_\epsilon(\omega))^{-1}e^{(\lambda_{j}-2\epsilon) m}\|\zeta_-^\bot\|
\label{unq-s}
\end{align}
since
$
\mathfrak{M}_\epsilon(\omega)\ge 8 C_\lambda K_\epsilon(\omega)^3
$
(i.e., $0<\delta(\omega)\le (8 C_\lambda K_\epsilon(\omega)^3)^{-1}$) by (\ref{MM-p}). This shows the second inequality of (\ref{ED-md}) and (\ref{ED-md}) is proved.

Next, we know from the first inequality of (\ref{c-off-RS}) that
\begin{align*}
\|DF(\omega,x)-DF(\omega,y)\|
\le \|DF(\omega,x)-\Lambda(\omega)\|+\|DF(\omega,y)-\Lambda(\omega)\|
\le 2\delta(\omega),\quad \forall x,y\in \mathbb{R}^d,
\end{align*}
and from (\ref{int-m-d}) that
\begin{align}
\|DF(m,\omega,x)\|\le \max_{j=1,...,p}\|DF(m,\omega,x)|_{E_j(\omega,x)}\|\le K_\epsilon(\omega)e^{\lambda_{\max}m},
\qquad \forall m\ge 0,
\label{DF-max}
\end{align}
where
$
\lambda_{\rm max}:=\max\{2\lambda_1,-2\lambda_{{p}}\}>0.
$
Then, we have
\begin{align}
&\|DF(\theta^m\omega,F(m,\omega,x))-
DF(\theta^m\omega,F(m,\omega,y))\|
\nonumber\\
&=\|DF(\theta^m\omega,F(m,\omega,x))-
DF(\theta^m\omega,F(m,\omega,y))\|^{1/2+1/2}
\nonumber\\
&\le (2\delta(\theta^m\omega))^{1/2}\{M_\epsilon(\theta^m\omega)
\|F(m,\omega,x)-F(m,\omega,x)\|\}^{1/2}
\nonumber\\
&\le (2\delta(\theta^m\omega))^{1/2}\{M_\epsilon(\theta^m\omega)
K_\epsilon(\omega)e^{\lambda_{\max}m}\|x-y\|\}^{1/2}
\nonumber\\
&\le \{2\delta(\omega)M_\epsilon(\omega)
K_\epsilon(\omega)\}^{1/2}
e^{(1/2)(\lambda_{\max}+2\epsilon) m}\|x-y\|^{1/2}
\label{DF-DF-1}
\end{align}
by the second inequality of (\ref{c-off-RS}) and
\begin{align}
\|DF(\theta^m\omega,F(m,\omega,x))-
DF(\theta^m\omega,F(m,\omega,y))\|
\le 2\delta(\theta^m\omega)\le 2e^{\epsilon m}\delta(\omega).
\label{DF-DF-2}
\end{align}
In the above (\ref{DF-DF-1})-(\ref{DF-DF-2}), we know that
$
e^{-\epsilon |k|}\delta(\theta^k\omega)\le\delta(\omega),
$
as known from (\ref{def-d}).

Notice the fact that
\begin{align*}
a_{m-1}...a_0-b_{m-1}...b_0
&=a_{m-1}...a_1\,(a_0-b_0)
+(a_{m-1}-b_{m-1})\,b_{m-2}...b_0
\\
&=\sum_{i=0}^{m-1}a_{m-1}...a_{i+1}\,(a_i-b_i)\,b_{i-1}...b_1.
\end{align*}
It follows from (\ref{DF-DF-1}) that
\begin{align}
&\|DF(m,\omega,x)-DF(m,\omega,y)\|
\nonumber\\
&=\|D{F}(\theta^{m-1}\omega,{F}(m-1, \omega, x))\circ \cdots \circ D{F}(\omega, {F}(0, \omega, x))
\nonumber\\
&\quad- D{F}(\theta^{m-1}\omega,{F}(m-1, \omega, y))\circ \cdots \circ D{F}(\omega, {F}(0, \omega, y))\|
\nonumber\\
&\le \sum_{i=0}^{m-1}\|DF(m-i-1,\theta^{i+1}\omega,F(i+1,\omega,x))\|
\nonumber\\
&\quad\cdot\|DF(\theta^i\omega,F(i,\omega,x))-
DF(\theta^i\omega,F(i,\omega,y))\|\,\|DF(i,\omega,y)\|
\nonumber\\
&\le \sum_{i=0}^{m-1}K_\epsilon(\theta^{i+1}\omega)e^{\lambda_{\max}(m-i-1)}
\{2\delta(\omega)M_\epsilon(\omega)
K_\epsilon(\omega)\}^{1/2}
e^{(1/2)(\lambda_{\max}+2\epsilon) i}
\|x-y\|^{1/2}K_\epsilon(\omega)e^{\lambda_{\max}i}
\nonumber\\
&\le \sum_{i=0}^{m-1}
e^{\epsilon(i+1)}e^{\lambda_{\max}(m-1)}
e^{(1/2)(\lambda_{\max}+2\epsilon) i}\{2\delta(\omega)M_\epsilon(\omega)
K_\epsilon(\omega)^5\}^{1/2}
\|x-y\|^{1/2}
\nonumber\\
&\le e^{2\lambda_{\max}m}\delta(\omega)^{1/4}\|x-y\|^{1/2}
\label{DFmm-2}
\end{align}
since
$
\mathfrak{M}_\epsilon(\omega)\ge (2 M_\epsilon(\omega) K_\epsilon(\omega)^5)^{2}
$
(i.e., $\delta(\omega)\le (2 M_\epsilon(\omega)K_\epsilon(\omega)^5)^{-2}$) by (\ref{MM-p}),
where
\[
me^{\epsilon(i+1)}
e^{(1/2)(\lambda_{\max}+2\epsilon) i}\le e^{\lambda_{\max}m}
\]
as $0\le i\le m-1$ and $\epsilon>0$ is small. Similarly, it follows from (\ref{DF-DF-2}) that
\begin{align}
&\|DF(m,\omega,x)-DF(m,\omega,y)\|
\nonumber\\
&\le \sum_{i=0}^{m-1}\|DF(m-i-1,\theta^{i+1}\omega,F(i+1,\omega,x))\|
\nonumber\\
&\quad\cdot\|DF(\theta^i\omega,F(i,\omega,x))-
DF(\theta^i\omega,F(i,\omega,y))\|\,\|DF(i,\omega,y)\|
\nonumber\\
&\le \sum_{i=0}^{m-1}K_\epsilon(\theta^{i+1}\omega)e^{\lambda_{\max}(m-i-1)}
2e^{\epsilon m}\delta(\omega)K_\epsilon(\omega)e^{\lambda_{\max}i}
\nonumber\\
&\le \sum_{i=0}^{m-1}e^{\epsilon(i+1)}e^{\lambda_{\max}(m-1)}
2e^{\epsilon m}\delta(\omega)K_\epsilon(\omega)^2
\le e^{2\lambda_{\max}m}\delta(\omega)^{3/4}
\label{DFmm-1}
\end{align}
since $\mathfrak{M}_\epsilon(\omega)\ge (2 K_\epsilon(\omega)^2)^{4}$ (i.e., $\delta(\omega)\le (2 K_\epsilon(\omega)^2)^{-4}$) by (\ref{MM-p}).

Hence, using (\ref{ED-md}) and (\ref{DFmm-2})-(\ref{DFmm-1}),
we apply \cite[Lemma 3.10]{Pesin-book} for every fixed $\omega\in\Omega$, which is given exactly as follows (in notations used in \cite{Pesin-book}).
\begin{lm}\label{Pesin-lm}
Let $(A_n)_{n\in\mathbb{N}}$ and $(B_n)_{n\in\mathbb{N}}$ be two sequences of real $N\times N$ matrices such that for some $\Delta\in(0,1)$ and $a>1$,
$$
\|A_n-B_n\|\le \Delta a^n
$$
for all positive integers $n$. Assume that there exist subspaces $E_A, E_B\in \mathbb{R}^N$
and numbers $0<\lambda<\mu$ and $C>1$ such that $\lambda<a$ and for each $n\ge 0$,
\begin{align*}
&\|A_n v\|\le C\lambda^n ~ if ~ v\in E_A;\qquad \|A_nw\|\ge C^{-1}\mu^n\|w\|~ if ~ w\in {E_A}^\bot;
\\
&\|B_n v\|\le C\lambda^n ~ if ~ v\in E_B;\qquad \|B_nw\|\ge C^{-1}\mu^n\|w\|~ if ~ w\in {E_B}^\bot.
\end{align*}
Then
$$
{\rm dist}(E_A,E_B)\le 3C^2\frac{\mu}{\lambda}\Delta^{\frac{\log\mu-\log\lambda}
{\log a-\log\lambda}}.
$$
\end{lm}
Notice that if $\|x-y\|<1$ then the notations of \cite[Lemma 3.10]{Pesin-book} and ours have the correspondence:
\begin{align*}
&n:=m, \quad A_m:=DF(m,\omega,x), \quad B_m:=DF(m,\omega,y),
\\
&\Delta:=\delta(\omega)^{1/4}\|x-y\|^{1/2} ~{\rm or}~ \delta(\omega)^{3/4},\quad a:=e^{2\lambda_{\max}},\quad \mbox{(see \eqref{DFmm-2} or \eqref{DFmm-1})}
\\
&E_A:=E_-(\omega,x),\quad E_B:=E_-(\omega,y),\quad \mbox{(see \eqref{ED-md})}
\\
&C:=2K_\epsilon(\omega),\quad \lambda:=e^{\lambda_{j+1}+\epsilon},\quad \mu:=e^{\lambda_{j}-2\epsilon}.\quad \mbox{(see \eqref{ED-md})}
\end{align*}
Then, all conditions of Lemma \ref{Pesin-lm} are verified and we conclude that for all $x,y$ satisfying $\|x-y\|<1$
\begin{align*}
&{\rm dist}(E_-(\omega,x),E_-(\omega,y))
\\
&\le 12K_\epsilon(\omega)^2 e^{\lambda_{j}-\lambda_{j+1}-3\epsilon}\big(\delta(\omega)^{1/4}
\|x-y\|^{1/2}\big)^{(\lambda_{j}-\lambda_{j+1}-3\epsilon)
/(2\lambda_{\max}-\lambda_{j+1}-\epsilon)}
\\
&\le \|x-y\|^{(1/2)(\lambda_{j}-\lambda_{j+1}-3\epsilon)
/(2\lambda_{\max}-\lambda_{j+1}-\epsilon)}
\end{align*}
since $\mathfrak{M}_\epsilon(\omega)\ge \big\{12K_\epsilon(\omega)^2 e^{\lambda_{j}-\lambda_{j+1}-3\epsilon}\big\}
^{4(2\lambda_{\max}-\lambda_{j+1}-\epsilon)
/(\lambda_{j}-\lambda_{j+1}-3\epsilon)}$ by (\ref{MM-p}), i.e.,
$$
0<\delta(\omega)\le \big\{12K_\epsilon(\omega)^2 e^{\lambda_{j}-\lambda_{j+1}-3\epsilon}\big\}
^{-4(2\lambda_{\max}-\lambda_{j+1}-\epsilon)
/(\lambda_{j}-\lambda_{j+1}-3\epsilon)},
$$
and that
\begin{align*}
{\rm dist}(E_-(\omega,x),E_-(\omega,y))
&\le 12K_\epsilon(\omega)^2 e^{\lambda_{j}-\lambda_{j+1}-3\epsilon}\delta(\omega)^{(3/4)(\lambda_{j}-\lambda_{j+1}-3\epsilon)
/(2\lambda_{\max}-\lambda_{j+1}-\epsilon)}
\\
&\le \delta(\omega)^{(1/2)(\lambda_{j}-\lambda_{j+1}-3\epsilon)
/(2\lambda_{\max}-\lambda_{j+1}-\epsilon)}.
\end{align*}
A similar conclusion holds in the subspaces $E_+(\omega,x)$. Then, by choosing different $j\in \{1,...,p\}$ and intersecting appropriate subspaces, when $\|x-y\|<1$ we prove (\ref{DB-us}) with
$$
\beta_{E}:=\min_{\ell= 1,...,p-1}\frac{\lambda_{\ell}-\lambda_{\ell+1}-3\epsilon}
{6\lambda_{\max}},\qquad \delta_E(\omega):=\delta(\omega)^{\beta_{E}},
$$
where $\beta_{E}\in (0,1)$, as known from (\ref{def-ep}). When $\|x-y\|\ge 1$, the first inequality of (\ref{DB-us}) can be indicated by the second inequality of (\ref{DB-us}) as $\delta_E(\omega)<1$. This completes the proof of Theorem \ref{thm-HD}. \qquad $\Box$

The original idea of the above proof comes from \cite{BS} for the invariant distributions of Anosov diffeomorphims, i.e., hyperbolic diffeomorphims defined on manifolds, from which we see that the classical method only gives H\"older continuity of distributions in the base point $x$. In fact, it is widely believed that the distributions are not $C^1$ (not even Lipschitzian) smooth with respect to $x$ in general without the bunching condition (see e.g. \cite{Hass-ETDS94,KB-Book95,Pesin-book}). However, in the following {Theorem \ref{thm-globC1}}, we will show that in the local case the stable and unstable distributions are $C^{1,\beta_\alpha}$ under some mild assumptions
(i.e., the assumptions of Theorem \ref{Main-ThmC1-RS}).

In order to prove {Theorem \ref{thm-globC1}}, we need to study the intermediate foliations in section 4 and then decompose a cohomological equation for the distribution in every intermediate foliation in section 5, where two lemmas concerning random normal form and the decomposed cohomological equation will be proved in sections 6-8.

\section{Random invariant foliation}
\setcounter{equation}{0}
In this section, we present the results on various invariant foliations for random dynamical systems according to the splits of Lyapunov exponents. Some of results are borrowed  from \cite{Arn98, LiLu2,LL, SLZ}.

First, from \cite[pp.362-363]{Arn98} and \cite[Theorem 3.8 and 3.9]{LiLu2},   there are the {\it random unstable foliation} $\{{\mathcal W}_u(\omega,x)\}_{x\in \mathbb{R}^d}$
and the {\it random stable foliation} $\{{\mathcal W}_s(\omega,x)\}_{x\in \mathbb{R}^d}$ of $\mathbb{R}^d$ with
{\it leaves} ${\mathcal W}_u(\omega,x)$ and ${\mathcal W}_s(\omega,x)$, which are
measurable with respect to $\omega$ and continuous with respect to $x$, such that
$\mathbb{R}^d=\cup_{x\in \mathbb{R}^d}{\mathcal W}_u(\omega,x)=\cup_{x\in \mathbb{R}^d}{\mathcal W}_s(\omega,x)$
and
\begin{align}
F(\omega,{\mathcal W}_u(\omega, x))\subset {\mathcal W}_u(\theta\omega,F( \omega, x)),\quad
F(\omega,{\mathcal W}_s(\omega, x))\subset {\mathcal W}_s(\theta\omega,F( \omega, x)).
\label{inv-fol}
\end{align}
Notice that (\ref{inv-fol}) means the invariance of foliations under the action of $F(\omega, \cdot)$. Let $\pi_u:\mathbb{R}^d\to X_u$ and $\pi_s:\mathbb{R}^d\to X_s$ be the projections associated with $\mathbb{R}^d=X_u\oplus X_s$ and let
\begin{align*}
\Lambda_* (\omega):=\Lambda(\omega) |_{X_*},\quad \Lambda_* (n,\omega) :=\Lambda (n,\omega)|_{X_*}\qquad {\rm for} ~~ *=u,s.
\end{align*}
Consider the following random Lyapunov-Perron equations (abbreviated as the L-P equations) for constructing the unstable and stable foliations respectively.
\begin{align}
\begin{split}
p_n(\omega,x,y_u)
&=\Lambda_u(n,\omega)(y_u-\pi_u x)
-\sum_{k=n}^{-1}\Lambda_u({n-k-1},\theta^{k+1}\omega)
\,\Xi_{k}^{u}(p_k)(\omega,x,y_u)
\\
&\quad+\sum_{k=-\infty}^{n-1}\Lambda_s({n-k-1},\theta^{k+1}\omega)
\,\Xi_{k}^{s}(p_k)(\omega,x,y_u),\qquad  \forall n< 0,
\\
p_0(\omega,x,y_u)
&=y_u-\pi_u x+\sum_{k=-\infty}^{-1}\Lambda_s({-k-1},\theta^{k+1}\omega)
\,\Xi_{k}^{s}(p_k)(\omega,x,y_u),\qquad  n=0,
\end{split}
\label{eqns-foli-U-RS}
\end{align}
where
$
\Xi_{k}^{*}(p_k)(\omega,x,y_u)
:=\pi_*f\big(\theta^k \omega,p_k(\omega, x,y_u)+F(k,\omega,x)\big)
-\pi_*f\big(\theta^k\omega,F(k,\omega,x)\big)
$
for $*=u,s$, and
\begin{align}
\begin{split}
q_n(\omega,x,y_s)
&=\Lambda_s(n,\omega)(y_s-\pi_s x)
+\sum_{k=0}^{n-1}\Lambda_s({n-k-1},\theta^{k+1}\omega)
\,\Xi_{k}^{s}(q_k)(\omega,x,y_s)
\\
&\quad-\sum_{k=n}^{\infty}\Lambda_u({n-k-1},\theta^{k+1}\omega)
\,\Xi_{k}^{u}(q_k)(\omega,x,y_s), \qquad  \forall n> 0,
\\
q_0(\omega,x,y_s)
&=y_s-\pi_s x
-\sum_{k=0}^{\infty}\Lambda_u({-k-1},\theta^{k+1}\omega)
\,\Xi_{k}^{u}(q_k)(\omega,x,y_s), \qquad  n=0,
\end{split}
\label{eqns-foli-S-RS}
\end{align}
where
$
\Xi_{k}^{*}(q_k)(\omega,x,y_s)
:=\pi_*f\big(\theta^k \omega,q_k(\omega, x,y_s)+F(k,\omega,x)\big)
-\pi_*f\big(\theta^k\omega,F(k,\omega,x)\big)
$
for $*=u,s$. It follows from \cite[Theorem 3.8 and 3.9]{LiLu2} that there is a constant $C_\lambda>0$ depending only on the Lyapunov exponents such that if
\begin{equation}\label{Contraction}
\delta(\omega)\le C_\lambda/(K(\theta\omega)K(\omega)),
\end{equation} then both equation \eqref{eqns-foli-U-RS} and equation \eqref{eqns-foli-S-RS} have continuous solutions $(p_n(\omega,\cdot,\cdot))_{n\le 0}$ and $(q_n(\omega,\cdot,\cdot))_{n\ge 0}$ respectively. Notice that the condition \eqref{Contraction} holds since (\ref{MM-p}) and (\ref{def-d}), which insures that the associated Lyapunov-Perron operators are contraction. Furthermore, the leaves of random unstable foliation and random stable foliation are given by
\begin{align}
\begin{split}
{\mathcal W}_u(\omega,x):=\{x+p_0(\omega,x,y_u):y_u\in X_u\},
\quad {\mathcal W}_s(\omega,x):=\{x+q_0(\omega,x,y_s):y_s\in X_s\},
\end{split}
\label{def-leaf}
\end{align}
respectively and pass through the point $x\in\mathbb{R}^d$. This point $x$ is called the {\it base point} of the leaves.

Next, we consider random pseudo-unstable foliations and  pseudo-unstable foliations.
As we did in section 3, we decompose $\mathbb{R}^d=X_+\oplus X_-$ for $1\leq j<d$, where
\begin{align*}
X_+:=X_1 \oplus \cdots \oplus X_{j},\quad
X_-:=X_{j+1}\oplus \cdots \oplus X_{{p}}.
\end{align*}  Let $\pi_\pm$ be the associated projections. Set
\begin{align*}
\Lambda_\pm(\omega) :=\Lambda(\omega)|_{X_\pm}, \quad
\Lambda_\pm(n,\omega) :=\Lambda(n,\omega)|_{X_\pm},\quad
f_\pm:=\pi_\pm f.
\end{align*}

Consider the following Lyapunov-Perron equations associated with the pseudo-unstable and pseudo-stable invariant foliations respectively,
\begin{align}
\begin{split}
p_{n,j}(\omega,x,y_+)
&=\Lambda_+(n,\omega)(y_+-\pi_+ x)
-\sum_{k=n}^{-1}\Lambda_+({n-k-1},\theta^{k+1}\omega)
\,\Xi_{k}^{+}(p_{k,j})(\omega,x,y_+)
\\
&\quad+\sum_{k=-\infty}^{n-1}\Lambda_-({n-k-1},\theta^{k+1}\omega)
\,\Xi_{k}^{-}(p_{k,j})(\omega,x,y_+),\qquad  \forall n< 0,
\\
p_{0,j}(\omega,x,y_+)
&=y_+-\pi_+ x
+\sum_{k=-\infty}^{-1}\Lambda_-({-k-1},\theta^{k+1}\omega)
\,\Xi_{k}^{-}(p_{k,j})(\omega,x,y_+),\qquad  n=0,
\end{split}
\label{eqns-foli-+}
\end{align}
where
$
\Xi_{k}^{\pm}(p_{k,j})(\omega,x,y_+)
:=f_\pm\big(\theta^k \omega,p_{k,j}(\omega, x,y_+)+F(k,\omega,x)\big)
-f_\pm\big(\theta^k\omega,F(k,\omega,x)\big),
$
and
\begin{align}
\begin{split}
q_{n,j}(\omega,x,y_-)
&=\Lambda_-(n,\omega)(y_--\pi_- x)
+\sum_{k=0}^{n-1}\Lambda_-({n-k-1},\theta^{k+1}\omega)
\,\Xi_{k}^{-}(q_{k,j})(\omega,x,y_-)
\\
&\quad-\sum_{k=n}^{\infty}\Lambda_+({n-k-1},\theta^{k+1}\omega)
\,\Xi_{k}^{+}(q_{k,j})(\omega,x,y_-), \qquad  \forall n\ge 0,
\\
q_{0,j}(\omega,x,y_-)
&=y_--\pi_- x
-\sum_{k=0}^{\infty}\Lambda_+({-k-1},\theta^{k+1}\omega)
\,\Xi_{k}^{+}(q_{k,j})(\omega,x,y_-), \qquad  n=0,
\end{split}
\label{eqns-foli--}
\end{align}
where
$
\Xi_{k}^{\pm}(q_{k,j})(\omega,x,y_-)
:=f_\pm\big(\theta^k \omega,q_{k,j}(\omega, x,y_-)+F(k,\omega,x)\big)
-f_\pm\big(\theta^k\omega,F(k,\omega,x)\big).
$ Using the same arguments used in \cite[Lemma 4.5]{SLZ},
one has that  equations  \eqref{eqns-foli-+} and  \eqref{eqns-foli--} have continuous solutions
\[(p_{n,j}(\omega,\cdot,\cdot))_{n\le 0}\in C^-_\varrho:=\big\{\{p_n\}_{n\leq 0}:\; p_n\in \mathbb{R}^d, \;\;\sup_{n\leq 0}\{e^{-\varrho n}\|p_n\|\}<\infty \big\}
\] and
\[(q_{n,j}(\omega,\cdot,\cdot))_{n\ge 0}\in C^+_\varrho:=\big\{\{q_n\}_{n\ge 0}:\; q_n\in \mathbb{R}^d,\;\; \sup_{n\ge 0}\{e^{-\varrho n}\|q_n\|\}<\infty\big\}
\]respectively, where $\varrho\in (\lambda_{j+1}+3\epsilon,\lambda_{j}-3\epsilon)$ and condition \eqref{Contraction} is used to show that the associated Lyapunov-Perron operators are contraction in  the Banach spaces $C^-_\varrho$ and $C^+_\varrho$ respectively.

Using these solutions, one obtains the following:
\begin{itemize}
\item[$\bullet$]
The strong-unstable (including only unstable directions) foliations $\{{\mathcal W}_{\le j}^{uu}(\omega,x)\}_{x\in\mathbb{R}^d}$ for $j\le \tau$.

\item[$\bullet$]
The pseudo-unstable (including unstable directions and stable ones with larger Lyapunov exponents) foliations $\{{\mathcal W}_{\le j}^{pu}(\omega,x)\}_{x\in\mathbb{R}^d}$ for $j> \tau$.

\item[$\bullet$]
The pseudo-stable (including stable directions and unstable ones with smaller Lyapunov exponents) foliations $\{{\mathcal W}_{\ge j+1}^{ps}(\omega,x)\}_{x\in\mathbb{R}^d}$ for $j< \tau$.

\item[$\bullet$]
The strong-stable (including only stable directions) foliations $\{{\mathcal W}_{\ge j+1}^{ss}(\omega,x)\}_{x\in\mathbb{R}^d}$ for $j\ge \tau$.
\end{itemize}
The leaves of these foliations are given  by
\begin{align}
\begin{split}
&{\mathcal W}_{\le j}^{uu}(\omega,x):=\{x+p_{0,j}(\omega,x,y_+)\in \mathbb{R}^d: y_+\in X_+\},\quad\forall j\le \tau,
\\
&{\mathcal W}_{\le j}^{pu}(\omega, x):=\{x+p_{0,j}(\omega,x,y_+)\in \mathbb{R}^d: y_+\in X_+\},\quad\forall j>\tau,
\\
&{\mathcal W}_{\ge j+1}^{ps}(\omega, x):=\{x+q_{0,j}(\omega, x,y_-)\in \mathbb{R}^d: y_-\in X_-\},\quad\forall j<\tau,
\\
&{\mathcal W}_{\ge j+1}^{ss}(\omega, x):=\{x+q_{0,j}(\omega, x,y_-)\in \mathbb{R}^d: y_-\in X_-\},\quad\forall j\ge \tau,
\end{split}
\label{pqpq}
\end{align}
respectively.  Note that the case $j=\tau$ are the random unstable and stable foliations. The results on random strong unstable foliations and strong stable foliations can also be found in \cite{LL}.

The following theorem is on the regularity of the leaves and their continuity with respect to the base point $x$.

\begin{thm}\label{thm-PSU}
Suppose that the random diffeomorphism $F:\Omega \times \mathbb{R}^d\to \mathbb{R}^d$ satisfies {\rm (\ref{c-off-RS})}. Then, the measurable mappings $p_{0,j}(\omega,x,\cdot)$ {\rm (}$j\le \tau${\rm )} and $q_{0,j}(\omega,x,\cdot)$ {\rm (}$j\ge \tau${\rm )}
are tempered $C^{2,\alpha}$ ${\rm(}\alpha\in[0,1]{\rm)}$ and $p_{0,j}(\omega,x,\cdot)$ {\rm (}$j> \tau${\rm )} and $q_{0,j}(\omega,x,\cdot)$ {\rm (}$j< \tau${\rm )}
are tempered $C^{1,\beta}$ with a constant $\beta\in (0,1]$.
All the above derivatives and H\"older coefficients are uniformly bounded with respect to $x,y_-,y_+$, and
\begin{align}
\begin{split}
&\|\partial_{y_+}p_{0,j}(\omega, x,y_+)-{\rm id}_+\|
\le C_\lambda K_\epsilon(\omega)
C_\epsilon(\omega)\delta(\omega),
\\
&\|\partial_{y_-}q_{0,j}(\omega, x,y_-)-{\rm id}_-\|
\le C_\lambda K_\epsilon(\omega)
C_\epsilon(\omega)\delta(\omega)
\end{split}
\label{D-qq12-tild}
\end{align}
for all $j=1,...,p$ and all $(\omega, x,y_\pm)\in \Omega\times \mathbb{R}^d\times X_\pm$, where
$C_\lambda>0$ is a constant, ${\rm id}_\pm$ is the identity in $X_\pm$, and $C_\epsilon(\omega):=4(K_\epsilon(\omega))^2M_\epsilon(\omega)$.
Moreover, the following results hold:
\begin{itemize}
\item[(A1)]
For every $\omega\in\Omega$, the mappings $x\mapsto p_{0,j}(\omega,x,\cdot)$ and $x\mapsto q_{0,j}(\omega,x,\cdot)$ are uniformly continuous in the $C^{1,\beta}$-topology on $U_1$;

\item[(A2)]
the mappings $x\mapsto \partial^2_{y_+^2}p_{0,j}(\omega,x,y_+)$ {\rm (}$j\le \tau${\rm )} and $x\mapsto \partial^2_{y_-^2}q_{0,j}(\omega,x,y_-)$ {\rm (}$j\ge \tau${\rm )}
are H\"older continuous with uniformly bounded H\"older coefficients with respect to $x,y_+,y_-$ on $U_1$,
\end{itemize}
where $U_1:=\{x\in\mathbb{R}^d: \|x\|< 1/2\}$.
\end{thm}

\begin{re}
{\rm (i)} In the theorem above, we used the following simpler notations:
\[
\partial_{y_+}p_{0,j}(\omega, x,y_+):=\frac{\partial p_{0,j}(\omega, x,y_+)}{\partial {y_+}},\qquad
\partial^2_{y_+^2}p_{0,j}(\omega,x,y_+)
:=\frac{\partial^2p_{0,j}(\omega,x,y_+)}{\partial y_+^2}.
\]
\begin{itemize}
\item[(ii)] The most part of this theorem for pseudo-stable foliations was established in \cite{SLZ} except estimate \eqref{D-qq12-tild} and the uniform boundedness of derivatives with respect to $x,y_+,y_-$.
\item[(iii)] In both pseudo-unstable and pseudo-stable cases, this theorem gives that every leaf of the foliations is $C^{1,\tilde \beta}$ for all $\tilde \beta\in (0,\beta]$, where $\beta\in (0,1]$ is determined by the spectral gap between $\lambda_j$ and $\lambda_{j+1}$. In fact, according to \cite[Theorem 3.1]{SLZ}, the $\beta$ can be chosen as
\begin{align}
\beta=\min\Bigg\{\min_{j=1,...,p-1}\frac{\lambda_{j}-\lambda_{j+1}}
{\lambda_{\max}},
~
\min_{\mbox{\tiny
$\begin{array}{c}
i\in\{1,...,\tau\},\kappa\in\{\tau+1,...,{p}\},
\\
j\in\{1,...,{p}\}
\end{array}$}}
\frac{|\lambda_i+\lambda_\kappa-\lambda_j|}{2\lambda_{\max}}\Bigg\}
\label{def-bbb}
\end{align}
with $\lambda_{\max}:=\max\{2\lambda_1,-2\lambda_p\}$.
\item[(iv)] Mapping $x\mapsto h(x,\cdot)$ is continuous in the $C^{1,\beta}$-topology means that $h(x,y)$ is $C^1$ with respect to $y$ such that
\begin{align*}
&\|h(x,\cdot)-h(\tilde x,\cdot)\|_{C^{1,\beta}}
\\
&:=\|h(x,\cdot)-h(\tilde x,\cdot)\|_{C^1}+\sup_{y,\tilde y\in \mathbb{R}^d}\frac{\|\{\partial_yh(x,y)-\partial_yh(\tilde x,y)\}-\{\partial_yh(x,\tilde y)-\partial_yh(\tilde x,\tilde y)\}\|}{\|y-\tilde y\|^\beta}\to 0
\end{align*}
as $\|x-\tilde x\|\to 0$.
\end{itemize}
\end{re}

{\bf Proof of Theorem \ref{thm-PSU}}. We prove the theorem only for  $q_{0,j}$ since the proof for $p_{0,j}$ are almost the same.   First, using (\ref{RS-ED}),  we have the following exponential dichotomy:
\begin{align}
\begin{split}
&\|\Lambda(n,\omega)|_{X_+}\| \leq K_\epsilon(\omega)
e^{(\lambda_j-\epsilon)n},
~~ \forall n< 0,
\\
&\|\Lambda(n,\omega)|_{X_-}\| \leq  K_\epsilon(\omega)
e^{(\lambda_{j+1}+\epsilon)n},
~~\forall n\geq 0,
\end{split}
\label{ED+-}
\end{align}
where $K_\epsilon(\omega)\ge K(\omega)$.
Then, using \cite[Lemma 4.5]{SLZ}, we obtain that equation
(\ref{eqns-foli--}) has a unique solution
$(q_{n,j}(\omega,x,y_-))_{n\ge 0}$, measurable with respect to $\omega$ and tempered $C^{2,\alpha}$ with respect to $y_-$ when $j\ge \tau$ (or tempered $C^{1,\beta}$ when $j<\tau$), such that
\begin{align*}
&\sup_{n\ge 0}\{e^{-\varrho n}\|q_{n,j}(\omega,x,y_-)\|\}<\infty,
\qquad
\sup_{n\ge 0}\{e^{-\varrho n}\|\partial_{y_-}q_{n,j}(\omega,x,y_-)\|\}<\infty,
\\
&\sup_{n\ge 0}\Bigg\{e^{-\varrho n}\frac{\|\partial_{y_-}q_{n,j}(\omega,x,y_-)-\partial_{y_-}q_{n,j}(\omega,x,\tilde y_-)\|}{\|y_--\tilde y_-\|^\beta}\Bigg\}<\infty.
\end{align*}

Furthermore,
differentiating both sides of equation (\ref{eqns-foli--}) with respect to $y_-$, we get
\begin{align}
\begin{split}
&\partial_{y_-}q_{n,j}(\omega, x, y_-)
\\
&=\Lambda_-(n,\omega)
+\sum_{k=0}^{n-1}\Lambda_-(n-k-1,\theta^{k+1}\omega)Df_-(\theta^k \omega,q_{k,j}(\omega, x, y_-)+F(k,\omega,x))\,\partial_{y_-}q_{k,j}(\omega, x, y_-)
\\
&\quad-\sum_{k=n}^{\infty}\Lambda_+(n-k-1,\theta^{k+1}\omega)Df_+(\theta^k \omega,q_{k,j}(\omega, x, y_-)+F(k,\omega,x))\,\partial_{y_-}q_{k,j}(\omega, x, y_-),\quad \forall n>0,
\\
&\partial_{y_-}q_{0,j}(\omega, x, y_-)
\\
&={\rm id}_-
-\sum_{k=0}^{\infty}\Lambda_+(-k-1,\theta^{k+1}\omega)Df_+(\theta^k \omega,q_{k,j}(\omega, x, y_-)+F(k,\omega,x))\,\partial_{y_-}q_{k,j}(\omega, x, y_-),\quad n=0.
\end{split}
\label{L-P-SD}
\end{align}
Then, for an arbitrarily given constant $\varrho\in (\lambda_{j+1}+3\epsilon,\lambda_{j}-3\epsilon)$, we have
\begin{align*}
&e^{-\varrho n}\|\partial_{y_-}q_{n,j}(\omega, x, y_-)\|
\\
&\le e^{-\varrho n}\|\Lambda_-(n,\omega)\|
+e^{-\varrho }\sum_{k=0}^{n-1}e^{-\varrho (n-k-1)}\|\Lambda_-(n-k-1,\theta^{k+1}\omega)\|\,
\delta(\theta^k \omega)e^{-\varrho  k}\|\partial_{y_-}q_{k,j}(\omega, x, y_-)\|
\\
&\quad-e^{-\varrho}\sum_{k=n}^{\infty}e^{-\varrho (n-k-1)}\|\Lambda_+(n-k-1,\theta^{k+1}\omega)\|
\delta(\theta^k \omega)e^{-\varrho  k}\|\partial_{y_-}q_{k,j}(\omega, x, y_-)\|
\\
&\le K_\epsilon(\omega)+C_\lambda K_\epsilon(\theta^{k+1}\omega)\delta(\theta^k\omega)\sup_{k\ge 0}\big\{e^{-\varrho  k}\|\partial_{y_-}q_{k,j}(\omega, x, y_-)\|\big\}
\\
&\le K_\epsilon(\omega)+\frac{1}{2}\sup_{k\ge 0}\big\{e^{-\varrho  k}\|\partial_{y_-}q_{k,j}(\omega, x, y_-)\|\big\},
\end{align*}
where (\ref{c-off-RS}) and (\ref{ED+-}),  the fact that $C_\lambda K_\epsilon(\theta^{k+1}\omega)\delta(\theta^k\omega)\le 1/2$ (following from (\ref{MM-p})) are used.
Hence,
$$
\sup_{n\ge 0}\big\{e^{-\varrho  n}\|\partial_{y_-}q_{n,j}(\omega, x, y_-)\|\big\}
\le 2K_\epsilon(\omega).
$$
Similarly, we have that
\[
\sup_{n\ge 0}\big\{e^{-\varrho n}\|q_{n,j}(\omega, x, y_-)\|\big\}
\le 2K_\epsilon(\omega)(\|x\|+\|y_-\|)
\]
and
\[
\sup_{n\ge 0}\big\{e^{-(\varrho+2\epsilon)  n}\|\partial_{y_-}q_{n,j}(\omega,x,y_-)
-\partial_{y_-}q_{n,j}(\omega,x,\tilde y_-)\|\big\}
\le 4(K_\epsilon(\omega))^2M_\epsilon(\omega)\|y_--\tilde y_-\|^\beta.
\]

Summarizing the above discussion gives
\begin{align}
\begin{split}
&q_{0,j}(\omega,x,\pi_- x)=0, \qquad \partial_{y_-}q_{0,j}(\omega,0,0)={\rm id}_-,
\\
&\sup_{n\ge 0}\big\{e^{-\varrho n}\|q_{n,j}(\omega,x,y_-)\|\big\}\le C_\epsilon(\omega)(\|x\|+\|y_-\|),\quad \sup_{n\ge 0}\big\{e^{-\varrho n}\|\partial_{y_-}q_{n,j}(\omega,x,y_-)\|\big\}\le C_\epsilon(\omega),
\\
&\|\partial_{y_-}q_{0,j}(\omega,x,y_-)-\partial_{y_-}q_{0,j}(\omega,x,\tilde y_-)\|\le C_\epsilon(\omega)\|y_--\tilde y_-\|^\beta,
\end{split}
\label{Wgt-Nm-DD}
\end{align}
where $C_\epsilon(\omega):=4(K_\epsilon(\omega))^2M_\epsilon(\omega)>0$.


To see that \eqref{D-qq12-tild} holds, we consider the case of $n=0$ in (\ref{L-P-SD}). It follows from (\ref{ED+-}), (\ref{c-off-RS}) and (\ref{Wgt-Nm-DD}) that
\begin{align*}
&\|\partial_{y_-}q_{0,j}(\omega,x,y_-)-{\rm id}_-\|
=e^{(\varrho+2\epsilon)\,0 }\|\partial_{y_-}q_{0,j}(\omega,x,y_-)
-\partial_{y_-}q_{0,j}(\omega,{{0}},{{0}})\|
\\
&\le e^{-(\varrho+\epsilon)}\sum_{k=0}^{\infty}
e^{-\epsilon (k+1)}e^{(\varrho+2\epsilon)(k+1)}
\|\Lambda_+(-k-1,\theta^{k+1}\omega)\|
 e^{-\epsilon k}\delta(\theta^k\omega) e^{-\varrho k}\|\partial_{y_-}q_{k,j}(\omega, x, y_-)\|
\\
&\le e^{-(\varrho+\epsilon)}\sum_{k=0}^{\infty} e^{-\epsilon (k+1)}K_\epsilon(\theta^{k+1}\omega)e^{(\varrho-\lambda_{j}+3\epsilon)(k+1)}
e^{-\epsilon k}\delta(\theta^k\omega)e^{-\varrho k}\|\partial_{y_-}q_{k,j}(\omega, x, y_-)\|
\\
&\le C_\lambda K_\epsilon(\omega) C_\epsilon(\omega)\delta(\omega),
\end{align*}
where the facts that $\varrho-\lambda_{j}+3\epsilon<0$ and
$
e^{-\epsilon |k|}\delta(\theta^k\omega)\le\delta(\omega)
$
(by (\ref{def-d})) are used. This proves the second inequality of (\ref{D-qq12-tild}).

Next, we prove (A1) for the mapping $x\mapsto q_{0,j}(\omega,x,\cdot)$. For the purpose,
we claim that
\begin{align}
\sup_{n\ge 0}\Big\{e^{-\lambda_{\max} n}\|q_{n,j}(\omega,x,y_-)-q_{n,j}(\omega,\tilde x,y_-)\|\Big\}
\le L_\epsilon(\omega)\|x-\tilde x\|^\sigma
\label{Hd-qq-1}
\end{align}
for all $x,\tilde x, y_-, \tilde y_-\in U_1$, where $L_\epsilon(\omega)>0$ is an $\epsilon$-slowly varying random variable and $\sigma\in (0,1)$ is a constant such that $\sigma<(\lambda_{j}-\lambda_{j+1})/(4\lambda_{\max})$. In order to prove (\ref{Hd-qq-1}), we notice from (\ref{DF-max}) that
\begin{align}
e^{-\lambda_{\max}k}\|
F(k,\omega,x)-F(k,\omega,\tilde x)\|
&\le
e^{-\lambda_{\max}k}\sup_{\xi\in\mathbb{R}^d}\|DF(k,\omega,\xi)\|
\,\|x-\tilde x\|
\nonumber\\
&\le K_\epsilon(\omega)\|x-\tilde x\|,\qquad \forall k\ge 0.
\label{Fmax}
\end{align}
Using the notation $[h(\cdot)]_{\tilde x}^x:=h(x)-h(\tilde x)$ for short, we have
\begin{align}
&\big[f(\theta^k \omega,q_{k,j}(\omega,\cdot,y_-)+F(k,\omega,\cdot))-f(\theta^k \omega,F(k,\omega,\cdot))\big]_{\tilde x}^{ x}
\nonumber\\
&=\big[f(\theta^k \omega,q_{k,j}(\omega,\cdot,y_-)+F(k,\omega,x))\big]_{\tilde x}^{ x}
+\big[f(\theta^k \omega,q_{k,j}(\omega,\tilde x,y_-)+F(k,\omega,\cdot))\big]_{\tilde x}^{ x}
\nonumber\\
&\quad -\big[f(\theta^k\omega,F(k,\omega,\cdot))\big]_{\tilde x}^x
\nonumber\\
&=\big[f(\theta^k \omega,q_{k,j}(\omega,\cdot,y_-)+F(k,\omega,x))\big]_{\tilde x}^{ x}
\nonumber\\
&\quad+\int_{0}^{1}\big[Df(\theta^k \omega,tq_{k,j}(\omega,\tilde x,y_-)+F(k,\omega,\cdot))q_{k,j}(\omega,\tilde x,y_-)\big]_{\tilde x}^xdt.
\label{defrnc}
\end{align}
Then, setting
$$
{\tilde \lambda}:=\lambda_{j+1}+(\lambda_{j}-\lambda_{j+1})/2,
\qquad
\|h(\cdot)\|_{\tilde x}^x:=\|h(x)-h(\tilde x)\|
$$
for short, we obtain from (\ref{defrnc}), (\ref{c-off-RS}), (\ref{Wgt-Nm-DD}) and (\ref{Fmax}) that
\begin{align}
&e^{-{\tilde \lambda} k-(\lambda_{\max}\sigma+\epsilon) k}\big\|f(\theta^k \omega,q_{k,j}(\omega,\cdot,y_-)+F(k,\omega,\cdot))-f(\theta^k \omega,F(k,\omega,\cdot))\big\|_{\tilde x}^{ x}
\nonumber\\
&\le\sup_{\xi\in \mathbb{R}^d}\|Df(\theta^k \omega,\xi)\|\,
e^{-{\tilde \lambda} k-(\lambda_{\max}\sigma+\epsilon) k}\big\|q_{k,j}(\omega,\cdot,y_-)\big\|_{\tilde x}^{x}
\nonumber\\
&\quad+e^{-(\lambda_{\max}\sigma+\epsilon) k}\Big(\big\|Df(\theta^k \omega,tq_{k,j}(\omega,\tilde x,y_-)+F(k,\omega,\cdot))\big\|
_{\tilde x}^{x}\Big)^\sigma e^{-{\tilde \lambda} k}\| q_{k,j}(\omega,\tilde x,y_-)\|
\nonumber\\
&\le \delta(\theta^k\omega)e^{-{\tilde \lambda} k-(\lambda_{\max}\sigma+\epsilon) k}\big\|q_{k,j}(\omega,\cdot,y_-)\big\|_{\tilde x}^{x}
+
e^{-\epsilon k}M_\epsilon(\theta^k\omega)^\sigma
\Big(e^{-\lambda_{\max}k}\|
F(k,\omega,\cdot)\big\|_{\tilde x}^{x}\Big)^{\sigma}
C_\epsilon(\omega)(\|\tilde x\|+\|y_-\|)
\nonumber\\
&\le \delta(\theta^k\omega)e^{-{\tilde \lambda} k-(\lambda_{\max}\sigma+\epsilon) k}\big\|q_{k,j}(\omega,\cdot,y_-)\big\|_{\tilde x}^{x}\,+\tilde M_\epsilon(\omega)\|x-\tilde x\|^\sigma
\label{HH-L}
\end{align}
for all $k\ge 0$ and all $\tilde x,y_-\in U_1$, where $\tilde M_\epsilon(\omega):=2M_\epsilon(\omega)^\sigma
K_\epsilon(\omega)^\sigma C_\epsilon(\omega)>0$ and the following estimate is used:
\[
\big\|Df(\theta^k \omega,tq_{k,j}(\omega,\tilde x,y_-)+F(k,\omega,\cdot))\big\|
_{\tilde x}^{x}
\le
\Big(\big\|Df(\theta^k \omega,tq_{k,j}(\omega,\tilde x,y_-)+F(k,\omega,\cdot))\big\|
_{\tilde x}^{x}\Big)^\sigma
\]
for $\sigma\in (0,1)$ since
$\|Df(\theta^k \omega,x)\|\le \delta(\theta^k\omega)<1/2$ due to (\ref{c-off-RS}).

Notice that ${\tilde \lambda} +\lambda_{\max}\sigma+2\epsilon\in (\lambda_{j+1}+3\epsilon,\lambda_{j}-3\epsilon)$ with small $\epsilon$ due to the above choices of ${\tilde \lambda}$ and $\sigma$, which implies that
$$
\sup_{k\ge 0}\Big\{e^{-{\tilde \lambda} k-(\lambda_{\max}\sigma+2\epsilon) k}\big\|q_{k,j}(\omega,\cdot,y_-)\big\|_{\tilde x}^{x}\Big\}<\infty
$$
by (\ref{Wgt-Nm-DD}).
Hence, we see from (\ref{eqns-foli--}), (\ref{ED+-}) and (\ref{HH-L}) that
for $n\ge 0$
\begin{align}
&e^{-{\tilde \lambda} n-(\lambda_{\max}\sigma+2\epsilon) n}\|q_{n,j}(\omega,x,y_-)-q_{n,j}(\omega,\tilde x,y_-)\|
\nonumber\\
&\le e^{-{\tilde \lambda} n-(\lambda_{\max}\sigma+2\epsilon) n}\|\Lambda_-(n,\omega)\|\|x-\tilde x\|
\nonumber\\
&\quad+e^{-{\tilde \lambda}-(\lambda_{\max}\sigma+2\epsilon)}
\sum_{k=0}^{n-1}e^{-{\tilde \lambda} (n-k-1)-(\lambda_{\max}\sigma+2\epsilon) (n-k-1)}\|\Lambda_-(n-k-1,\theta^{k+1}\omega)\|
\nonumber\\
&\quad
\cdot e^{-{\tilde \lambda} k-(\lambda_{\max}\sigma+2\epsilon) k}\big\|f_-(\theta^k \omega,q_{k,j}(\omega,\cdot,y_-)+F(k,\omega,\cdot))-f_-(\theta^k \omega,F(k,\omega,\cdot))\big\|_{\tilde x}^{ x}
\nonumber\\
&\quad+e^{-{\tilde \lambda}-(\lambda_{\max}\sigma+2\epsilon)}
\sum_{k=n}^{\infty}e^{-{\tilde \lambda} (n-k-1)-(\lambda_{\max}\sigma+2\epsilon) (n-k-1)}\|\Lambda_+(n-k-1,\theta^{k+1}\omega)\|
\nonumber\\
&\quad\cdot e^{-{\tilde \lambda} k-(\lambda_{\max}\sigma+2\epsilon) k}\big\|f_+(\theta^k \omega,q_{k,j}(\omega,\cdot,y_-)+F(k,\omega,\cdot))-f_+(\theta^k \omega,F(k,\omega,\cdot))\big\|_{\tilde x}^{ x}
\nonumber\\
&\le K_{\epsilon}(\omega)e^{(\lambda_{j+1}-{\tilde \lambda}-
\lambda_{\max}\sigma-\epsilon) n}\|x-\tilde x\|^\sigma
\nonumber\\
&\quad+e^{-{\tilde \lambda}-(\lambda_{\max}\sigma+2\epsilon)}
\bigg(\sum_{k=0}^{n-1}K_\epsilon(\theta^{k+1}\omega)
e^{(\lambda_{j+1}-{\tilde \lambda}-\lambda_{\max}\sigma-\epsilon) (n-k-1)}
\nonumber\\
&\quad+
\sum_{k=n}^{\infty}K_\epsilon(\theta^{k+1}\omega)
e^{(\lambda_{j}-{\tilde \lambda}-\lambda_{\max}\sigma-3\epsilon) (n-k-1)}\bigg)
\nonumber\\
&\quad \cdot \sup_{k\ge 0}\Big\{\delta(\theta^k\omega)e^{-{\tilde \lambda} k-(\lambda_{\max}\sigma+2\epsilon) k}\big\|q_{k,j}(\omega,\cdot,y_-)\big\|_{\tilde x}^{x}\,+e^{-\epsilon k}\tilde M_\epsilon(\omega)\|x-\tilde x\|^\sigma\Big\}
\nonumber\\
&\le K_{\epsilon}(\omega)\|x-\tilde x\|^\sigma
\nonumber\\
&\quad +C_\lambda K_{\epsilon}(\theta^{k+1}\omega)\sup_{k\ge 0}\Big\{\delta(\theta^k\omega)e^{-{\tilde \lambda} k-(\lambda_{\max}\sigma+2\epsilon) k}\big\|q_{k,j}(\omega,\cdot,y_-)\big\|_{\tilde x}^{x}\,+e^{-\epsilon k}\tilde M_\epsilon(\omega)\|x-\tilde x\|^\sigma\Big\}
\nonumber\\
&\le 2C_\lambda K_\epsilon(\omega)\tilde M_\epsilon(\omega)\|x-\tilde x\|^\sigma+\frac{1}{2} \sup_{k\ge 0}\Big\{e^{-{\tilde \lambda} k-(\lambda_{\max}\sigma+2\epsilon) k}\big\|q_{k,j}(\omega,\cdot,y_-)\big\|_{\tilde x}^{x}\Big\}
\label{LL-F}
\end{align}
when $\|x-\tilde x\|<1$. Here,
the facts $\mathfrak{M}_\epsilon(\omega)\ge 2C_\lambda K_\epsilon(\theta\omega)$ (i.e., $\delta(\theta^k\omega)\le (2C_\lambda K_\epsilon(\theta^{k+1}\omega))^{-1}$ by (\ref{MM-p})) and $\lambda_{j+1}-{\tilde \lambda}-\lambda_{\max}\sigma-\epsilon<0$, $\lambda_{j}-{\tilde \lambda}-\lambda_{\max}\sigma-3\epsilon>0$ (see the above choices of ${\tilde \lambda}$ and $\sigma$) are used. Hence, we have
\begin{align*}
\sup_{n\ge 0}\Big\{e^{-{\tilde \lambda} n-(\lambda_{\max}\sigma+2\epsilon) n}\|q_{n,j}(\omega,x,y_-)-q_{n,j}(\omega,\tilde x,y_-)\|\Big\}
\le 4C_\lambda K_\epsilon(\omega)\tilde M_\epsilon(\omega)\|x-\tilde x\|^\sigma,
\end{align*}
which proves (\ref{Hd-qq-1}) by setting an $\epsilon$-slowly varying random variable $L_\epsilon(\omega)\ge 4C_\lambda K_\epsilon(\omega)\tilde M_\epsilon(\omega)>0$
since $e^{-\lambda_{\max} n}\le e^{-{\tilde \lambda} n-(\lambda_{\max}\sigma+2\epsilon) n}$.

Next, we prove that
\begin{align}
\sup_{n\ge 0}\Big\{e^{-{\tilde \lambda} n-(\lambda_{\max}\sigma+2\epsilon) n}\big\|\partial_{y_-}q_{n,j}(\omega, x, y_-)-
\partial_{y_-}q_{n,j}(\omega, \tilde x, y_-)\big\|\Big\}
\le \tilde L_\epsilon(\omega)\|x-\tilde x\|^{\sigma^2},
\label{Hd-qq-2}
\end{align}
where $\tilde L_\epsilon(\omega):= 2C_\lambda K_\epsilon(\omega)M_\epsilon(\omega)
(L_\epsilon(\omega)+K_\epsilon(\omega))^\sigma C_\epsilon(\omega)>0$. In fact,
notice that by (\ref{Hd-qq-1}) and (\ref{Fmax}) we have
\begin{align}
&e^{-(\lambda_{\max}\sigma +\epsilon)k}\big\|Df(\theta^k \omega,q_{k,j}(\omega, \cdot, y_-)+F(k,\omega,\cdot))\big\|^{x}_{\tilde x}
\nonumber\\
&=
\Big(\big\|Df(\theta^k \omega,q_{k,j}(\omega, \cdot, y_-)+F(k,\omega,\cdot))\big\|^{x}_{\tilde x}\Big)^{1-\sigma}
\nonumber\\
&\quad\cdot e^{-\epsilon k}\Big(e^{-\lambda_{\max} k}\big\|Df(\theta^k \omega,q_{k,j}(\omega, \cdot, y_-)+F(k,\omega,\cdot))\big\|^{x}_{\tilde x}\Big)^\sigma
\nonumber\\
&\le (2\delta(\theta^k\omega))^{1-\sigma}e^{-\epsilon k}M_\epsilon(\theta^k\omega)^{\sigma}
\Big(e^{-\lambda_{\rm max} k}\big\|q_{k,j}(\omega, \cdot, y_-)\big\|^{x}_{\tilde x}+e^{-\lambda_{\rm max} k}\big\|F(k,\omega,\cdot)\big\|^{x}_{\tilde x}\Big)^\sigma
\nonumber\\
&\le  M_\epsilon(\omega)
(L_\epsilon(\omega)+K_\epsilon(\omega))^\sigma\|x-\tilde x\|^{\sigma^2}
\label{Df-Df}
\end{align}
when $\|x-\tilde x\|<1$.
Then, using similar arguments to (\ref{LL-F}), we obtain from (\ref{L-P-SD}) that
\begin{align}
&e^{-{\tilde \lambda} n-(\lambda_{\max}\sigma+2\epsilon) n}\|\partial_{y_-}q_{n,j}(\omega, x, y_-)-\partial_{y_-}q_{n,j}(\omega, \tilde x, y_-)\|
\nonumber\\
&= C_\lambda K_\epsilon(\omega)\sup_{k\ge 0}\bigg\{ e^{-(\lambda_{\max}\sigma +2\epsilon) k}\big\|Df(\theta^k \omega,q_{k,j}(\omega, \cdot, y_-)+F(k,\omega,\cdot))\big\|^{x}_{\tilde x}
\nonumber\\
&\quad\cdot e^{-{\tilde \lambda} k}\|\partial_{y_-}q_{k,j}(\omega, x, y_-)\|
+\|Df(\theta^k \omega,q_{k,j}(\omega, \tilde x, y_-)+F(k,\omega,\tilde x))\|
\nonumber\\
&\quad \cdot e^{-{\tilde \lambda} k-(\lambda_{\max}+2\epsilon) k}\big\|\partial_{y_-}q_{k,j}(\omega, \cdot, y_-)\big\|^{x}_{\tilde x}\bigg\}
\nonumber\\
&\le C_\lambda K_\epsilon(\omega)M_\epsilon(\omega)
(L_\epsilon(\omega)+K_\epsilon(\omega))^\sigma C_\epsilon(\omega)\|x-\tilde x\|^{\sigma^2}
\nonumber\\
&\quad+\frac{1}{2} \sup_{k\ge 0}\Big\{e^{-{\tilde \lambda} k-(\lambda_{\max}+2\epsilon) k}\big\|\partial_{y_-}q_{k,j}(\omega, x, y_-)-
\partial_{y_-}q_{k,j}(\omega, \tilde x, y_-)\big\|\Big\},
\label{JCA}
\end{align}
which yields (\ref{Hd-qq-2}).

Setting $n=0$ in (\ref{Hd-qq-2}) we obtain that
\begin{align}
\big\|\partial_{y_-}q_{0,j}(\omega, x, y_-)-
\partial_{y_-}q_{0,j}(\omega, \tilde x, y_-)\big\|
\le \tilde L_\epsilon(\omega)\|x-\tilde x\|^{\sigma^2}.
\label{Hd-qq-22}
\end{align}
Therefore for any small constant $\varepsilon>0$, we obtain
\begin{align}
&\frac{\big\|\{\partial_{y_-}q_{0,j}(\omega, x, y_-)-\partial_{y_-}q_{0,j}(\omega, \tilde x, y_-)\}-\{\partial_{y_-}q_{0,j}(\omega, x, \tilde y_-)
-\partial_{y_-}q_{0,j}(\omega, \tilde x, \tilde y_-)\}\big\|}
{\|y_--\tilde y_-\|^{(1-\varepsilon)\beta}\|x-\tilde x\|^{\varepsilon\sigma^2}}
\nonumber\\
&\le \bigg(\frac{\|\partial_{y_-}q_{0,j}(\omega, x, y_-)-\partial_{y_-}q_{0,j}(\omega, x, \tilde y_-)\|+\|\partial_{y_-}q_{0,j}(\omega, \tilde x, y_-)
-\partial_{y_-}q_{0,j}(\omega, \tilde x, \tilde y_-)\|}
{\|y_--\tilde y_-\|^\beta}\bigg)^{1-\varepsilon}
\nonumber\\
&\quad\cdot \bigg(\frac{\|\partial_{y_-}q_{0,j}(\omega, x, y_-)-\partial_{y_-}q_{0,j}(\omega, \tilde x, y_-)\|+\|\partial_{y_-}q_{0,j}(\omega, x, \tilde y_-)
-\partial_{y_-}q_{0,j}(\omega, \tilde x, \tilde y_-)\|}{\|x-\tilde x\|^{\sigma^2}}\bigg)^{\varepsilon}
\nonumber\\
&\le 2C_\epsilon(\omega)\tilde L_\epsilon(\omega),
\label{Hd-qq-33}
\end{align}
due to (\ref{Wgt-Nm-DD}). Combining (\ref{Hd-qq-1}) (setting $n=0$),
(\ref{Hd-qq-22}) with (\ref{Hd-qq-33}) together gives that
\begin{align*}
&\|q_{0,j}(\omega, x, \cdot)-q_{0,j}(\omega, \tilde x, \cdot)\|_{C^{1,(1-\varepsilon)\beta}}
\\
&=\|q_{0,j}(\omega, x, \cdot)-q_{0,j}(\omega, \tilde x, \cdot)\|_{C^1}
\\
&\quad+\sup_{y_-,\tilde y_-\in X_-}\frac{\big\|\{\partial_{y_-}q_{0,j}(\omega, x, y_-)-\partial_{y_-}q_{0,j}(\omega, \tilde x, y_-)\}-\{\partial_{y_-}q_{0,j}(\omega, x, \tilde y_-)
-\partial_{y_-}q_{0,j}(\omega, \tilde x, \tilde y_-)\}\big\|}
{\|y_--\tilde y_-\|^{(1-\varepsilon)\beta}}
\\
&\le \{L_\epsilon(\omega)+\tilde L_\epsilon(\omega)+2C_\epsilon(\omega)\tilde L_\epsilon(\omega)\}\|x-\tilde x\|^{\varepsilon\sigma^2}
\end{align*}
when $\|x-\tilde x\|<1$.
Since $\beta$ is in an open interval determined by the spectral gaps (see (\ref{def-bbb})) and $\varepsilon>0$ can be sufficiently small, it is reasonable to use $\beta$ to denote $(1-\varepsilon)\beta$ and therefore the mapping $x\mapsto q_{0,j}(\omega,x,\cdot)$ is uniformly continuous in the $C^{1,\beta}$-topology on $U_1$.
This proves (A1) for the mapping $x\mapsto q_{0,j}(\omega,x,\cdot)$.

Similarly to (\ref{Hd-qq-22}), in the case of $j\ge \tau$, we show that
the mapping $x\mapsto \partial^2_{y_-^2}q_{0,j}(\omega,x,y_-)$
is H\"older continuous with uniformly bounded H\"older coefficients with respect to $x,y_-$ on $U_1$. This proves (A2) for the mapping $x\mapsto q_{0,j}(\omega,x,\cdot)$ by using a similar discussion to the one given below (\ref{Hd-qq-22}).
The proof of this theorem is completed. \qquad$\Box$

For the sake of convenient representations,
we use notations
\begin{align}
{\mathcal W}_{\le j}^{u}(\omega, x):=
\left\{
  \begin{array}{lll}
{\mathcal W}_{\le j}^{uu}(\omega, x),&\!j\!<\! \tau,
\vspace{1ex}\\
{\mathcal W}_u(\omega,x),&\!j\!=\!\tau,
\vspace{1ex}\\
{\mathcal W}_{\le j}^{pu}(\omega, x),&\! j\!>\!\tau,
\end{array}
\right.
\qquad
{\mathcal W}_{\ge j+1}^{s}(\omega, x):=
\left\{
  \begin{array}{lll}
{\mathcal W}_{\ge j+1}^{ps}(\omega, x),&\! j\!<\!\tau,
\vspace{1ex}\\
{\mathcal W}_s(\omega,x),&\!j\!=\!\tau,
\vspace{1ex}\\
{\mathcal W}_{\ge j+1}^{ss}(\omega, x),&\!j\!>\! \tau.
\end{array}
\right.
\label{uuss}
\end{align}
Notice that these leaves are graphs of the mappings
$$
y_+\mapsto p_{0,j}(\omega,x,y_+)+\pi_+x-y_+,\qquad
y_-\mapsto q_{0,j}(\omega,x,y_-)+\pi_-x-y_-
$$
with Lipschitz constants $<1$ by (\ref{D-qq12-tild}). Then, corresponding to each $X_j$
we introduce the leaf
\begin{align}
{\mathcal W}_j(\omega, x):=
\left\{
  \begin{array}{lll}
    {\mathcal W}_{\le 1}^{u}(\omega, x),&\quad \mbox{if } j=1,
    \vspace{1ex}\\
    {\mathcal W}_{\le j}^u(\omega, x)\cap {\mathcal W}_{\ge j}^{s}(\omega, x),&\quad \mbox{if } j=2,...,{{p}}-1,
     \vspace{1ex}\\
     {\mathcal W}_{\ge {{p}}}^{s}(\omega, x),&\quad \mbox{if }j={{p}}.
  \end{array}
\right.
\label{inter-foli}
\end{align}
Clearly, for each $j\in \{1,...,p\}$ it defines a foliation
${\mathcal W}_j(\omega):=\{{\mathcal W}_j(\omega, x)\}_{x\in \mathbb{R}^d}$,
called the {\it random intermediate foliation}, and it is the graph of a mapping $\check q_j(\omega, x, \cdot):X_j\to X_1\oplus \cdots\oplus X_{j-1}\oplus X_{j+1}\oplus \cdots\oplus X_p$ such that $\mathbb{R}^d=\cup_{x\in\mathbb{R}^d}{\mathcal W}_j(\omega, x)$ and the following corollary holds.
\begin{cor}\label{cor-PSU}
Suppose that the random diffeomorphism $F:\Omega \times \mathbb{R}^d\to \mathbb{R}^d$ satisfies {\rm (\ref{c-off-RS})}. Then for each $j=1,...,p$
the measurable mapping $\check q_j(\omega,x,\cdot)$
is tempered $C^{1,\beta}$ with a constant $\beta\in (0,1]$
such that
\begin{align}
\|\partial_{y_j}\check q_{j}(\omega, x,y_j)-{\rm id}_j\|
\le C_\lambda K_\epsilon(\omega)
C_\epsilon(\omega)\delta(\omega)
\label{D-qq12-chek}
\end{align}
for all $(\omega, x,y_j)\in \Omega\times \mathbb{R}^d\times X_j$, where ${\rm id}_j$ is the identity in $X_j$. Moreover, the mapping $x\mapsto \check q_{j}(\omega,x,\cdot)$ is uniformly continuous in the $C^{1,\beta}$-topology on $U_1:=\{x\in\mathbb{R}^d: \|x\|< 1/2\}$.
\end{cor}
Concerning the relation between the intermediate distribution
$E_j(\omega)$ and the intermediate foliation ${\mathcal W}_j(\omega)$, we further show that
\begin{align}
E_j(\omega, x)&=\{\partial_{y_j}\check q(\omega, x,\pi_j x)\,z_j\in \mathbb{R}^d:~z_j\in X_j\},\qquad \forall x\in \mathbb{R}^d,
\label{EWjj}
\end{align}
(see (\ref{EEsu}) below in the stable case for more details), where $\pi_j:\mathbb{R}^d\to X_j$ is defined by $\pi_j(x_1+\cdots +x_p):=x_j$. It means that every fiber of
$E_j(\omega)$ at $x$ is the tangent space of the leaf of ${\mathcal W}_j(\omega)$ at $x$.

\section{Decomposition of a cohomological equation}
\setcounter{equation}{0}

In order to show that the stable and unstable distributions obtained in section 3 are $C^{1,\beta_\alpha}$ with a constant $\beta_\alpha\in [0,\alpha]$, we will prove that
there are $C^{1,\beta_\alpha}$ vectors fields $v^s_j(\omega, x)$ and $v^u_j(\omega, x)$ that span $E_s(\omega, x)$ and $E_u(\omega, x)$ at every point $x\in \mathbb{R}^d$ respectively (see \cite[p.41]{War}). This will be done in sections 5-9 mainly for $E_s(\omega, x)$, as the proof for $E_u(\omega, x)$ is similar.

In this section, we first derive a cohomological equation
for the vector field $v^s_j(\omega, x)$. Then, we solve the cohomological equation in 6 steps to show that,
for any $\bar x\in\mathbb{R}^d$ and $i=1,...,\tau$, there is a basis (depending on $\bar x$ and $i$) of $E_s(\omega, x)$ which is differentiable at $\bar x$ along every leaf ${\mathcal W}_i(\omega,\bar x)$ of the intermediate foliation with a (H\"older) continuous differential.
Letting $T_{\bar x}
{\mathcal W}_i(\omega,\bar x)$ be the tangent space of ${\mathcal W}_i(\omega,\bar x)$ at $\bar x$, we have the following theorem.

\begin{thm}\label{Thm-decop}
Suppose that the random diffeomorphism $F:\Omega \times \mathbb{R}^d\to \mathbb{R}^d$ satisfies {\rm (\ref{c-off-RS})} and that the Lyapunov exponents satisfy {\rm (\ref{EC-RS})}.
Then, for any point $\bar x\in \mathbb{R}^d$ and $i= 1,...,\tau$, the fibers $E_s(\omega, x)$ have bases
$\{{v_{\bar x, i}^{\kappa,\iota}}(\omega, x)\}_{\kappa={\tau}+1,...,{p},~\iota=1,...,d_\kappa}$ {for all $x\in {\mathcal W}_i(\omega,\bar x)$ near $\bar x$}, whose differential {\rm(}with respect to $x${\rm)}
\begin{align}
dv_{\bar x, i}^{\kappa,\iota}(\omega, \bar x):T_{\bar x}
{\mathcal W}_i(\omega,\bar x)\to \mathbb{R}^d ~~(= T_{v_{\bar x, i}^{\kappa,\iota}(\omega, \bar x)}\mathbb{R}^d\,)
\label{dif-tial}
\end{align}
at the point $\bar x$ is $\beta_v$-H\"older continuous with respect to $\bar x$ for every $\kappa={\tau}+1,...,{p}$ and $\iota=1,...,d_\kappa$, where $\beta_v:=\min\{\beta_E, \epsilon/\lambda_{\max},\alpha\}$. Moreover, the differential and the H\"older coefficient are measurable with respect to $\omega$ and are uniformly bounded with respect to $\bar x$ and $i$.
\end{thm}

\begin{re} {\rm (i)} The $\beta_v$-H\"older continuity of $dv_{\bar x, i}^{\kappa,\iota}(\omega, \bar x)$ with respect to $\bar x$ means that
\begin{align*}
&\|dv_{\bar x, i}^{\kappa,\iota}(\omega, \bar x)\circ \hbar(\omega, \bar x)
-dv_{\bar y, i}^{\kappa,\iota}(\omega, \bar y)\circ \hbar(\omega, \bar y)\|\le L(\omega)\|\bar x-\bar y\|^{\beta_v}
\end{align*}
when $\|\bar x-\bar y\|<1$, where $\hbar(\omega, \cdot)$ is a homeomorphism such that $\hbar^{\pm 1}(\omega, \cdot)$ are $\beta_v$-H\"older continuous and $\hbar^{-1}(\omega, \cdot)$ trivializes the bundle
$
{\mathcal E}:=\bigcup_{\bar x\in \mathbb{R}^d}T_{\bar x}
{\mathcal W}_i(\omega,\bar x)
$0
with the base $\mathbb{R}^d$.
If $\beta_v=0$ {\rm (}when $\alpha=0${\rm )} then $\beta_v$-H\"older continuity merely means continuity, and the differential is uniformly bounded with respect to $\bar x$ and $i$.

{\rm (ii)} The bases of $E_s(\omega, x)$ are not well defined for all $x\in \mathbb{R}^d$ since $\{{v_{\bar x, i}^{\kappa,\iota}}(\omega, x)\}_{\kappa={\tau}+1,...,{p},~\iota=1,...,d_\kappa}$ depends on $\bar x$ and $i$. Therefore, it is possible that $x\in {\mathcal W}_i(\omega, \bar x)\cap {\mathcal W}_j(\omega, \bar y)$
for certain $\bar x, \bar y\in \mathbb{R}^d$ and $i,j\in \{1,...\tau\}$, but $v_{\bar x, i}^{\kappa,\iota}(\omega, x)\ne v_{\bar y, j}^{\kappa,\iota}(\omega, x)$. We will solve this problem in section 9 by finding a ``canonical" basis
of $E_s(\omega, x)$ not depending on $\bar x, i$ and giving the relation between the ``canonical" basis and $\{{v_{\bar x, i}^{\kappa,\iota}}(\omega, x)\}_{\kappa={\tau}+1,...,{{p}},~\iota=1,...,d_\kappa}$.
\end{re}

{\bf Proof of Theorem \ref{Thm-decop}}. We consider the {\it random cohomological equation} for $v(\omega, x)\in E_s(\omega, x)$, i.e.,
\begin{align}
DF(\omega,x)v(\omega,x)=v(\theta\omega, F(\omega,x)),\quad \forall \omega\in \Omega,\quad \forall x\in \mathbb{R}^d.
\label{ch-eqn}
\end{align}
In what follows, our goal is to find out the solutions $v(\omega, x)$ of equation (\ref{ch-eqn}), which form the basis of $E_s(\omega, x)$ and are differentiable at every point $x\in \mathbb{R}^d$ with (H\"older) continuous derivatives. We will show this in 6 steps.

{\bf Step 1}. {\it Given any point $\bar x\in \mathbb{R}^d$, transform the orbit $(F(n,\omega,\bar x))_{n\in\mathbb{Z}}$ to the origin $0$ to get a new system {\rm (\ref{def-hatF})}, and then reduce to its normal form {\rm (\ref{bar-hat-0})}.}

In order to study the differentiability of $v(\omega,\cdot)$ along the orbit $(F(n,\omega,\bar x))_{n\in\mathbb{Z}}$, we let
$$
\varpi:=(\omega,\bar x)\in \bar {\bf \Omega}:=\Omega\times \mathbb{R}^d,\quad \vartheta\varpi:=(\theta\omega,F(\omega,\bar x)),\quad {\bf \Lambda}(\varpi):=DF(\omega, \bar x).
$$
Then, we have
$$
\vartheta^n\varpi=(\theta^n\omega, F(n,\omega,\bar x)),\quad
{\bf \Lambda}(\vartheta^n\varpi)=DF(\theta^n\omega, F(n,\omega,\bar x)),
\quad \forall n\in \mathbb{Z}.
$$
Note that the random matrix ${\bf \Lambda}(\varpi)$ generates the linear cocycle
\begin{align}
{\bf \Lambda}(n, \varpi):=
\begin{cases} {\bf \Lambda}(\vartheta^{n-1}\varpi) \cdots {\bf \Lambda}(\varpi), & n >0, \\
{\rm id}, & n=0, \\
{\bf \Lambda}^{-1}(\vartheta^n\varpi) \cdots {\bf \Lambda}^{-1}(\vartheta^{-1}\varpi),& n<0.
\end{cases}
\label{RS-CYC}
\end{align}
It is clear that ${\bf \Lambda}(n, \varpi)=DF(n,\omega,\bar x)$. Therefore, by (\ref{int-m-d}) we have
\begin{align*}
\begin{split}
&\|{\bf \Lambda}(n, \varpi)|_{E_j(\varpi)}\| \leq K_\epsilon(\omega)
e^{(\lambda_j-\epsilon)n},
\qquad \forall n< 0,
\\
&\|{\bf \Lambda}(n, \varpi)|_{E_j(\varpi)}\| \leq  K_\epsilon(\omega)
e^{(\lambda_j+\epsilon)n},
\qquad \forall n\geq 0.
\end{split}
\end{align*}
Then, we use \cite[Corollary 4.3.12]{Arn98} to block-diagonalize ${\bf \Lambda}(n, \varpi)$, namely, there is a random matrix $P(\varpi)\in Gl(d,\mathbb{R})$ such that
\begin{align}
\bar{\bf \Lambda}(\varpi):=P(\vartheta\varpi){\bf \Lambda}(\varpi)P^{-1}(\varpi)={\rm diag}(\bar{\bf \Lambda}_{1}(\varpi),...,\bar{\bf \Lambda}_{p}(\varpi)),
\label{blo-diag}
\end{align}
and
\begin{align}
\begin{split}
&\|\bar{\bf \Lambda}_{j}(n,\varpi)\| \leq  K_\epsilon(\omega)
e^{(\lambda_j-\epsilon)n},
\qquad \forall n< 0,
\\
&\|\bar{\bf \Lambda}_{j}(n,\varpi)\| \leq K_\epsilon(\omega)
e^{(\lambda_j+\epsilon)n},
\qquad \forall n\geq 0.
\end{split}
\label{ED-tmp}
\end{align}
{Moreover precisely, $P(\varpi):=(\vec{e}_{1,1}(\varpi)
,...,\vec{e}_{p,d_p}(\varpi))$ with measurable unit vectors $\vec{e}_{j,\iota}(\varpi)$ for all $j=1,...,p$ and $\iota=1,...,d_j$ is a transition matrix from the standard basis $\mathfrak{C}$ of $\mathbb{R}^d$ to a basis $\mathfrak{C}(\varpi):=\{\vec{e}_{j,\iota}(\varpi)\}
_{j=1,...,p,\,\iota=1,...,d_j}$ associated with the measurable decompositions
$$
\mathbb{R}^d=E_1(\varpi)\oplus \cdots\oplus E_{p}(\varpi).
$$
By the second inequality of \eqref{DB-us}, we can choose $\vec{e}_{j,\iota}(\varpi)\in \mathfrak{C}(\varpi)$ that is uniformly (with respect to $\varpi$) close to the corresponding standard unit vector $\vec{e}_{j,\iota}\in \mathfrak{C}$. Thus,
for all $\varpi:=(\omega,\bar x)$ and $\tilde \varpi:=(\omega,\bar y)\in \bar {\bf \Omega}$,
\begin{align*}
\|\vec{e}_{j,\iota}(\varpi)
-\vec{e}_{j,\iota}(\tilde\varpi)\|\le 2{\rm dist}(E_j(\varpi),E_j(\tilde\varpi))
\end{align*}
by (\ref{dis-subsp}). Then, it follows from {Theorem \ref{thm-HD}} that
\begin{align}
\begin{split}
&\|P^{\pm 1}(\varpi)-P^{\pm 1}(\tilde \varpi)\|\le \|\bar x-\bar y\|^{\beta_E},
\quad \forall \varpi:=(\omega,\bar x),~ \tilde \varpi:=(\omega,\bar y)\in \bar {\bf \Omega},
\\
&\|P^{\pm 1}(\varpi)-{\rm id}\|
=\|P^{\pm 1}(\omega,\bar x)-P^{\pm 1}(\omega,0)\|
\le \delta_E(\omega),
\end{split}
\label{PPP}
\end{align}
where the constant $2$ is omitted without loss of generality.
}

Now, for every $\varpi=(\omega,\bar x)$ set
\begin{align}
\hat {\bf F}(\varpi,x)
:=P(\vartheta\varpi)F(\omega, P^{-1}(\varpi)x+\bar x)-P(\vartheta\varpi)F(\omega,\bar x),\quad x\in \mathbb{R}^d.
\label{def-hatF}
\end{align}
Then,
\begin{align*}
\hat {\bf F}(\vartheta^n\varpi,x)
&=P(\theta^{n+1}\omega, F(n+1,\omega,\bar x))F(\theta^{n}\omega, P^{-1}(\theta^n\omega, F(n,\omega,\bar x))x+ F(n,\omega,\bar x))
\\
&\quad-P(\theta^{n+1}\omega, F(n+1,\omega,\bar x))F(\theta^{n}\omega, F(n,\omega,\bar x)),
\end{align*}
which means that the problem of $F(\theta^n\omega,\cdot)$ near its orbit $F(n,\omega,\bar x)$ is converted to the problem of
$\hat {\bf F}(\varpi,x)$ near the origin $0$.
It is straightforward to compute that
\begin{align*}
\hat {\bf F}(\varpi,{{0}})={{0}},\qquad
D\hat {\bf F}(\varpi,{{0}})=P(\vartheta\varpi){\bf \Lambda}(\varpi)
P^{-1}(\varpi)=\bar{\bf \Lambda}(\varpi),
\end{align*}
and
\begin{align}
\begin{split}
\|D^2\hat {\bf F}(\varpi,x)\|\le 8M_{\epsilon}(\omega),
\quad \forall x\in \mathbb{R}^d,
\label{c-off-RS-hat}
\end{split}
\end{align}
due to (\ref{c-off-RS}). To show \eqref{c-off-RS-hat}, we also need to note that
\begin{align}
\|P^{\pm 1}(\varpi)\|\le 1+\delta_E(\omega)=
1+\delta(\omega)^{\beta_E} = 1+ 1/ \mathfrak{M}_\epsilon(\omega)^{\beta_E}\le 2
\label{P2}
\end{align}
by (\ref{PPP}) and (\ref{MM-p}).
For the sake of simplicity, we still use $M_{\epsilon}(\omega)$ to denote $8M_{\epsilon}(\omega)$.

Next, we state a result on  the normal form of
$\hat {\bf F}(\varpi,x)$, which will be proved in section 6.

\begin{prop}\label{prop-NF}
There is a $C^{2,\alpha}$ random diffeomorphism $N(\varpi,\cdot):\mathbb{R}^d\to \mathbb{R}^d$ transforming $\hat {\bf F}(\varpi,\cdot)$ to $\bar {\bf F}(\varpi,\cdot)$, i.e.,
\begin{align}
\bar {\bf F}(\varpi,x):=N(\vartheta\varpi,\cdot)\circ \hat {\bf F}(\varpi,\cdot) \circ N^{-1}(\varpi,x),
\label{bar-hat-0}
\end{align}
such that
\begin{align}
\bar {\bf F}(\varpi,{{0}})={{0}},\quad
D\bar {\bf F}(\varpi,{{0}})=\bar{\bf \Lambda}(\varpi),\quad
\partial^2_{x_u x_s}\bar {\bf F}(\varpi,{{0}})={\bf 0}~~~ \mbox{\rm (the zero mapping)},
\label{NF-1-0}
\end{align}
and
\begin{align}
\begin{split}
&N^{\pm 1}(\varpi,0)=0,\qquad DN^{\pm 1}(\varpi,0)={\rm id},\quad \|DN^{\pm 1}(\varpi,x)-{\rm id}\|\le \delta(\omega)/K_\epsilon(\omega),
\\
&\|D^2N^{\pm 1}(\varpi,x)\|\le \check M_\epsilon(\omega),\qquad
\|D^2N^{\pm 1}(\varpi,x)-D^2N^{\pm 1}(\varpi,y)\|\le \check L_{\epsilon, \delta}(\omega)\|x-y\|^\alpha
\end{split}
\label{estm-N}
\end{align}
for all $x,y\in\mathbb{R}^d$, where $\check M_{\epsilon}(\omega)>0$ is an $\epsilon$-slowly varying random variable depending on $K_\epsilon(\omega), M_\epsilon(\omega)$ and
$\check L_{\epsilon, \delta}(\omega)>0$ is an $\epsilon$-slowly varying random variable depending on $K_\epsilon(\omega)$, $M_\epsilon(\omega)$,
$\delta(\omega)$.

Moreover, $D^2N^{\pm 1}(\varpi,0)=D^2N^{\pm 1}(\omega,\bar x, 0)$ are measurable with respect to $\omega$ and are continuous with respect to $\bar x$ such that
\begin{align}
&\|D^2N^{\pm 1}(\omega,\bar x,0)-D^2N^{\pm 1}(\omega,\bar y,0)\|
\le C_\lambda K_\epsilon(\omega)^3\hat L_{\epsilon, \rho}(\omega)\|\bar x-\bar y\|^{\beta_N},
\label{estm-N-bar}
\end{align}
when $\|\bar x-\bar y\|< 1$, where
$
\beta_N:={\min\{\epsilon/\lambda_{\max},\alpha\}}
$
and $\hat L_{\epsilon, \rho}(\omega)>0$ is an $\epsilon$-slowly varying random variable depending on $K_\epsilon(\omega)$, $M_\epsilon(\omega)$ and $L_{\epsilon,\rho}(\omega)$.
\end{prop}

Remark that in the above definition of $\beta_N$, the choice of $\epsilon>0$ depends on the Lyapunov exponents only (see (\ref{def-ep})). Moreover, a key part of the above result is to obtain the third equality of (\ref{NF-1-0}), which will be used to estimate the nonlinear part of a cohomological equation deduced from (\ref{ch-eqn}), as will be seen in (\ref{ff00})-(\ref{EM-B1}) below.

Using (\ref{def-hatF}) and (\ref{bar-hat-0}), we have that
\begin{align}
\bar {\bf F}(\varpi,x)
={\Psi}(\vartheta\varpi, \cdot)\circ
F(\omega,\cdot)\circ {\Psi}^{-1}(\varpi, x),\quad \forall x\in\mathbb{R}^d,
\label{def-zzF}
\end{align}
where
$
{\Psi}(\varpi, x):=N(\varpi,\cdot)\circ P(\varpi)(x-\bar x).
$
It follows from
(\ref{c-off-RS}) and (\ref{estm-N}) that
\begin{align}
\begin{split}
&\|D\bar {\bf F}(\varpi,x)-\bar {\bf \Lambda}(\varpi)\|\le C_\lambda \delta(\omega),\qquad \|D^2\bar {\bf F}(\varpi,x)\|\le \bar M_\epsilon(\omega),
\\
&\|D^2\bar {\bf F}(\varpi,x)-D^2\bar {\bf F}(\varpi,y)\|\le \bar L_{\epsilon,\delta}(\omega)\|x-y\|^\alpha,\quad \forall x,y\in\mathbb{R}^d,
\end{split}
\label{con-barF}
\end{align}
and $D^2\bar {\bf F}(\varpi,\cdot)$ is $\epsilon$-slowly continuous when $\alpha=0$,
where $\bar M_{\epsilon}(\omega)>0$ is an $\epsilon$-slowly varying random variable depending on $K_\epsilon(\omega), M_\epsilon(\omega)$, and
$\bar L_{\epsilon,\delta}(\omega)>0$ is an $\epsilon$-slowly varying random variable depending on $K_\epsilon(\omega),  M_\epsilon(\omega), \delta(\omega)$. As we saw in section 3 (below (\ref{c-off-RS})), the third inequality of (\ref{con-barF}) holds for all H\"older exponents $\tilde \alpha\in [0,\alpha]$.

{\bf Step 2}. {\it Convert equation {\rm (\ref{ch-eqn})} to a new equation {\rm (\ref{InvDis-eqn-bar})} with respect to $\bar {\bf F}(\varpi,x)$, whose solution is represented by $\bar {\bf v}(\varpi,x)$, and establish the relation between $v(\omega,x)$ and $\bar {\bf v}(\varpi,x)$.}

For any $v(\omega,x)$ satisfying equation (\ref{ch-eqn}), we  define $\bar {\bf v}(\varpi, \cdot):\mathbb{R}^d\to \mathbb{R}^d$ by
\begin{align}
\begin{split}
\bar {\bf v}(\varpi, x)
:=D{\Psi}(\varpi,{\Psi}^{-1}(\varpi, x))\,v(\omega, {\Psi}^{-1}(\varpi, x)).
\end{split}
\label{vvv-1}
\end{align}
Then,  $\bar {\bf v}(\varpi, x)$ satisfies
\begin{align}
D\bar {\bf F}(\varpi, x)\bar {\bf v}(\varpi, x)
=\bar {\bf v}(\vartheta\varpi, \bar {\bf F}(\varpi, x)),\quad \forall x\in \mathbb{R}^d.
\label{InvDis-eqn-bar}
\end{align}
In fact, using (\ref{ch-eqn}) and (\ref{def-zzF}), we obtain that
\begin{align*}
&D\bar {\bf F}(\varpi, x)\bar {\bf v}(\varpi, x)
\\
&
=D\bar {\bf F}(\varpi, x)
\,
D{\Psi}(\varpi,{\Psi}^{-1}(\varpi, x))
\,
v(\omega, {\Psi}^{-1}(\varpi, x))
\\
&=D{\Psi}(\vartheta\varpi, F(\omega,\cdot)\circ {\Psi}^{-1}(\varpi, x))
\,
DF(\omega,{\Psi}^{-1}(\varpi, x))
\,
D{\Psi}^{-1}(\varpi, x)
\\
&\quad \cdot
D{\Psi}(\varpi,{\Psi}^{-1}(\varpi, x))
\,
v(\omega, {\Psi}^{-1}(\varpi, x))
\\
&=D{\Psi}(\vartheta\varpi, F(\omega,\cdot)\circ {\Psi}^{-1}(\varpi, x))
\,
DF(\omega,{\Psi}^{-1}(\varpi, x))
\,
v(\omega,{\Psi}^{-1}(\varpi, x))
\\
&=D{\Psi}(\vartheta\varpi, F(\omega,\cdot)\circ {\Psi}^{-1}(\varpi, x))
\,
v(\theta\omega, F(\omega,\cdot)\circ{\Psi}^{-1}(\varpi, x))
\\
&=D{\Psi}(\vartheta\varpi, {\Psi}^{-1}(\vartheta\varpi, \cdot)\circ \bar {\bf F}(\varpi, x))
\,
v(\theta\omega, {\Psi}^{-1}(\vartheta\varpi, \cdot)\circ \bar {\bf F}(\varpi, x))
\\
&=\bar {\bf v}(\vartheta\varpi, \bar {\bf F}(\varpi, x)).
\end{align*}
The converse is also true. Thus, the problem of finding the solution $v(\omega,x)$ of equation (\ref{ch-eqn}) is equivalent to finding the solution $\bar {\bf v}(\varpi, x)$ of equation (\ref{InvDis-eqn-bar}).

Let ${\mathcal W}_i(\varpi):={\mathcal W}_i(\omega,\bar x)$ be the leaf of an intermediate foliation of $F(\omega,x)$ (see (\ref{inter-foli})) passing through $\bar x$ for every $i=1,...,\tau$. Then,
\begin{align}
\bar {\mathcal N}_i(\varpi)
:=\Psi(\varpi,{\mathcal W}_i(\varpi)),
\quad \forall n\in \mathbb{Z},
\label{NNMM}
\end{align}
is an invariant manifold (i.e., the leaf passing through $0$) of $\bar {\bf F}(\varpi,x)$ corresponding to the subspace $X_i$.
In fact, using (\ref{def-zzF}), we have that
\begin{align*}
\bar {\bf F}(\varpi,\bar {\mathcal N}_i(\varpi))
&=\bar {\bf F}(\varpi,\cdot)\circ \Psi(\varpi,{\mathcal W}_i(\varpi))
=\Psi(\vartheta\varpi,\cdot)\circ F(\omega,{\mathcal W}_i(\varpi))
\\
&\subset\Psi(\vartheta\varpi,{\mathcal W}_i(\vartheta\varpi))=\bar {\mathcal N}_i(\vartheta\varpi).
\end{align*}

{\bf Step 3}. {\it Linearize $\bar {\bf F}(\varpi,x)$ along $\bar {\mathcal N}_i(\varpi)$ for $i=1,...,\tau$.}

By Corollary \ref{cor-PSU} (using (\ref{con-barF}) instead of (\ref{c-off-RS})), we see that $\bar {\mathcal N}_i(\varpi)$ is the graph of a mapping, i.e,
\begin{align}
\bar {\mathcal N}_i(\varpi)
=\{x_i+{\gamma}_i(\varpi,x_i):x_i\in X_i\},\quad i=1,...,\tau,
\label{def-Nh}
\end{align}
where ${\gamma}_i(\varpi,\cdot): X_i\to X_1\oplus\cdots\oplus X_{i-1}\oplus X_{i+1}\oplus \cdots \oplus X_{p}$ is a
$C^{1,{\beta}}$ mapping such that
\begin{align}
\begin{split}
&{\gamma}_i(\varpi,{{0}})={{0}},\quad D{\gamma}_i(\varpi,{{0}})={\bf 0},\quad\|D{\gamma}_i(\varpi,x_i)\|\le C_\lambda K_\epsilon(\omega)
C_\epsilon(\omega)\delta(\omega),
\\
&\|D{\gamma}_i(\varpi,x_i)-D{\gamma}_i(\varpi,y_i)\|\le L_{\epsilon, {\gamma}}(\omega)\|x_i-y_i\|^{\beta},\quad \forall x_i,y_i\in X_i,
\label{hhOO}
\end{split}
\end{align}
with an $\epsilon$-slowly varying tempered random variable $L_{\epsilon, \gamma}(\omega)>0$.
Using the transformation
$$
\varphi_i(\varpi, x):=x-{\gamma}_i(\varpi,x_i),
$$
we straighten up the invariant manifold $\bar {\mathcal N}_i({\varpi})$. Then,
we have a $C^{1,{\beta}}$ mapping $\bar {\bf F}_{i}(\varpi,\cdot):X_i\to X_i$ by
\begin{align}
&\bar {\bf F}_{i}(\varpi,x_i)
:=\varphi_i(\vartheta\varpi, \cdot)\circ
\bar {\bf F}(\varpi, \cdot)\circ \varphi^{-1}_i(\varpi, x_i),
\label{F*12}
\end{align}
which satisfies that
\begin{align}
\begin{split}
&\bar {\bf F}_{i}(\varpi,{{0}})={{0}},\quad D\bar {\bf F}_{i}(\varpi,{{0}})=\bar{\bf \Lambda}_{i}(\varpi),
\quad \|D\bar {\bf F}_{i}(\varpi,x_i)-\bar{\bf \Lambda}_{i}(\varpi)\|
\le \delta^*(\omega),
\\
&\|D\bar {\bf F}_{i}^{\pm 1}(\varpi,x_i)-D\bar {\bf F}_{i}^{\pm 1}(\varpi,y_i)\|
\le L_\epsilon^*(\omega)\|x_i-y_i\|^{\beta}, \qquad \forall x_i,y_i\in X_i,
\end{split}
\label{FFOO-RS}
\end{align}
by (\ref{con-barF}) and (\ref{hhOO}). Here, $\bar{\bf \Lambda}_{i}(\varpi)$ is given in (\ref{blo-diag}), $\delta^*(\omega):=C_\lambda K_\epsilon(\omega)^2
C_\epsilon(\omega)\delta(\omega)>0$ and
$
L_\epsilon^*(\omega)>0
$
depends on $L_{\epsilon, {\gamma}}(\omega)$ and $\bar M_\epsilon(\omega)$.

In order to solve equation (\ref{InvDis-eqn-bar}), we will use the Lyapunov norms to overcome the nonuniformity (see Lemma \ref{lm-LO} below together with the remark given at the end of section 7). Since the Lyapunov norms is applied only to the linear system $\bar{\bf \Lambda}_{ i}(\varpi)$, we need to linearize $\bar {\bf F}_{i}(\varpi,x_i)$.

The following lemma states that $\bar {\bf F}_{i}(\varpi,x_i)$ is $C^{1,\beta}$ with $\beta\in (0,1]$ given in (\ref{def-bbb}) conjugate to its linear part, which will be proved in section 7.

\begin{lm}\label{lm-lin}
There is a $C^1$ diffeomorphism $\varphi_{i,*}(\varpi,\cdot):X_i\to X_i$ such that
\begin{align}
\bar{\bf \Lambda}_{i}(\varpi)\varphi_{i,*}(\varpi,x_i)
=\varphi_{i,*}(\vartheta\varpi,\cdot)\circ \bar {\bf F}_{i}(\varpi,x_i),\quad \forall i=1,...,\tau,
\label{Con-eqn*}
\end{align}
and
\begin{align}
\begin{split}
D\varphi_{i,*}^{\pm 1}(\varpi,{{0}})={\rm id}_i,
\quad\|D\varphi_{i,*}^{\pm 1}(\varpi,x_i)-{\rm id}_i\|\le C^*_{4\epsilon}(\omega) \|x_i\|^{\beta},
\quad \forall x_i\in X_i,
\end{split}
\label{phi-Bd}
\end{align}
where $C^*_{4\epsilon}(\omega)>0$ is a $4\epsilon$-slowly varying tempered random variable and ${\rm id}_i$ the identity mapping in $X_i$.
\end{lm}

{\bf Step 4}. {\it Decompose equation {\rm (\ref{InvDis-eqn-bar})}, obtained in step {\rm 2}, along every $\bar {\mathcal N}_i(\varpi)$ for $i=1,...,\tau$ by using the conjugacy obtained in step {\rm 3}, and solve the decomposed equation {\rm (\ref{Distb-RRS})}.}

Recall that $\varphi_i(\varpi, \cdot):\mathbb{R}^d\to \mathbb{R}^d$
straightens up $\bar {\mathcal N}_i({\varpi})$ which is tangent to $X_i$ at ${{0}}$. Let
\begin{align}
\bar \varphi_{i}(\varpi,\cdot):=\varphi_{i,*}(\varpi,\cdot)\circ \varphi_{i}(\varpi,\cdot)|_{\bar {\mathcal N}_i({\varpi})}: \bar {\mathcal N}_i({\varpi})\to X_i.
\label{def-p*}
\end{align}
Then, we have
\begin{align*}
\bar \varphi_i(\vartheta\varpi, \cdot)\circ
\bar {\bf F}(\varpi, \cdot)\circ \bar \varphi^{-1}_i(\varpi, x_i)
=\bar{\bf \Lambda}_{i}(\varpi)
\end{align*}
by (\ref{F*12}) and (\ref{Con-eqn*}).
Replacing $x$ in (\ref{InvDis-eqn-bar}) with $\bar \varphi^{-1}_i(\varpi, x_i)$, we obtain an equation
\begin{align}
D\bar {\bf F}(\varpi,\bar \varphi_i^{-1}(\varpi,x_i))\,\bar {\bf v}_{i}(\varpi,x_i)
=\bar {\bf v}_{i}(\vartheta\varpi,\bar{\bf \Lambda}_{i}(\varpi)x_i)
\label{Distb-RRS}
\end{align}
for $i=1,...,\tau$,
where $\bar {\bf v}_{i}(\varpi,\cdot):X_i\to \mathbb{R}^d$ is defined by
\begin{align}
\bar {\bf v}_{i}(\varpi,x_i):=\bar {\bf v}(\varpi, \bar \varphi_i^{-1}(\varpi,x_i)).
\label{vvxb}
\end{align}
For $\kappa={\tau}+1,...,{p}$ and $\iota= 1,...,d_\kappa$, where $d_\kappa$ is the dimension of $X_\kappa$, let $\vec{e}_{\kappa,\iota}:=(0,...,1,...,0)^T$ be the standard  unit vectors in $X_\kappa$, i.e., the $(d_1+\cdots+d_{\kappa-1}+\iota)$-th component is 1 and others
are 0.

The following proposition is on the existence of solutions of equation (\ref{Distb-RRS}), which will be proved in sections 7 and 8 for the cases of $\alpha\in (0,1]$ and $\alpha=0$ respectively.

\begin{prop}\label{prop-ID-RS}
Given a $\kappa={\tau}+1,...,{p}$ and $d_{\kappa}$ standard unit vectors
$
\vec{e}_{\kappa,\iota}\in X_\kappa$ {\rm (}$\iota= 1,...,d_\kappa${\rm )}, for every $\varpi\in \bar {\bf \Omega}$, there exist sets $\{{\bar {\bf v}_{i}^{\kappa,\iota}}(\vartheta^n\varpi,x_i)\}_{n\in\mathbb{Z}}$
{\rm (}$i=1,...,\tau${\rm )} of the solutions of equation {\rm (\ref{Distb-RRS})}
such that
$
{\bar {\bf v}_{i}^{\kappa,\iota}}(\vartheta^n\varpi,{{0}})=\bar{\bf \Lambda}_{\kappa}(n,\varpi)\,\vec{e}_{\kappa,\iota}\in X_\kappa,
$
\begin{align}
\frac{\|{\bar {\bf v}_{i}^{\kappa,\iota}}(\varpi,x_i)-{\bar {\bf v}_{i}^{\kappa,\iota}}(\varpi,{{0}})\|}{\|x_i\|}
\to 0 \qquad {\rm as} ~\|x_i\|\to 0,
\label{v-1-sig-XX}
\end{align}
and
\begin{align}
&\sup_{n\in\mathbb{Z}}\sup_{x_i\in X_i\backslash\{{{0}}\}}\bigg\{e^{-\lambda_\kappa n-12\epsilon |n|} \frac{\|{\bar {\bf v}_{i}^{\kappa,\iota}}(\vartheta^n\varpi,x_i)-{\bar {\bf v}_{i}^{\kappa,\iota}}(\vartheta^n\varpi,{{0}})\|}{\|x_i\|^{\varsigma}}
\bigg\}
\le M_{v}(\omega),
\label{v-sig-RS}
\end{align}
where $M_{v}(\omega)>0$ is a tempered random variable and $\varsigma\in (0,1)$ is a constant satisfying
\begin{align}
\begin{split}
0<\varsigma
<\frac{1}{10}\min\Bigg\{\frac{\lambda_{\tau}}{\lambda_{\max}}, \frac{-\lambda_{\tau+1}}{\lambda_{\max}}, \min_{\mbox{\tiny
$\begin{array}{c}
{\kappa\in\{\tau+1,...,p\},j\in\{1,...,p\}}
\\
{\rm such ~ that}~ \lambda_\kappa-\lambda_j\ne 0
\end{array}$}}
\frac{|\lambda_\kappa-\lambda_j|}
{\lambda_1}\Bigg\}
\end{split}
\label{sg12}
\end{align}
with
$\lambda_{\rm max}:=\max\{2\lambda_1,-2\lambda_{p}\}>0$.
\end{prop}

{\bf Step 5}. {\it Use $\{\bar {\bf v}_{i}^{\kappa,\iota}(\theta^n\varpi,x_i)\}_{n\in\mathbb{Z}}$
to get a set $\{v_{\bar x, i}^{\kappa,\iota}(\theta^n\omega,x)\}_{n\in\mathbb{Z}}$ of solutions
of equation {\rm (\ref{ch-eqn})} for every $\omega\in \Omega$ on the leaves $\{{\mathcal W}_i(\theta^n\omega,F(n,\omega, \bar x))\}_{n\in\mathbb{Z}}$ through {\rm (\ref{vvv-1})}. Then, show that $v_{\bar x, i}^{\kappa,\iota}(\omega,x)$ is differentiable at $\bar x$ along ${\mathcal W}_i(\omega,\bar x)$
such that the differential is {\rm (}H\"older{\rm )} continuous with respect to $\bar x$}.

For each point $x\in {\mathcal W}_i(\varpi)$ for $i=1,...,\tau$, by the invariance of foliation, we have that
$F(n,\omega,x)\in {\mathcal W}_i(\vartheta^n\varpi)$. Then, using (\ref{NNMM}) and  the fact that
$\bar \varphi_i(\vartheta^n\varpi,\cdot)$ straightens up $\bar {\mathcal N}_i(\vartheta^n\varpi)$, we get
\begin{align}
\bar \varphi_i(\vartheta^n\varpi,\cdot)\circ
\Psi(\vartheta^n\varpi,F(n,\omega,x))\in X_i,
\qquad \forall n\in\mathbb{Z},
\label{phiX}
\end{align}
It implies that
\[
{\bar {\bf v}_{i}^{\kappa,\iota}}(\vartheta^n\varpi, \bar \varphi_i(\vartheta^n\varpi,\cdot)\circ{\Psi}(\vartheta^n\varpi, F(n,\omega,x))
\]
is well-defined. Then, by (\ref{vvv-1}) (implying
$v(\omega, x)
=\{D{\Psi}(\varpi,x)\}^{-1}\bar {\bf v}(\varpi, {\Psi}(\varpi, x))$) and by (\ref{vvxb}),
we obtain that
\begin{align}
&{v_{\bar x, i}^{\kappa,\iota}}(\theta^n\omega,F(n,\omega,x))
\nonumber\\
&:=\{D{\Psi}(\vartheta^n\varpi,F(n,\omega,x))\}^{-1}\,{\bar {\bf v}_{i}^{\kappa,\iota}}(\vartheta^n\varpi,
\bar \varphi_i(\vartheta^n\varpi,\cdot)\circ{\Psi}(\vartheta^n\varpi, F(n,\omega,x))),\quad \forall n\in \mathbb{Z},
\label{def-vkl}
\end{align}
satisfies (\ref{ch-eqn}) at $F(n,\omega,x)$, i.e.,
\begin{align}
DF(\theta^n\omega,F(n,\omega,x))\,{v_{\bar x, i}^{\kappa,\iota}}(\theta^n\omega,F(n,\omega,x))
={v_{\bar x, i}^{\kappa,\iota}}(\theta^{n+1}\omega,F(n+1,\omega,x)).
\label{ch-eqn-Fn}
\end{align}
Since
$
{\Psi}^{-1}(\varpi, x)=P^{-1}(\varpi)\circ N^{-1}(\varpi,x)+\bar x
$
due to the definition of ${\Psi}$ given below (\ref{def-zzF}), we have
$$
{\Psi}^{-1}(\vartheta^n\varpi, x)=P^{-1}(\vartheta^n\varpi)\circ N^{-1}(\vartheta^n\varpi,x)+F(n,\omega,\bar x)
$$
by replacing $\varpi=(\omega,\bar x)$ and $\bar x$ with $\vartheta^n\varpi=(\theta^n\omega,F(n,\omega,\bar x))$
and $F(n,\omega,\bar x)$ respectively.
Therefore,
\begin{align*}
\{D{\Psi}(\vartheta^n\varpi,F(n,\omega,x))\}^{-1}
&=D{\Psi}^{-1}(\vartheta^n\varpi,{\Psi}(\vartheta^n\varpi,F(n,\omega,x)))
\\
&=P^{-1}(\vartheta^n\varpi)\,DN^{-1}(\vartheta^n\varpi,{\Psi}
(\vartheta^n\varpi,F(n,\omega,x))),
\end{align*}
which together with (\ref{def-vkl}) yields
\begin{align}
&{v_{\bar x, i}^{\kappa,\iota}}(\theta^n\omega, F(n,\omega,x))
\nonumber\\
&=P^{-1}(\vartheta^n\varpi)\,DN^{-1}(\vartheta^n\varpi,{\Psi}(\vartheta^n\varpi,F(n,\omega,x)))\,{\bar {\bf v}_{i}^{\kappa,\iota}}(\vartheta^n\varpi, \bar \varphi_i(\vartheta^n\varpi,\cdot)\circ{\Psi}(\vartheta^n\varpi, F(n,\omega,x))).
\label{def-vkl-D}
\end{align}

Next, we show the differentiability of ${v_{\bar x, i}^{\kappa,\iota}}(\omega, x)$ at $x=\bar x$ along
${\mathcal W}_i(\omega, \bar x)={\mathcal W}_i(\varpi)$. In view of (\ref{phiX}) with $n=0$, we  let
\begin{align*}
x_i:=\bar \varphi_i(\varpi,\cdot)\circ\Psi(\varpi,x)\in X_i.
\end{align*}
Thus, $x\mapsto x_i(x)$ is a $C^{1}$ diffeomorphism from ${\mathcal W}_i(\varpi)$ to $X_i$ with the inverse $x_i\mapsto x(x_i):=\Psi^{-1}(\varpi,\cdot)\circ\bar \varphi_i^{-1}(\varpi,x_i)$ such that $x({{0}})=\bar x$.
Then, in view of (\ref{def-vkl-D}) with $n=0$ and with $x$ being replaced by $x(x_i):=\Psi^{-1}(\varpi,\cdot)\circ\bar \varphi_i^{-1}(\varpi,x_i)$, we have that
\begin{align*}
{v_{\bar x, i}^{\kappa,\iota}}(\omega, x(x_i))
=P^{-1}(\varpi)\,DN^{-1}(\varpi,\bar\varphi_i^{-1}(\varpi,x_i))\,{\bar {\bf v}_{i}^{\kappa,\iota}}(\varpi,x_i),\quad \forall x_i\in X_i.
\end{align*}
Therefore,
\begin{align}
\frac{dv_{\bar x, i}^{\kappa,\iota}(\omega, x(x_i))}{dx_i}\Big|_{x_i={{0}}}
&=P^{-1}(\varpi)\{D^2N^{-1}(\varpi, {{0}}){\rm id}_{i}\}\,\vec{e}_{\kappa,\iota}
\nonumber\\
&=P^{-1}(\omega,\bar x)\{D^2N^{-1}(\omega,\bar x, {{0}}){\rm id}_{i}\}\,\vec{e}_{\kappa,\iota}\in {\mathcal L}(X_i,\mathbb{R}^d)
\label{Dv-WU}
\end{align}
since $D\bar\varphi_i^{-1}(\varpi,{{0}})={\rm id}_{i}$ (see ({\ref{hhOO}}) and (\ref{phi-Bd})), ${\bar {\bf v}_{i}^{\kappa,\iota}}(\varpi,{{0}})
=\vec{e}_{\kappa,\iota}$ (see {Proposition \ref{prop-ID-RS}}) and $D{\bar {\bf v}_{i}^{\kappa,\iota}}(\varpi,{{0}})=0$ (see (\ref{v-1-sig-XX})).
Then (\ref{Dv-WU}) determines the differential (\ref{dif-tial}) of $v_{\bar x, i}^{\kappa,\iota}(\omega, \cdot):{\mathcal W}_i(\omega,\bar x)\to \mathbb{R}^d$ along ${\mathcal W}_i(\omega,\bar x)$ at $\bar x$ (\cite[p.16]{War}).

More precisely, since
$$
\frac{dx({{0}})}{dx_i}
=\frac{\Psi^{-1}(\varpi,\cdot)\circ\bar \varphi_i^{-1}(\varpi,x_i)}{dx_i}\Big|_{x_i=0}
=P^{-1}(\varpi)\partial_{x_i}N^{-1}(\varpi,0)=P^{-1}(\omega, \bar x){\rm id}_i\in {\mathcal L}(X_i,\mathbb{R}^d)
$$
by (\ref{estm-N}), we have that
$\{P^{-1}(\omega, \bar x)\vec{e}_{i,\hat\iota}\}_{\hat\iota=1,...,d_i}$ is a basis of $T_{\bar x}
{\mathcal W}_i(\omega,\bar x)$, where $\{\vec{e}_{i,\hat\iota}\}_{\hat\iota=1,...,d_i}$
is the standard basis of $X_i$. Notice that every
$P^{-1}(\omega, \bar x)\vec{e}_{i,\hat\iota}$ corresponds (one-to-one) to a linear mapping
\begin{align}\label{def-deriv}
h_{\mathcal W_i}\mapsto \frac{\partial (h_{\mathcal W_i}\circ x(x_i))}{\partial x_{i,\hat \iota}}\Big|_{x_i=0}=Dh_{\mathcal W_i}(\bar x)P^{-1}(\omega, \bar x){\rm id}_{i,\hat\iota}\in {\mathcal L}(X_{i,\hat\iota},\mathbb{R})
\end{align}
for all smooth functions $h_{\mathcal W_i}: {\mathcal W}_i(\omega, \bar x)\to \mathbb{R}$, where ${\rm id}_{i,\hat\iota}$ denotes the identity in the 1-dimensional space $X_{i,\hat\iota}$. {This observation enables us to use $P^{-1}(\omega, \bar x)\vec{e}_{i,\hat\iota}$ here to denote the tangent vector for convenience, although (\ref{def-deriv}) is exactly the definition of tangent vector given in
\cite[Definition 1.19]{War}.}
Then, we have that
for any given ${\bf a}_{\omega,\bar x}\in T_{\bar x}
{\mathcal W}_i(\omega,\bar x)~(\cong X_i)$,
there are
$
c_1,\cdots,c_{d_i}\in\mathbb{R}
$
such that
\begin{align}
{\bf a}_{\omega,\bar x}=\sum_{\hat \iota=1,...,d_i}c_{\hat \iota}P^{-1}(\omega,\bar x)\vec{e}_{i,\hat\iota}.
\label{aaa-1}
\end{align}
This enables us to compute the differential of $v_{\bar x, i}^{\kappa,\iota}(\omega, \cdot)$
along ${\mathcal W}_i(\omega,\bar x)$ at $\bar x$ by (\ref{Dv-WU}), that is,
\begin{align}
dv_{\bar x, i}^{\kappa,\iota}(\omega, \bar x)\,{\bf a}_{\omega,\bar x}
&=\sum_{\hat \iota=1,...,d_i}c_{\hat \iota}P^{-1}(\omega, \bar x)\{D^2N^{-1}(\omega, \bar x, {{0}})\vec{e}_{i,\hat\iota}\}\,\vec{e}_{\kappa,\iota}
\nonumber\\
&=\sum_{\hat \iota=1,...,d_i}c_{\hat \iota}P^{-1}(\omega,\bar x)D^2N^{-1}(\omega,\bar x, {{0}})(\vec{e}_{i,\hat\iota},\vec{e}_{\kappa,\iota})
\label{diff-v}
\end{align}
(see \cite[p.17, 1.23(a)]{War}). Note that $dv_{\bar x, i}^{\kappa,\iota}(\omega, \bar x):T_{\bar x}
{\mathcal W}_i(\omega,\bar x)\to \mathbb{R}^d  ~ (=T_{v_{\bar x, i}^{\kappa,\iota}(\omega, \bar x)}\mathbb{R}^d)$.

Regarding the H\"older continuity of $dv_{\bar x, i}^{\kappa,\iota}(\omega, \bar x)$ with respect to $\bar x$,
we define the bundle
$$
{\mathcal E}:=\bigcup_{\bar x\in \mathbb{R}^d}T_{\bar x}
{\mathcal W}_i(\omega,\bar x),
$$
which is a subbundle of $T\mathbb{R}^d(=\mathbb{R}^{2d})$ (therefore, in ${\mathcal E}$, one can use the same norm as the one in $\mathbb{R}^{2d}$).
In order to trivialize ${\mathcal E}$, let
$
{\bf a}:=\sum_{\hat \iota=1,...,d_i}c_{\hat \iota}\vec{e}_{i,\hat\iota}\in X_i
$
and define a mapping $\hbar:\Omega\times \mathbb{R}^d\times X_i\to {\mathcal E}$ as
\begin{align*}
\hbar(\omega, \bar x, {\bf a})= {\bf a}_{\omega,\bar x}.
\end{align*}
Since $\mathbb{R}^d\times X_i$ is a subspace of $\mathbb{R}^{2d}$, we regard $\hbar(\omega,\cdot)$ as a mapping from a subset of $\mathbb{R}^{2d}$ into $\mathbb{R}^{2d}$, which is clearly one-to-one. Thus, in view of (\ref{aaa-1}),
we see that $\hbar(\omega,\cdot)$ and $\hbar^{-1}(\omega,\cdot)$ are measurable with respect to $\omega$, $\beta_E$-H\"older continuous with respect to $\bar x$ with a uniformly bounded H\"older coefficient with respect to $\bar x$ and $i$ by (\ref{PPP}). Moreover,
we have that $\hbar(\omega,\bar x, {\bf a})$ and $\hbar^{-1}(\omega,\bar x, {\bf a})$ are linear with respect to ${\bf a}$ and that $\hbar^{-1}(\omega,\cdot)$ trivializes $T_{\mathbb{R}^d}
{\mathcal W}_i(\omega,\bar x)$.

In what follows, we rewrite $\hbar(\omega, \bar x, {\bf a})$ as $\hbar(\omega, \bar x) {\bf a}$ for convenience. Then,
for any ${\bf a}\in X_i$ such that $\|{\bf a}\|=1$, using (\ref{diff-v}), we have that
\begin{align*}
&\|dv_{\bar x, i}^{\kappa,\iota}(\omega, \bar x)\circ \hbar(\omega, \bar x){\bf a}
-dv_{\bar y, i}^{\kappa,\iota}(\omega, \bar y)\circ \hbar(\omega, \bar y) {\bf a}\|
=\|dv_{\bar x, i}^{\kappa,\iota}(\omega, \bar x){\bf a}_{\omega,\bar x}
-dv_{\bar y, i}^{\kappa,\iota}(\omega, \bar y){\bf a}_{\omega,\bar y}\|
\\
&\le \sum_{\hat \iota=1,...,d_i}\big\|c_{\hat \iota}P^{-1}(\omega,\bar x)D^2N^{-1}(\omega,\bar x, {{0}})(\vec{e}_{i,\hat\iota},\vec{e}_{\kappa,\iota})
-c_{\hat \iota}P^{-1}(\omega,\bar y)D^2N^{-1}(\omega,\bar y, {{0}})(\vec{e}_{i,\hat\iota},\vec{e}_{\kappa,\iota})
\big\|
\\
&\le 4\check M_\epsilon(\omega)\|\bar x-\bar y\|^{\beta_E}+4C_\lambda K_\epsilon(\omega)^3\hat L_{\epsilon,\rho}(\omega)\|\bar x-\bar y\|^{\beta_N}
\\
&\le 4\{\check M_\epsilon(\omega)+C_\lambda K_\epsilon(\omega)^3\hat L_{\epsilon,\rho}(\omega)\}\|\bar x-\bar y\|^{\beta_v}
\end{align*}
when $\|\bar x-\bar y\|<1$. Here, (\ref{PPP}) and (\ref{estm-N})-(\ref{estm-N-bar}) are also used to obtain the above estimates, and
$$
\beta_v:=\min\{\beta_E,\beta_N\}
=\min\{\beta_E,\epsilon/\lambda_{\max},\alpha\}\ge 0.
$$
Hence,
\begin{align*}
&\|dv_{\bar x, i}^{\kappa,\iota}(\omega, \bar x)\circ \hbar(\omega, \bar x)
-dv_{\bar y, i}^{\kappa,\iota}(\omega, \bar y)\circ \hbar(\omega, \bar y)\|\le 4\{\check M_\epsilon(\omega)+C_\lambda K_\epsilon(\omega)^3\hat L_{\epsilon,\rho}(\omega)\}\|\bar x-\bar y\|^{\beta_v},
\end{align*}
which means that
the differential $dv_{\bar x, i}^{\kappa,\iota}(\omega, \bar x)$ is $\beta_v$-H\"older continuous with respect to $\bar x$ with a uniformly bounded H\"older coefficient with respect to $\bar x$ and $i$.

{\bf Step 6}. {\it Show that when $x\in {\mathcal W}_i(\omega, \bar x)$ is close to $\bar x$, $v_{\bar x, i}^{\kappa,\iota}(\omega,x)$
 belong to $E_s(\omega, x)$ for all $\kappa={\tau}+1,...,{p}$ and $\iota=1,...,d_\kappa$, and form a basis of $E_s(\omega, x)$.
}

Next, we show that for every $x\in {\mathcal W}_i(\omega, \bar x)$ near $\bar x$ the set $\{{v_{\bar x, i}^{\kappa,\iota}}(\omega, x)\}_{\kappa={\tau}+1,...,{p},~\iota=1,...,d_\kappa}$ is a basis of $E_s(\omega, x)$. Notice that (\ref{def-vkl-D}) gives
\begin{align}
{v_{\bar x, i}^{\kappa,\iota}}(\theta^n\omega, F(n,\omega,\bar x))
&=P^{-1}(\vartheta^n\varpi)\,{\bar {\bf v}_{i}^{\kappa,\iota}}(\vartheta^n\varpi, {{0}})=P^{-1}(\vartheta^n\varpi)\bar{\bf \Lambda}_{\kappa}(n,\varpi)\,\vec{e}_{\kappa,\iota}
\label{def-vkl-Dbar}
\end{align}
since $DN^{-1}(\vartheta^n\varpi,{{0}})={\rm id}$ and ${\Psi}(\vartheta^n\varpi,F(n,\omega,\bar x))={{0}}$. Then, in view of (\ref{def-vkl-D}) and (\ref{def-vkl-Dbar}),
we get
\begin{align}
&e^{-\lambda_\kappa n-16\epsilon|n|}\|{v_{\bar x, i}^{\kappa,\iota}}(\theta^n
\omega,F(n,\omega,x))-{v_{\bar x, i}^{\kappa,\iota}}(\theta^n\omega,F(n,\omega,\bar x))\|
\nonumber\\
&\le e^{-\lambda_\kappa n-16\epsilon|n|}\|P^{-1}(\vartheta^n\varpi)\|
\nonumber\\
&\quad\cdot
\|DN^{-1}(\vartheta^n\varpi,{\Psi}(\vartheta^n\varpi,F(n,\omega,x)))\,{\bar {\bf v}_{i}^{\kappa,\iota}}(\vartheta^n\varpi, \bar \varphi_i(\vartheta^n\varpi,\cdot)\circ{\Psi}(\vartheta^n\varpi, F(n,\omega,x)))-{\bar {\bf v}_{i}^{\kappa,\iota}}(\vartheta^n\varpi, {{0}})\|
\nonumber\\
&\le 2e^{-\lambda_\kappa n-16\epsilon|n|}\Big\{\|DN^{-1}(\vartheta^n\varpi,{\Psi}(\vartheta^n\varpi,F(n,\omega,x)))\,{\bar {\bf v}_{i}^{\kappa,\iota}}(\vartheta^n\varpi, \bar \varphi_i(\vartheta^n\varpi,\cdot)\circ{\Psi}(\vartheta^n\varpi, F(n,\omega,x)))
\nonumber\\
&\quad-DN^{-1}(\vartheta^n\varpi,{\Psi}(\vartheta^n\varpi,F(n,\omega,x)))\,{\bar {\bf v}_{i}^{\kappa,\iota}}(\vartheta^n\varpi, {{0}})\|
\nonumber\\
&\quad+\|DN^{-1}(\vartheta^n\varpi,{\Psi}(\vartheta^n\varpi,F(n,\omega,x)))\,{\bar {\bf v}_{i}^{\kappa,\iota}}(\vartheta^n\varpi, {{0}})
-{\bar {\bf v}_{i}^{\kappa,\iota}}(\vartheta^n\varpi, {{0}})\|\Big\}
\nonumber\\
&\le 4e^{-\lambda_\kappa n-16\epsilon|n|}\Big\{\|{\bar {\bf v}_{i}^{\kappa,\iota}}(\vartheta^n\varpi, \bar \varphi_i(\vartheta^n\varpi,\cdot)\circ{\Psi}(\vartheta^n\varpi, F(n,\omega,x)))-
{\bar {\bf v}_{i}^{\kappa,\iota}}(\vartheta^n\varpi, {{0}})\|
\nonumber\\
&\quad+\|DN^{-1}(\vartheta^n\varpi,{\Psi}(\vartheta^n\varpi,F(n,\omega,x)))
-{\rm id}\|^\varsigma\|\bar{\bf \Lambda}_{\kappa}(n,\varpi)\|\Big\},
 \label{v-difference}
\end{align}
where $\varsigma\in (0,1)$ is given in (\ref{sg12}). Here, we also used $\|P^{-1}(\vartheta^n\varpi)\|\le 2$ by (\ref{P2}) and
\[
\|DN^{-1}(\varpi,x)-{\rm id}\|\le \delta(\omega)/K_\epsilon(\omega)
=1/(K_\epsilon(\omega)
\mathfrak{M}_\epsilon(\omega))\le 1
\]
from (\ref{MM-p}), which implies that
\[
\|DN^{-1}(\vartheta^n\varpi,x)\|\le 2,\quad \|DN^{-1}(\vartheta^n\varpi,x)
-{\rm id}\|\le \|DN^{-1}(\vartheta^n\varpi,
x)
-{\rm id}\|^\varsigma,\quad \forall x\in \mathbb{R}^d.
\]

{In order to further estimate (\ref{v-difference}), for all $z_i\in\bar {\mathcal N}_i(\vartheta^n\varpi)$, we notice that
\begin{align*}
\|\varphi_{i}(\vartheta^n\varpi,z_i)-z_i\|
&=\|\gamma_i(\vartheta^n\varpi,z_i)\|\le \sup_{c\in(0,1)}\|D\gamma_{i}(\vartheta^n\varpi,cz_i)\|
\,\|z_i\|
\\
&\le C_\lambda K_\epsilon(\theta^n\omega)
C_\epsilon(\theta^n\omega)\delta(\theta^n\omega)\,\|z_i\|\le
C_\lambda K_\epsilon(\theta^n\omega)
C_\epsilon(\theta^n\omega)/\mathfrak{M}_\epsilon(\theta^n\omega)\|z_i\|
\\
&\le \|z_i\|
\end{align*}
by the first and the third inequalities of (\ref{hhOO}) and by (\ref{MM-p}), implying $\|\varphi_{i}(\vartheta^n\varpi,z_i)\|\le 2\|z_i\|$.
Moreover, the second inequality of (\ref{phi-Bd}) gives
$$
\|\varphi_{i,*}(\vartheta^n\varpi,x_i)
-x_i\|\le \sup_{c\in(0,1)}\|D\varphi_{i,*}(\vartheta^n\varpi,cx_i)
-{\rm id}_i\|\,\|x_i\|\le C^*_{4\epsilon}(\theta^n\omega)\|x_i\|^{1+\beta}
$$
for all $x_i\in X_i$, implying $\|\varphi_{i,*}(\vartheta^n\varpi,x_i)\|\le \|x_i\|+C^*_{4\epsilon}(\theta^n\omega)\|x_i\|^{1+\beta}$.
Thus, combining (\ref{def-p*}) with the above facts,
we get
\begin{align}
\|\bar \varphi_{i}(\vartheta^n\varpi,z_i)\|
&= \|\varphi_{i,*}(\vartheta^n\varpi,\cdot)\circ \varphi_{i}(\vartheta^n\varpi,z_i)\|
\le \|\varphi_{i}(\vartheta^n\varpi,z_i)\|+C^*_{4\epsilon}(\theta^n\omega)\|\varphi_{i}(\vartheta^n\varpi,z_i)\|^{1+\beta}
\nonumber\\
&\le 2\|z_i\|+4e^{4\epsilon|n|}C^*_{4\epsilon}(\omega)
\|z_i\|^{1+\beta}.
\label{phi**}
\end{align}
}
We also remind that
$
{\Psi}(\vartheta^n\varpi, x):=N(\vartheta^n\varpi,\cdot)\circ P(\vartheta^n\varpi)(x-F(n,\omega,\bar x)),
$
which implies that
\begin{align}\label{es-PPSS}
\|{\Psi}(\vartheta^n\varpi, F(n,\omega,x))\|
&=
\|N(\vartheta^n\varpi,\cdot)\circ P(\vartheta^n\varpi)(F(n,\omega,x)-F(n,\omega,\bar x))\|
\nonumber\\
&\le \sup_{\xi\in\mathbb{R}^d}\|DN(\vartheta^n\varpi,\xi)\|
\,\|P(\vartheta^n\varpi)\|\,\|F(n,\omega,x)-F(n,\omega,\bar x)\|
\nonumber\\
&\le 4\|F(n,\omega,x)-F(n,\omega,\bar x)\|.
\end{align}

Now, using (\ref{v-sig-RS}), {(\ref{estm-N})}, (\ref{ED-tmp}) and (\ref{phi**})-(\ref{es-PPSS}), we continue to estimate \eqref{v-difference} and obtain that
\begin{align}
&e^{-\lambda_\kappa n-16\epsilon|n|}\|{v_{\bar x, i}^{\kappa,\iota}}(\theta^n
\omega,F(n,\omega,x))-{v_{\bar x, i}^{\kappa,\iota}}(\theta^n\omega,F(n,\omega,\bar x))\|
\nonumber\\
&\le 4e^{-4\epsilon|n|}M_{v}(\omega)\|\bar \varphi_i(\vartheta^n\varpi,\cdot)\circ{\Psi}(\vartheta^n\varpi, F(n,\omega,x))\|^\varsigma
\nonumber\\
&\quad+{ 4\sup_{\xi\in\mathbb{R}^d}
\|D^2N^{-1}(\vartheta^n\varpi,\xi)\|^\varsigma
\|{\Psi}(\vartheta^n\varpi,F(n,\omega,x))\|^\varsigma
e^{-\lambda_\kappa n-16\epsilon|n|}\|\bar{\bf \Lambda}_{\kappa}(n,\varpi)\|}
\nonumber\\
&\le 4M_{v}(\omega)
\Big\{2\|{\Psi}(\vartheta^n\varpi, F(n,\omega,x))\|
+C^*_{4\epsilon}(\omega)
\|{\Psi}(\vartheta^n\varpi, F(n,\omega,x))\|^{1+\beta}\Big\}^\varsigma
\nonumber\\
&\quad+4K_\epsilon(\omega)\check M_\epsilon(\omega)\|{\Psi}(\vartheta^n\varpi, F(n,\omega,x))\|^\varsigma
\nonumber\\
&\le 4M_{v}(\omega)
\Big\{8\|F(n,\omega,x)-F(n,\omega,\bar x)\|
+16C^*_{4\epsilon}(\omega)
\|F(n,\omega,x)-F(n,\omega,\bar x)\|^{1+\beta}\Big\}^\varsigma
\nonumber\\
&\quad+ 16K_\epsilon(\omega)\check M_\epsilon(\omega)\|F(n,\omega,x)-F(n,\omega,\bar x)\|^\varsigma
\nonumber\\
&
{\le
\tilde M_{v}(\omega)
\Big\{\sup_{\xi\in\mathbb{R}^d}\|DF(n,\omega,\xi)\|\,\|x-\bar x\|
+
\sup_{\xi\in\mathbb{R}^d}\|DF(n,\omega,\xi)\|^{1+\beta}\|x-\bar x\|^{1+\beta}\Big\}^\varsigma}
\nonumber\\
&\quad{ + \tilde M_{v}(\omega)\sup_{\xi\in\mathbb{R}^d}\|DF(n,\omega,\xi)\|^\varsigma
\|x-\bar x\|^\varsigma}
\nonumber\\
&\le 3\tilde M_{v}(\omega)\sup_{\xi\in\mathbb{R}^d}\|DF(n,\omega,\xi)\|^{2\varsigma}\|x-\bar x\|^\varsigma,\qquad \forall n\ge n(\omega)>0,
\label{long-est}
\end{align}
with a large integer $n(\omega)>0$ due to the facts that $\sup_{\xi\in\mathbb{R}^d}\|DF(n,\omega,\xi)\|>1$
for all $n\ge n(\omega)$ (see (\ref{ED-ge})) and $\|x-\bar x\|<1$,
where
$
\tilde M_{v}(\omega):=64M_{v}(\omega)C^*_{4\epsilon}(\omega)+16K_\epsilon(\omega)\check M_\epsilon(\omega)>0.
$

On the other hand, we see from (\ref{ch-eqn-Fn}) that, for all ${n}\ge 0$,
\begin{align*}
DF({n},\omega,x){v_{\bar x, i}^{\kappa,\iota}}(\omega,x)
&=DF(\theta^{{n}-1}\omega,F({n}-1,\omega,x))\cdots
DF(\omega,x)
{v_{\bar x, i}^{\kappa,\iota}}(\omega,x)
\\
&={v_{\bar x, i}^{\kappa,\iota}}(\theta^{n}\omega,F({n},\omega,x)).
\end{align*}
Hence, using (\ref{long-est}), (\ref{def-vkl-Dbar}),(\ref{DF-max}), and (\ref{ED-tmp}), we obtain
\begin{align}
&\sup_{{n}\ge n(\omega)}\{e^{-\lambda_\kappa {n}-(16\epsilon+2\lambda_{\max}\varsigma){{n}}}\|DF({n},\omega,x){v_{\bar x, i}^{\kappa,\iota}}(\omega, x)\|\}
\nonumber\\
&=\sup_{{n}\ge n(\omega)} \{e^{-\lambda_\kappa {n}-(16\epsilon+2\lambda_{\max}\varsigma){{n}}}\|{v_{\bar x, i}^{\kappa,\iota}}(\theta^{n}\omega,F({n},\omega,x))\|\}
\nonumber\\
&\le
\sup_{{n}\ge n(\omega)} \{e^{-\lambda_\kappa {n}-16\epsilon{{n}}-2\lambda_{\rm max}{\varsigma}{{n}}}\|{v_{\bar x, i}^{\kappa,\iota}}(\theta^{n}\omega,F({n},\omega,x))-{v_{\bar x, i}^{\kappa,\iota}}(\theta^{n}\omega,F({n},\omega,\bar x))\|\}
\nonumber\\
&\quad+\sup_{{n}\ge n(\omega)} \{e^{-\lambda_\kappa {n}-\epsilon{{n}}}\|{v_{\bar x, i}^{\kappa,\iota}}(\theta^{n}\omega,F({n},\omega,\bar x))\|\}
\nonumber\\
&\le 3\tilde M_{v}(\omega)\sup_{{n}\ge 0} \{e^{-\lambda_{\rm max} {{n}}}\sup_{\xi\in\mathbb{R}^d}\|DF({n},\omega,\xi)\|\}^{2\varsigma}\|x-\bar x\|^\varsigma
\nonumber\\
&\quad+2\sup_{{n}\ge 0} \{e^{-\lambda_\kappa {n}-\epsilon{{n}}}\|\bar{\bf \Lambda}_{\kappa}({n},\varpi)\|\}
<\infty
\label{vvxx-1}
\end{align}
{for all $x\in {\mathcal W}_i(\omega,\bar x)$ such that $\|x-\bar x\|<1$.
It means that
\begin{align}
\sup_{{n}\ge n(\omega)}\|DF({n},\omega,x){v_{\bar x, i}^{\kappa,\iota}}(\omega, x)\|<\infty
\label{DF-finit}
\end{align}
since $\lambda_\kappa+16\epsilon+2\lambda_{\max}\varsigma<0$ for all $\kappa=\tau+1,...,p$ by (\ref{def-ep}) and (\ref{sg12}), which enables us to claim that
\begin{align}
{v_{\bar x, i}^{\kappa,\iota}}(\omega,x)\in E_s(\omega,x),
\quad \forall x\in {\mathcal W}_i(\omega,\bar x)~\mbox{near $\bar x$},
\label{vbE}
\end{align}
for all $\kappa={\tau}+1,...,{p}$ and $\iota=1,...,d_\kappa$.
In fact, let $\Pi_u(\omega,x)$ and $\Pi_s(\omega,x)$ be the projections associated with $\mathbb{R}^d=E_u(\omega,x)\oplus E_s(\omega,x)$.
Assume that ${v_{\bar x, i}^{\kappa,\iota}}(\omega,x)\notin E_s(\omega,x)$, which implies that
$$
\|\Pi_u(\omega,x)v_{\bar x, i}^{\kappa,\iota}(\omega,x)\|\ne 0.
$$
Then, similarly to (\ref{unq-s}) we have
\begin{align*}
&\|DF({n},\omega,x){v_{\bar x, i}^{\kappa,\iota}}(\omega,x)\|
\\
&\ge (1/K_\epsilon(\omega))e^{(\varrho_0-\epsilon) {n}}\|\Pi_u(\omega,x)v_{\bar x, i}^{\kappa,\iota}(\omega,x)\|
-K_\epsilon(\omega)e^{-\varrho_0  {n}}\|\Pi_s(\omega,x)v_{\bar x, i}^{\kappa,\iota}(\omega,x)\|
\\
&\to +\infty, \qquad {\rm as}\quad {n}\to +\infty,
\end{align*}
with a constant $\varrho_0>0$ given in (\ref{ED}), which contradicts (\ref{DF-finit}). Hence, (\ref{vbE}) is proved.
}

Moreover, we have that for $x\in {\mathcal W}_i(\omega,\bar x)$ near $\bar x$ the vectors ${v_{\bar x, i}^{\kappa,\iota}}(\omega,x)$ are close to ${v_{\bar x, i}^{\kappa,\iota}}(\omega,\bar x)
=P^{-1}(\omega, \bar x)\vec{e}_{\kappa,\iota}$ (by (\ref{def-vkl-Dbar}) with $n=0$) for all $\kappa={\tau}+1,...,{p}$ and $\iota=1,...,d_\kappa$.
It means that these ${v_{\bar x, i}^{\kappa,\iota}}(\omega,x)$ are also close to $\vec{e}_{\kappa,\iota}$ by (\ref{PPP}) and therefore
the set $\{{v_{\bar x, i}^{\kappa,\iota}}(\omega,x)\}_{\kappa={\tau}+1,...,{p},~\iota=1,...,d_\kappa}$ is a basis (depending on $\bar x$ and $i$) of $E_s(\omega,x)$ for $x\in {\mathcal W}_i(\omega,\bar x)$ near $\bar x$.
This completes the proof of Theorem \ref{Thm-decop}. \qquad $\Box$

{Remark that, in the above proof, Proposition \ref{prop-ID-RS} is an important result. In fact, (\ref{v-1-sig-XX}) plays a key role in proving the differentiability of $v_{\bar x, i}^{\kappa,\iota}(\omega, \cdot)$ along ${\mathcal W}_i(\omega,\bar x)$ at $\bar x$ in step 5 (see (\ref{Dv-WU})), while
(\ref{v-sig-RS}) plays a key role in proving that
$v_{\bar x, i}^{\kappa,\iota}(\omega,x)\in E_s(\omega, x)$ for all $x\in {\mathcal W}_i(\omega, \bar x)$ near $\bar x$ in step 6 (see the first inequality of (\ref{long-est})).
}

\section{Random normal form}
\setcounter{equation}{0}

In this section, we prove Proposition \ref{prop-NF}. A general theory of normal forms for random dynamical systems was established in \cite{Arn98}. However, to prove $C^{1,\beta_\alpha}$ linearization, we need a result on  $C^{2,\alpha}$ normal form with specific estimates, which was not studied in \cite{Arn98}.

{\bf Proof of {Proposition \ref{prop-NF}}}.
Let $x_j\in X_j$ for $j=1,...,p$ and ${\mathcal H}_d^2$ be the vector space of all homogeneous
polynomials of the second order with $d$ variables in $\mathbb{R}^d$. Set
$
{\mathcal H}_d^2=\bigoplus_{i,\kappa,j\in \{1,...,{p}\}}{\mathcal S}^2_{i,\kappa,j}
$
where
\begin{align*}
{\mathcal S}^2_{i,\kappa,j}:=\big\{c_{i,\kappa,j}x_ix_\kappa\in {\mathcal H}_d^2: c_{i,\kappa,j}\in{\mathcal L}(X_i, X_\kappa; X_j)\big\}.
\end{align*}
Here ${\mathcal L}(X_i, X_\kappa; X_j)$ denotes the set of all
bilinear mappings from $X_i\times X_\kappa$ into $X_j$.
Define
$$
N(\varpi,x):=x+a(\varpi)x^2,\quad \forall n\in \mathbb{Z},
$$
where $a(\varpi)\in {\mathcal L}(\mathbb{R}^d, \mathbb{R}^d; \mathbb{R}^d)$ is a bilinear mapping from $\mathbb{R}^d\times \mathbb{R}^d$ into $\mathbb{R}^d$. As in \cite[pp.30-32]{LZZ}, one can use the standard
argument as the proof for the Inverse Function Theorem
to show that if $a(\varpi)$ is tempered, then $N(\varpi,\cdot):\mathbb{R}^d\to \mathbb{R}^d$ is invertible and its inverse
$$
N^{-1}(\varpi,x)=x-a(\varpi)x^2+o(\|x\|^2)
$$
exists in a small tempered ball $\tilde V(\omega)\subset \mathbb{R}^d$. Thus, we compute that
\begin{align*}
&N(\vartheta\varpi,\cdot)\circ \hat {\bf F}(\varpi,\cdot)\circ N^{-1}(\varpi,x)
\\
&=\bar{\bf \Lambda}(\varpi) x+\Big\{a(\vartheta\varpi) (\bar{\bf \Lambda}(\varpi) x)^2-\bar{\bf \Lambda}(\varpi) (a(\varpi)x^2)
+\frac{1}{2}D^2\hat {\bf F}(\varpi,{{0}})x^2\Big\}
+o(\|x\|^2).
\end{align*}
On the other hand, let
$\Pi_{{\mathcal S}^2_{i,\kappa,j}}:{\mathcal H}_d^2\to {\mathcal S}^2_{i,\kappa,j}$ be the projection associated with the above direct sum and let
$a_{i,\kappa,j}(\varpi)\in {\mathcal L}(X_i,X_\kappa; X_j)$ be a bilinear mapping defined by
\begin{align*}
a_{i,\kappa,j}(\varpi)x_ix_\kappa:=\Pi_{{\mathcal S}^2_{i,\kappa,j}}\{a(\varpi)x^2\}.
\end{align*}
Noting that our goal is to eliminate the term $x_ux_s$, i.e., the terms $x_ix_\kappa$ for all $i=1,...,\tau$ and $\kappa=\tau+1,...,{p}$, we compute that
\begin{align*}
&\Pi_{{\mathcal S}^2_{i,\kappa,j}}\Big\{a({\vartheta\varpi}) (\bar{\bf \Lambda}(\varpi) x)^2-\bar{\bf \Lambda}(\varpi) (a(\varpi)x^2)
+\frac{1}{2}D^2\hat {\bf F}(\varpi,{{0}})x^2\Big\}
\\
&=a_{i,\kappa,j}(\vartheta\varpi) \{\bar{\bf \Lambda}_{i}(\varpi) x_i\}\{\bar{\bf \Lambda}_{\kappa}(\varpi) x_\kappa\}-\bar{\bf \Lambda}_{j}(\varpi)\{a_{i,\kappa,j}(\varpi)x_ix_\kappa\}
-c_{i,\kappa,j}(\varpi)x_ix_\kappa
\end{align*}
where $c_{i,\kappa,j}(\varpi):=-\frac{1}{2}\partial^2_{x_i x_\kappa}(\pi_j\hat {\bf F})(\varpi,{{0}})$.
Thus, in what follows, we consider the equation
\begin{align}
&a_{i,\kappa,j}(\vartheta\varpi) \{\bar{\bf \Lambda}_{i}(\varpi) x_i\}\{\bar{\bf \Lambda}_{\kappa}(\varpi) x_\kappa\}-\bar{\bf \Lambda}_{j}(\varpi)\{a_{i,\kappa,j}(\varpi)x_ix_\kappa\}
=c_{i,\kappa,j}(\varpi)x_ix_\kappa
\label{Eqn-NF}
\end{align}
with unknown $a_{i,\kappa,j}(\varpi)$ for $i=1,...,\tau$, $\kappa=\tau+1,...,{p}$ and $j=1,...,{p}$.

In order to solve equation (\ref{Eqn-NF}), for an  $a_{i,\kappa,j}(\varpi)\in {\mathcal L}(X_i, X_\kappa; X_j)$, we define its norm by
\begin{align*}
\|a_{i,\kappa,j}(\varpi)\|:=
\sup_
{\mbox{\tiny
$\begin{array}{c}
x_i\in X_i\backslash\{{{0}}\},
\\
x_\kappa\in X_\kappa\backslash\{{{0}}\}
\end{array}$}}
\frac{\|a_{i,\kappa,j}(\varpi)\,x_i x_\kappa\|}
{\|x_i\|\|x_\kappa\|}.
\end{align*}
Notice that, by the non-resonant condition (\ref{EC-RS}), we have two cases that $\lambda_i+\lambda_\kappa>\lambda_j$ and $\lambda_i+\lambda_\kappa<\lambda_j$.
In the case that $\lambda_i+\lambda_\kappa>\lambda_j$, consider a sequence $(\hat a^{(k)}_{i,\kappa,j}(\varpi))_{k\in \mathbb{N}}$ where $\hat a^{(0)}_{i,\kappa,j}(\varpi):={\bf 0}$, defined recursively by
\begin{align*}
&\hat a_{i,\kappa,j}^{(k+1)}(\vartheta\varpi) x_ix_\kappa
\\
&:=c_{i,\kappa,j}(\varpi)\{\bar{\bf \Lambda}_{i}^{-1}(\varpi) x_i\}\{\bar{\bf \Lambda}_{\bar x,\kappa}^{-1}(\varpi) x_\kappa\}+\bar{\bf \Lambda}_{j}(\varpi)\hat a_{i,\kappa,j}^{(k)}(\varpi)\{\bar{\bf \Lambda}_{i}^{-1}(\varpi) x_i\}\{\bar{\bf \Lambda}_{\bar x,\kappa}^{-1}(\varpi) x_\kappa\}
\end{align*}
via (\ref{Eqn-NF}). We
compute that for all $k\in \mathbb{N}$
\begin{align}
\hat a^{(k)}_{i,\kappa,j}(\varpi)\,x_i x_\kappa
= \sum_{m=1}^{k}\bar{\bf \Lambda}_{j}(m-1,\vartheta^{-m+1}\varpi)
 c_{i,\kappa,j}(\vartheta^{-m}\varpi)\{\bar{\bf \Lambda}_{i}(-m,\varpi)x_i\}\{\bar{\bf \Lambda}_{\kappa}(-m,\varpi)x_\kappa\}.
\label{YD}
\end{align}
For its estimation, we note that for all $n\in\mathbb{Z}$
\begin{align}
\sup_
   {\mbox{\tiny
   $\begin{array}{c}x_i\in X_i\backslash\{{{0}}\},
   \\
   x_\kappa\in X_\kappa\backslash\{{{0}}\}
   \end{array}$}}\frac{\|c_{i,\kappa,j}(\varpi)x_ix_\kappa\|}{\|x_i\|\|x_\kappa\|}
= \|c_{i,\kappa,j}(\varpi)\|\le \|\partial^2_{x_u x_s} \hat {\bf F}(\varpi,{{0}})\|\le  M_\epsilon(\omega)
\label{eesstt}
\end{align}
by (\ref{c-off-RS-hat}).
Then, for all $m\ge 1$
\begin{align}
&
\frac{\|\bar{\bf \Lambda}_{j}(m-1,\vartheta^{-m+1}\varpi)
 c_{i,\kappa,j}(\vartheta^{-m}\varpi)\{\bar{\bf \Lambda}_{i}(-m,\varpi)x_i\}\{\bar{\bf \Lambda}_{\kappa}(-m,\varpi)x_\kappa\}\|}
  {\|x_i\|\,\|x_\kappa\|}
\nonumber\\
&\le K_\epsilon(\vartheta^{-m+1}\varpi)e^{(\lambda_j+\epsilon)(m-1)}
\nonumber\\
&\quad\cdot\frac{\|c_{i,\kappa,j}(\vartheta^{-m}\varpi)\{\bar{\bf \Lambda}_{i}(-m,\varpi)x_i\}\{\bar{\bf \Lambda}_{\kappa}(-m,\varpi)x_\kappa\}\|
  }{\|\bar{\bf \Lambda}_{i}(-m,\varpi)x_i\|\,\|\bar{\bf \Lambda}_{\kappa}(-m,\varpi)x_\kappa\|}
  \,\frac{\|
    \bar{\bf \Lambda}_{i}(-m,\varpi)x_i\|\|\bar{\bf \Lambda}_{\kappa}(-m,\varpi)x_\kappa\|}
    {\|x_i\|\,\|x_\kappa\|}
\nonumber\\
&\le e^{-(\lambda_j+2\epsilon)}e^{\epsilon m} K_\epsilon(\omega)^3e^{(\lambda_j+\epsilon)m}
   e^{-(\lambda_i-\epsilon)m}e^{-(\lambda_\kappa-\epsilon)m}
   e^{2\epsilon m}
   \sup_
   {\mbox{\tiny
   $\begin{array}{c}x_i\in X_i\backslash\{{{0}}\},
   \\
   x_\kappa\in X_\kappa\backslash\{{{0}}\}
   \end{array}$}}
   e^{-2\epsilon m}
   \frac{\|c_{i,\kappa,j}(\vartheta^{-m}\varpi)x_ix_\kappa\|
   }{\|x_i\|\,\|x_\kappa\|}
\nonumber\\
&=e^{-(\lambda_j+2\epsilon)}K_\epsilon(\omega)^3 M_\epsilon(\omega)
   e^{(\lambda_j-\lambda_i-\lambda_\kappa+6\epsilon)m}
\label{RS-CNG}
\end{align}
due to (\ref{ED-tmp}).
It implies that the sequence $(\hat  a^{(k)}_{i,\kappa,j}(\varpi))_{k\in \mathbb{N}}
=(\hat  a^{(k)}_{i,\kappa,j}(\omega,\bar x))_{k\in \mathbb{N}}$ is uniformly convergent with respect to $\bar x$ since
$
\lambda_j-\lambda_i-\lambda_\kappa+6\epsilon<0
$
in the present case of $\lambda_i+\lambda_\kappa> \lambda_j$ by (\ref{def-ep}). Namely,
$
 \hat a_{i,\kappa,j}(\varpi):=\lim_{k\to +\infty} \hat a^{(k)}_{i,\kappa,j}(\varpi)
$
is a solution of equation (\ref{Eqn-NF}), which is measurable with respect to $\omega$ and continuous with respect to $\bar x$ such that
\begin{align}
\|\hat a_{i,\kappa,j}(\varpi)\|\le \sum_{m=1}^{\infty}e^{-(\lambda_j+2\epsilon)} K_\epsilon(\omega)^3 M_\epsilon(\omega)
   e^{(\lambda_j-\lambda_i-\lambda_\kappa+5\epsilon)m} \le C_\lambda  K_\epsilon(\omega)^3  M_\epsilon(\omega)
\label{D2N-1}
\end{align}
with a constant $C_\lambda>0$.

Concerning the H\"older continuity with respect to $\bar x$, for $\varpi:=(\omega,\bar x)$ and $\tilde \varpi:=(\omega,\bar y)$, we verify from (\ref{def-hatF}), (\ref{c-off-RS}), (\ref{PPP}) and (\ref{DF-max}) that
\begin{align}
&\|D^2\hat {\bf F}(\vartheta^n\varpi,0)-D^2\hat {\bf F}(\vartheta^n\tilde \varpi,0)\|
\nonumber\\
&\le \|P(\vartheta^{n+1}\varpi)D^2F(\theta^n\omega, F(n,\omega,\bar x))P^{-2}(\vartheta^n\varpi)-P( \vartheta^{n+1}\tilde\varpi)D^2F(\theta^n\omega, F(n,\omega,\bar y))P^{-2}(\vartheta^n\tilde \varpi)\|
\nonumber\\
&\le 4M_\epsilon(\theta^n\omega)
\|P(\vartheta^{n+1}\varpi)-P(\vartheta^{n+1}\tilde\varpi)\|
^{\epsilon/(\beta_E\lambda_{\max})}
\nonumber\\
&\quad+8\|D^2F(\theta^n\omega, F(n,\omega,\bar x))-D^2F(\theta^n\omega, F(n,\omega,\bar y))\|
\nonumber\\
&\quad+8M_\epsilon(\theta^n\omega)
\|P^{-1}(\vartheta^n\varpi)-P^{-1}
(\vartheta^n\tilde\varpi)\|
^{\epsilon/(\beta_E\lambda_{\max})}
\nonumber\\
&\le 4M_\epsilon(\theta^n\omega)\|F(n+1,\omega,\bar x)-F(n+1,\omega,\bar y)\|^{\epsilon/\lambda_{\max}}
\nonumber\\
&\quad+8L_{\epsilon,\rho}(\theta^n\omega)
\|F(n,\omega,\bar x)-F(n,\omega,\bar y)\|^{\min\{\epsilon/\lambda_{\max},\alpha\}}
+8M_\epsilon(\theta^n\omega)\|F(n,\omega,\bar x)-F(n,\omega,\bar y)\|^{\epsilon/\lambda_{\max}}
\nonumber\\
&\le 4M_\epsilon(\theta^n\omega)
K^{\epsilon/\lambda_{\max}}_\epsilon(\omega)
e^{\epsilon|n+1|}\|\bar x-\bar y\|^{\epsilon/\lambda_{\max}}
\nonumber\\
&\quad+8\{L_{\epsilon,\rho}(\theta^n\omega)
+M_\epsilon(\theta^n\omega)\}K^{\epsilon/\lambda_{\max}}_\epsilon(\omega)
e^{\epsilon|n|}\|\bar x-\bar y\|^{\min\{\epsilon/\lambda_{\max},\alpha\}}
\nonumber\\
&\le e^{2\epsilon|n|}\hat L_{\epsilon, \rho}(\omega)\|\bar x-\bar y\|^{\beta_N}
\label{DFxy}
\end{align}
for all $\bar x,\bar y\in \mathbb{R}^d$ such that $\|\bar x-\bar y\|< 1$, where
$$
\hat L_{\epsilon, \rho}(\omega)
:=\{4 e^{\epsilon}M_\epsilon(\omega)
+8L_{\epsilon,\rho}(\omega)
+8M_\epsilon(\omega)\}
K_\epsilon(\omega)>0,
\qquad \beta_N:={\min\{\epsilon/\lambda_{\max},\alpha\}}\ge 0.
$$
Notice that in the above formula,
it is reasonable to add small exponents $\epsilon/(\beta_E\lambda_{\max})\in (0,1)$ and ${\min\{\epsilon/\lambda_{\max},\alpha\}}\in (0,\alpha]$ because
$\|P^{\pm}(\varpi)-P^{\pm}(\tilde\varpi)\|\le 1$ by a similar argument to (\ref{P2}).
Thus, by (\ref{YD}), (\ref{RS-CNG}) and (\ref{DFxy}), we have
\begin{align}
&\|\hat a_{i,\kappa,j}(\varpi)-\hat a_{i,\kappa,j}(\tilde \varpi)\|
\nonumber\\
&\le C_\lambda K_\epsilon(\omega)^3
\sup_{m\ge 1}\sup_
   {\mbox{\tiny
   $\begin{array}{c}x_i\in X_i\backslash\{{{0}}\},
   \\
   x_\kappa\in X_\kappa\backslash\{{{0}}\}
   \end{array}$}}
   e^{-2\epsilon m}
   \frac{\|\{c_{i,\kappa,j}(\vartheta^{-m}\varpi)
   -c_{i,\kappa,j}(\vartheta^{-m}\tilde\varpi)\}x_ix_\kappa\|
   }{\|x_i\|\,\|x_\kappa\|}
\nonumber\\
&\le C_\lambda K_\epsilon(\omega)^3\sup_{m\ge 1}
e^{-2\epsilon m}\|\partial^2_{x_u x_s} \hat {\bf F}(\vartheta^{-m}\varpi,{{0}})-\partial^2_{x_u x_s} \hat {\bf F}_{\bar y}(\vartheta^{-m}\tilde\varpi,{{0}})\|
\nonumber\\
&\le C_\lambda K_\epsilon(\omega)^3\hat L_{\epsilon, \rho}(\omega)\|\bar x-\bar y\|^{\beta_N}
\label{D2N-22}
\end{align}
for all $\bar x,\bar y\in \mathbb{R}^d$ such that $\|\bar x-\bar y\|< 1$.

Next, in the case that $\lambda_i+\lambda_\kappa< \lambda_j$, the same conclusion can be obtained by considering a sequence $(\check a^{(k)}_{i,\kappa,j}(\varpi))_{k\in \mathbb{N}}$ where $\check a^{(0)}_{i,\kappa,j}(\varpi):={\bf 0}$, defined recursively by
\begin{align*}
\check a_{i,\kappa,j}^{(k+1)}(\varpi)x_ix_\kappa
:=\bar{\bf \Lambda}_{j}^{-1}(\varpi)\check a_{i,\kappa,j}^{(k)}(\vartheta\varpi) (\bar{\bf \Lambda}_{i}(\varpi) x_i)(\bar{\bf \Lambda}_{\kappa}(\varpi) x_\kappa)-\bar{\bf \Lambda}_{j}^{-1}(\varpi)c_{i,\kappa,j}(\varpi)x_ix_\kappa
\end{align*}
via (\ref{Eqn-NF}). We compute that for all $k\in \mathbb{N}$
\begin{align}
&\check a^{(k)}_{i,\kappa,j}(\varpi)\,x_i x_\kappa
= \sum_{m=1}^{k}\bar{\bf \Lambda}_{j}(-m,\vartheta^m\varpi)
c_{i,\kappa,j}(\vartheta^{m-1}\varpi)\{\bar{\bf \Lambda}_{i}(m-1,\varpi)x_i\}\{\bar{\bf \Lambda}_{\kappa}(m-1,\varpi)x_\kappa\}.
\label{YD-2}
\end{align}
A similar argument to (\ref{RS-CNG}) yields
\begin{align*}
&\frac{\|\bar{\bf \Lambda}_{j}(-m,\vartheta^m\varpi)
c_{i,\kappa,j}(\vartheta^{m-1}\varpi)\{\bar{\bf \Lambda}_{i}(m-1,\varpi)x_i\}\{\bar{\bf \Lambda}_{\kappa}(m-1,\varpi)x_\kappa\}\|}
{\|x_i\|\,\|x_\kappa\|}
\nonumber\\
&\le K_\epsilon(\theta^m\omega)e^{-(\lambda_j-\epsilon)m}
\nonumber\\
&\quad\cdot\frac{\|c_{i,\kappa,j}(\vartheta^{m-1}\varpi)\{\bar{\bf \Lambda}_{i}(m-1,\varpi))x_i\}\{\bar{\bf \Lambda}_{\kappa}(m-1,\varpi)x_\kappa\}\|
  }{\|\bar{\bf \Lambda}_{i}(m-1,\varpi)x_i\|\,\|\bar{\bf \Lambda}_{\kappa}(m-1,\varpi)x_\kappa\|}
  \, \frac{\|\bar{\bf \Lambda}_{i}(m-1,\varpi)x_i\|\,\|\bar{\bf \Lambda}_{\kappa}(m-1,\varpi)x_\kappa\|}
    {\|x_i\|\,\|x_\kappa\|}
\nonumber\\
&\le e^{-(\lambda_j-\epsilon)}e^{\epsilon m}K_\epsilon(\omega)^3 e^{-(\lambda_j-\epsilon)(m-1)}
   e^{(\lambda_i+\epsilon)(m-1)}e^{(\lambda_\kappa+\epsilon)(m-1)}
   e^{2\epsilon(m-1)}
   \nonumber\\
   &\quad\cdot\sup_
   {\mbox{\tiny
   $\begin{array}{c}x_i\in X_i\backslash\{{{0}}\},
   \\
   x_\kappa\in X_\kappa\backslash\{{{0}}\}
   \end{array}$}}e^{-2\epsilon(m-1)}\frac{\|c_{i,\kappa,j}(\vartheta^{m-1}\varpi)x_ix_\kappa\|
   }{\|x_i\|\,\|x_\kappa\|}
\nonumber\\
&\le e^{-(\lambda_j-\epsilon)} K_\epsilon(\omega)^3 M_\epsilon(\omega)
   e^{-(\lambda_j-\lambda_i-\lambda_\kappa-6\epsilon)(m-1)}.
\end{align*}
It implies that the sequence $(\check a^{(k)}_{i,\kappa,j}(\varpi))_{k\in \mathbb{N}}=(\check a^{(k)}_{i,\kappa,j}(\omega,\bar x))_{k\in \mathbb{N}}$ is uniformly convergent with respect to $\bar x$ since
$
\lambda_j-\lambda_i-\lambda_\kappa-6\epsilon>0
$
in the present case of $\lambda_i+\lambda_\kappa< \lambda_j$ by (\ref{def-ep}). Namely,
$
 \check a_{i,\kappa,j}(\varpi):=\lim_{k\to +\infty} \check a^{(k)}_{i,\kappa,j}(\varpi)
$
is a solution of equation (\ref{Eqn-NF}), which is measurable with respect to $\omega$ and continuous with respect to $\bar x$ such that
\begin{align}
\|\check a_{i,\kappa,j}(\varpi)\|\le C_\lambda K_\epsilon(\omega)^3  M_\epsilon(\omega)
\label{ab-2}
\end{align}
and
\begin{align}
\|\check a_{i,\kappa,j}(\varpi)-\check a_{i,\kappa,j}(\tilde \varpi)\|
\le C_\lambda K_\epsilon(\omega)^3\hat L_{\epsilon, \rho}(\omega)\|\bar x-\bar y\|^{\beta_N}
\label{ab-33}
\end{align}
for all $\bar x,\bar y\in \mathbb{R}^d$ such that $\|\bar x-\bar y\|< 1$ by similar arguments to (\ref{D2N-1}) and (\ref{D2N-22}).

Recall that $\Pi_{{\mathcal S}^2_{i,\kappa,j}}$ is the projection from ${\mathcal H}_d^2$ to the subspace ${\mathcal S}^2_{i,\kappa,j}$. Now, we define $a(\varpi)\in {\mathcal L}(X_i, X_\kappa; X_j)$ by
\begin{align*}
\Pi_{{\mathcal S}^2_{i,\kappa,j}}\{a(\varpi)x^2\}:=\left\{
  \begin{array}{lll}
   \hat a_{i,\kappa,j}(\varpi)x_ix_\kappa& ~~~{\rm if}~\lambda_i+\lambda_\kappa>\lambda_j,
   \\
   \check a_{i,\kappa,j}(\varpi)x_ix_\kappa& ~~~{\rm if}~\lambda_i+\lambda_\kappa<\lambda_j,
  \end{array}
\right.
\end{align*}
for all $(i,\kappa,j)\in \{1,...,\tau\}\times\{\tau+1,...,{p}\}\times\{1,...,{p}\}$, and
$
\Pi_{{\mathcal S}^2_{i,\kappa,j}}\{a(\varpi)x^2\}:={{0}}
$
for all other $(i,\kappa,j)$. Then, we know that
$a(\varpi)=a(\omega,\bar x)$ is measurable with respect to $\omega$ and continuous with respect to $\bar x$ such that
\begin{align}
\|a(\varpi)\|\le C_\lambda K_\epsilon(\omega)^3  M_\epsilon(\omega),\qquad
\|a(\varpi)-a(\tilde \varpi)\|
\le C_\lambda K_\epsilon(\omega)^3\hat L_{\epsilon, \rho}(\omega)\|\bar x-\bar y\|^{\beta_N}
\label{Est-a}
\end{align}
by (\ref{D2N-1}), (\ref{D2N-22}) and (\ref{ab-2})-(\ref{ab-33}). This enables us to define
a global random transformation (still denoted by $N(\varpi,\cdot):\mathbb{R}^d\to \mathbb{R}^d$) by
\begin{align*}
N(\varpi,x):=
\begin{cases}
x+u(\omega,x)a(\varpi)x^2,&\quad \text{for all}\; x\in \tilde V(\omega)\\
x, &\quad \text{for all}\; x\in \mathbb{R}^d\backslash \tilde V(\omega),
\end{cases}
\end{align*}
where $\tilde V(\omega)\subset \mathbb{R}^d$ is a neighborhood with a small diameter $\tilde \rho(\omega):=\delta(\omega)\{(3C_u+1)C_\lambda K_\epsilon(\omega)^4  M_\epsilon(\omega)\}^{-1}>0$ and
$u$ is given in (\ref{cutoff-w}) with $\tilde V(\omega)$ in place of $V(\omega)$. Using a similar argument to Lemma \ref{Extension} together with the first inequality of (\ref{Est-a}), we compute that
\begin{align*}
\begin{split}
&\|DN(\varpi,x)-{\rm id}\|\le \tilde \delta(\omega),\qquad \|D^2N(\varpi,x)\|\le (3C_u+1)C_\lambda K_\epsilon(\omega)^3  M_\epsilon(\omega),
\\
&\|D^2N(\varpi,x)-D^2N(\varpi,y)\|\le 7C_uC_\lambda K_\epsilon(\omega)^3  M_\epsilon(\omega)/\tilde \rho(\omega)\|x-y\|^\alpha, \quad \forall x,y\in \mathbb{R}^d,
\end{split}
\end{align*}
where $\tilde \delta(\omega):=(3C_u+1)C_\lambda K_\epsilon(\omega)^3  M_\epsilon(\omega)\tilde \rho(\omega)=\delta(\omega)/K_\epsilon(\omega)$.
Applying the Inverse Function Theorem, we see that
the inverse of $N(\varpi,\cdot)$ has the same smoothness as $N(\varpi,\cdot)$ itself, which proves (\ref{estm-N}).

It is clear that $N^{\pm}(\varpi,0)=0$ and $DN^{\pm}(\varpi,0)={\rm id}$, which show the first two equality of (\ref{NF-1-0}). Moreover, the discussion given in the first paragraph of this proof indicates that $N^{\pm}(\varpi,\cdot)$ can be used to eliminate the term $x_ux_s$, i.e., the terms $x_ix_\kappa$ for all $i=1,...,\tau$ and $\kappa=\tau+1,...,{p}$, such that the third equality of (\ref{NF-1-0}) holds.

Finally, by the second inequality of (\ref{Est-a}), we see that $D^2N^{\pm 1}(\varpi,0)=D^2N^{\pm 1}(\omega,\bar x,0)$ are measurable with respect to $\omega$ and are continuous with respect to $\bar x$ such that
\begin{align*}
\|D^2N^{\pm 1}(\varpi,0)-D^2N^{\pm 1}(\tilde\varpi,0)\|=\|a(\varpi)-a(\tilde \varpi)\|\le
C_\lambda K_\epsilon(\omega)^3\hat L_{\epsilon, \rho}(\omega)\|\bar x-\bar y\|^{\beta_N},
\end{align*}
which proves (\ref{estm-N-bar}).
The proof of Proposition \ref{prop-NF} is completed. \qquad $\Box$

\section{Solving the decomposed cohomological equation when $\alpha\in (0,1]$}
\setcounter{equation}{0}

The main task of this section is to solve equation (\ref{Distb-RRS}), i.e., proving {Proposition \ref{prop-ID-RS}}, in the case of $\alpha\in (0,1]$. Before that we prove Lemma \ref{lm-lin}. It helps us find $\varphi_{i,*}$ to defined $\bar \varphi_{i}$ in (\ref{def-p*}), which is used in equation (\ref{Distb-RRS}).

{\bf Proof of Lemma \ref{lm-lin}}.
By the same manner as (\ref{DD-DS}), we define $\bar {\bf F}_{i}(n,\varpi,x_i)$ for all $i=1,...,\tau$ and, similarly to
(\ref{int-m-d}) (using (\ref{ED-tmp}) and (\ref{FFOO-RS}) instead of (\ref{RS-ED}) and (\ref{c-off-RS}) respectively), we have that
\begin{align}
\begin{split}
&\|D\bar {\bf F}_{i}(n,\varpi,x_i)\| \leq K_\epsilon(\omega)
e^{(\lambda_i-\epsilon)n},
\qquad \forall n< 0,
\\
&\|D\bar {\bf F}_{i}(n,\varpi,x_i)\| \leq K_\epsilon(\omega)
e^{(\lambda_i+\epsilon)n},
\qquad \forall n\geq 0.
\end{split}
\label{BF-ED}
\end{align}
Define the mapping $\varphi_{i,*}(\varpi,\cdot):X_i\to X_i$ by
\begin{align}
\varphi_{i,*}(\varpi,x_i):=\lim_{k\to +\infty}\bar{\bf \Lambda}_{i}(k,\vartheta^{-k}\varpi)\bar {\bf F}_{i}(-k,\varpi,x_i).
\label{def-lin*}
\end{align}
Then, it is clear that $D\varphi_{i,*}^{\pm 1}(\varpi,{{0}})={\rm id}_i$. Moreover, the equalities
\begin{align}
\bar{\bf \Lambda}_{i}(\varpi)\varphi_{i,*}(\varpi,x_i)
&=\lim_{k\to +\infty}\bar{\bf \Lambda}_{i}(k+1,\vartheta^{-k}\varpi)\bar {\bf F}_{i}(-k-1,\vartheta\varpi,\cdot)\circ \bar {\bf F}_{i}(\varpi,x_i)
\nonumber\\
&=\varphi_{i,*}(\vartheta\varpi,\cdot)\circ \bar {\bf F}_{i}(\varpi,x_i)
\label{Con-eqn*-1}
\end{align}
imply that $\varphi_{i,*}(\varpi,x_i)$ linearizes $\bar {\bf F}_{i}(\varpi,x_i)$, which proves (\ref{Con-eqn*}).

In order to prove (\ref{phi-Bd}), we note that
\begin{align*}
&\lim_{k\to +\infty}\bar{\bf \Lambda}_{i}(k,\vartheta^{-k}\varpi)\bar {\bf F}_{i}(-k,\varpi,x_i)
-x_i
\\
&=\sum_{k=0}^{+\infty}\Big\{\bar{\bf \Lambda}_{i}(k+1,\vartheta^{-k-1}\varpi)\bar {\bf F}_{i}(-k-1,\varpi,x_i)-\bar{\bf \Lambda}_{i}(k,\vartheta^{-k}\varpi)\bar {\bf F}_{i}(-k,\varpi,x_i)\Big\}
\end{align*}
and that for all $k\ge 0$
\begin{align*}
&\frac{\|\bar{\bf \Lambda}_{i}(k+1,\vartheta^{-k-1}\varpi)D\bar {\bf F}_{i}(-k-1,\varpi,x_i)-\bar{\bf \Lambda}_{i}(k,\vartheta^{-k}\varpi)D\bar {\bf F}_{i}(-k,\varpi,x_i)\|}{\|x_i\|^{\beta}}
\\
&\le \frac{1}{\|x_i\|^{\beta}}~\|\bar{\bf \Lambda}_{i}(k+1,\vartheta^{-k-1}\varpi)\|
\\
&\quad\cdot\|D\bar {\bf F}_{i}^{-1}(\vartheta^{-k-1}\varpi,\cdot)\circ \bar {\bf F}_{i}(-k,\varpi,x_i)-\bar{\bf \Lambda}_{i}^{-1}(\vartheta^{-k-1}\varpi)\|\,\|D\bar {\bf F}_{i}(-k,\varpi,x_i)\|
\\
&\le \|\bar{\bf \Lambda}_{i}(k+1,\vartheta^{-k-1}\varpi)\|
L^*_\epsilon(\theta^{-k-1}\omega)\sup_{\xi_i\in X_i}\|D\bar {\bf F}_{i}(-k,\varpi,\xi_i)\|^{\beta}
\|D\bar {\bf F}_{i}(-k,\varpi,x_i)\|
\\
&\le K_\epsilon(\theta^{-k-1}\omega)e^{(\lambda_i+\epsilon)(k+1)}
L^*_\epsilon(\theta^{-k-1}\omega)
\{{K}_\epsilon(\omega)e^{-(\lambda_i-\epsilon)k}\}^{1+\beta}
\\
&\le e^{\lambda_i+3\epsilon} {K}_\epsilon(\omega)^3L^*_\epsilon(\omega)
e^{(-\beta \lambda_i+5\epsilon)k}
\end{align*}
by (\ref{ED-tmp}), (\ref{FFOO-RS}) and (\ref{BF-ED}).
Since $-\beta \lambda_i+5\epsilon\le -(\lambda_{j}-\lambda_{j+1})\lambda_i/\lambda_{\max}+5\epsilon<0$ for all $i=1,...\tau$ by (\ref{def-bbb}) and (\ref{def-ep}), we see that $\varphi_{i,*}(\varpi,\cdot)$ is $C^1$ and satisfies
\begin{align*}
\|D\varphi_{i,*}(\varpi,x_i)-{\rm id}_i\|
&\le\sum_{k=0}^{\infty}\big\|\bar{\bf \Lambda}_{i}(k+1,\vartheta^{-k-1}\varpi)D\bar {\bf F}_{i}(-k-1,\varpi,x_i)
\\
&\qquad\quad-\bar{\bf \Lambda}_{i}(k,\vartheta^{-k}\varpi)D\bar {\bf F}_{i}(-k,\varpi,x_i)\big\|
\\
&\le C^*_{4\epsilon}(\omega) \|x_i\|^{\beta},
\quad \forall x_i\in X_i,
\end{align*}
where $C^*_{4\epsilon}(\omega)
:=C_{\lambda}{K}_\epsilon(\omega)^3L^*_\epsilon(\omega)>0$.
As in \cite[pp.30-32]{LZZ}, one can use the standard
argument as the proof for the well-known Inverse Function Theorem to prove that the inverse
$
\varphi_{i,*}^{-1}(\varpi,x_i)
$
exists and satisfies that
$
\|D\varphi_{i,*}^{-1}(\varpi,x_i)-{\rm id}_i\|\le C^*_{4\epsilon}(\omega)\|x_i\|^{\beta}.
$
Thus, (\ref{phi-Bd}) is proved and the proof of Lemma \ref{lm-lin} is completed. \qquad $\Box$

Next, in order to solve equation (\ref{Distb-RRS}), we will see that the key step is to show the invertibility of a cohomological operator $T_i$ (see (\ref{eqn-jj}) below), deduced from the linear part of the equation, under some special norms (see Lemma \ref{lm-LO}). Then, we will show that another operator ${\mathcal B}_i$ deduced from the nonlinear part of the equation is contractive with a small { Lipschitz constant}. Finally, we obtain the solution by a recursive approximation.

{\bf Proof of {Proposition \ref{prop-ID-RS}}}.
For every $n\in\mathbb{Z}$,
let $\eta_{i}(\vartheta^n\varpi, \cdot): X_i\to \mathbb{R}^d$ be a mapping defined by
$$
\eta_{i}(\vartheta^n\varpi, x_i):=\sum_{j=1}^{{p}}\eta_{i,j}(\vartheta^n\varpi,x_i),\qquad \forall x_i\in X_i,
$$
where $\eta_{i,j}(\vartheta^n\varpi, \cdot):X_i\to X_j$,
and let
\begin{align}
\bar {\bf v}^{\kappa,\iota}_{i}(\vartheta^n\varpi, x_i):=\bar{\bf \Lambda}_{\kappa}(n,\varpi)\,\vec{e}_{\kappa,\iota}
+\eta_{i}(\vartheta^n\varpi, x_i),
\qquad \forall i=1,...,\tau,
\label{vkvk}
\end{align}
where $\vec{e}_{\kappa,\iota}$ are the standard unit vectors in the subspace $X_\kappa\subset X_s$ for $\kappa={\tau}+1,...,{p}$ and $\iota=1,...,d_{\kappa}$.

Putting $\bar {\bf f}(\varpi,\cdot):=\bar {\bf F}(\varpi,\cdot)-\bar {\bf \Lambda}(\varpi)$, we have
$D\bar {\bf F}(\varpi,x)=\bar {\bf \Lambda}(\varpi)+D\bar {\bf f}(\varpi,x)$, which enables us to rewrite equation (\ref{Distb-RRS}) ($\varpi$ is replaced with $\vartheta^n\varpi$) as
\begin{align*}
&\{\bar {\bf \Lambda}(\vartheta^n\varpi)+D\bar {\bf f}(\vartheta^n\varpi,\bar \varphi_i^{-1}(\vartheta^n\varpi,x_i))\}\,
\{\bar{\bf \Lambda}_{\kappa}(n, \varpi)\,\vec{e}_{\kappa,\iota}+\eta_{i}(\vartheta^n\varpi, x_i)\}
\\
&=\bar{\bf \Lambda}_{\kappa}(n+1, \varpi)\,\vec{e}_{\kappa,\iota}
+\eta_{i}(\vartheta^{n+1}\varpi, \bar{\bf \Lambda}_{i}(\vartheta^n\varpi)x_i).
\end{align*}
Since
$$
\bar {\bf \Lambda}(\vartheta^n\varpi)\bar{\bf \Lambda}_{\kappa}(n, \varpi)\vec{e}_{\kappa,\iota}=\bar {\bf \Lambda}_{\kappa}(\vartheta^n\varpi)\bar{\bf \Lambda}_{\kappa}(n, \varpi)\vec{e}_{\kappa,\iota}=\bar{\bf \Lambda}_{\kappa}(n+1, \varpi)\vec{e}_{\kappa,\iota},
$$
the above equation is equivalent to
\begin{align}
&
\eta_{i}(\vartheta^n\varpi, x_i)
-\bar {\bf \Lambda}^{-1}(\vartheta^n\varpi)\eta_{i}(\vartheta^{n+1}\varpi, \bar{\bf \Lambda}_{i}(\vartheta^n\varpi)x_i)
\nonumber\\
&=-\bar {\bf \Lambda}^{-1}(\vartheta^n\varpi)D\bar {\bf f}(\vartheta^n\varpi,\bar \varphi_i^{-1}(\vartheta^n\varpi,x_i))
\{\bar{\bf \Lambda}_{\kappa}(n,\varpi)\,\vec{e}_{\kappa,\iota}+\eta_{i}(\vartheta^n\varpi, x_i)\}
\label{eqnv}
\end{align}
with unknown set $\{\eta_{i}(\vartheta^n\varpi, \cdot)\}_{n\in \mathbb{Z}}$ of mappings.
Then, we consider the equation
\begin{align}
T_{i} \eta_{i}(\vartheta^n\varpi,x_i)={\mathcal B}_{i}\eta_{i}(\vartheta^n\varpi,x_i),\qquad \forall i=1,...,\tau,\quad\forall n\in\mathbb{Z},
\label{eqn-jj}
\end{align}
where $T_{i}:=\sum_{j=1}^{{p}}T_{i,j}$
and $T_{i,j}$ is defined by
\begin{align}
T_{i,j}\eta_{i,j}(\vartheta^n\varpi,x_i):=\eta_{i,j}(\vartheta^n\varpi,x_i)- \bar{\bf \Lambda}_{j}^{-1}(\vartheta^n\varpi)\eta_{i,j}(\vartheta^{n+1}\varpi,\bar{\bf \Lambda}_{i}(\vartheta^n\varpi)x_i),\quad\forall j=1,...,{p},\frac{}{}
\label{iijj-1}
\end{align}
and ${\mathcal B}_i$ is defined by
\begin{align}
{\mathcal B}_i\eta_{i}(\vartheta^n\varpi,x_i)
:=-\bar {\bf \Lambda}^{-1}(\vartheta^n\varpi)D\bar {\bf f}(\vartheta^n\varpi,\bar\varphi_i^{-1}(\vartheta^n\varpi,x_i))
\{\bar{\bf \Lambda}_{\kappa}(n, \varpi)\,\vec{e}_{\kappa,\iota}+\eta_{i}(\vartheta^n\varpi, x_i)\}.
\label{iijj-2}
\end{align}
Due to \cite[Definition 8.2.3]{Arn98}, we call $T_i$ and $T_{i,j}$ the {\it random cohomological operators}.

{\bf Step 1}. {\it The invertibility of $T_i$}.

By \cite[Lemma 3]{BaVa} (see also \cite{Arn98}), for $j=1,...,p$, we introduce the {\it Lyapunov norm}
$$
\|x_j\|_{\varpi,j}:=\sum_{n\in\mathbb{Z}}\|\bar{\bf \Lambda}_{j}(n,\varpi)x_j\|e^{-\lambda_jn-\epsilon|n|},\qquad \forall x_j\in X_j.
$$
It satisfies that
$
\|x_j\|\le \|x_j\|_{\varpi,j}\le \hat K_\epsilon(\omega)\|x_j\|,
$
where $\hat K_\epsilon(\omega):=C_\lambda K_\epsilon(\omega)>0$, implying
\begin{align}
\|x_j\|\le \|x_j\|_{\vartheta^n\varpi, j}\le \hat K_\epsilon(\theta^n\omega)\|x_j\|,\qquad \forall n\in\mathbb{Z}.
\label{Lp-N-2}
\end{align}
Then, \cite[Lemma 1]{BaVa} tells that
\begin{align}
\begin{split}
&\|\bar{\bf \Lambda}_{j}(n,\varpi)x_j\|_{\vartheta^n\varpi,j} \leq
e^{(\lambda_j-\epsilon)n}\|x_j\|_{\varpi,j},
\qquad \forall n< 0, \\
&\|\bar{\bf \Lambda}_{j}(n,\varpi)x_j\|_{\vartheta^n\varpi,j} \leq
e^{(\lambda_j+\epsilon)n}\|x_j\|_{\varpi,j},
\qquad \forall n\geq 0,
\end{split}
\label{ED-tmp-LN}
\end{align}
for all $x_j\in X_j$. The necessity of using the Lyapunov norms will be illustrated at the end of this section.

Next, define
\begin{align*}
\Upsilon^{b,n}_{i,j}
:=\bigg\{h_{i,j}(\vartheta^n\varpi,\cdot)\in C^0(X_i,X_j): h_{i,j}(\vartheta^n\varpi,{{0}})={{0}},\sup_{x_i\in X_i\backslash\{{{0}}\}}\frac{\|h_{i,j}(\vartheta^n\varpi, x_i)\|
_{\vartheta^n\varpi, j}}{\|x_i\|_{\vartheta^n\varpi, i}^b}<\infty\bigg\},
\end{align*}
with a real $b\in\mathbb{R}$ for $i=1,...,\tau$ and $j=1,...,{p}$, which is a Banach space equipped with the norm
\begin{align*}
\|h_{i,j}(\vartheta^n\varpi,\cdot)\|_{\Upsilon^{b,n}_{i,j}}:=\sup_{x_i\in X_i\backslash\{{{0}}\}}\frac{\|h_{i,j}(\vartheta^n\varpi, x_i)\|_{\vartheta^n\varpi, j}}{\|x_i\|_{\vartheta^n\varpi, i}^b}.
\end{align*}
Remind that we assume that $\alpha\in (0,1]$ in this section. Let $\beta>0$ be the constant obtained in {Theorem \ref{thm-PSU}}
and assume that
\begin{align}
0<\beta\le \alpha
\label{a=b}
\end{align}
without loss of generality because, otherwise, we  use $\alpha$ instead of $\beta$.
Let $\varsigma>0$ be given in (\ref{sg12}). Then,
$
\Upsilon^n_{i,j}:=\Upsilon^{1+\beta,n}_{i,j}\cap \Upsilon^{\varsigma,n}_{i,j}
$
is a Banach space equipped with the norm
\begin{align}
\|h_{i,j}(\vartheta^n\varpi,\cdot)\|_{\Upsilon^n_{i,j}}:=\max\big\{
\|h_{i,j}(\vartheta^n\varpi,\cdot)\|_{\Upsilon^{1+\beta,n}_{i,j}},
~\|h_{i,j}(\vartheta^n\varpi,\cdot)\|_{\Upsilon^{\varsigma,n}_{i,j}}\big\}.
\label{sp-nor}
\end{align}
This enables us to define
$
\Upsilon^n_{i}:=\bigoplus_{j=1}^{{p}}\Upsilon^n_{i,j},
$
which is also a Banach space equipped with the norm
$
\|h_{i}\|_{\Upsilon^n_{i}}:=\max_{j=1,...,{p}}\|h_{i,j}\|_{\Upsilon^n_{i,j}}.
$
We have the following lemma concerning the invertibility of $T_{i,j}$.
\begin{lm}\label{lm-LO}
For any $\{{h}_{i,j}(\vartheta^n\varpi, \cdot)\}_{n\in\mathbb{Z}}\subset \Upsilon^n_{i,j}$ {\rm (}$i=1,...,\tau$, $j=1,...{p}${\rm )} such that
\begin{align}
\sup_{n\in\mathbb{Z}}\{e^{-\lambda_\kappa n-10\epsilon |n|}\|{h}_{i,j}(\vartheta^n\varpi,\cdot)\|
_{\Upsilon^n_{i,j}}\}<\infty,
\label{etbd}
\end{align}
the linear operator $T_{i,j}$ defined by {\rm (\ref{iijj-1})} has an inverse $T_{i,j}^{-1}$ such that
$T_{i,j}^{-1}{h}_{i,j}(\vartheta^n\varpi,\cdot)\in \Upsilon^n_{i,j}$ for $n\in\mathbb{Z}$ and
\begin{align}
&\sup_{n\in\mathbb{Z}}\{e^{-\lambda_\kappa n-10\epsilon|n|}\|T_{i,j}^{-1}{h}_{i,j}(\vartheta^n\varpi,\cdot)\|
_{\Upsilon^n_{i,j}}\}
\le C_\lambda \sup_{n\in\mathbb{Z}}\{e^{-\lambda_\kappa n-10\epsilon |n|}\|{h}_{i,j}(\vartheta^n\varpi,\cdot)\|
_{\Upsilon^n_{i,j}}\}
\label{T-1}
\end{align}
for $i=1,...,\tau$ and $j=1,...,p$, where $C_\lambda>0$ is a constant.
\end{lm}
This lemma will be proved after the completion of this proof.
Recalling that $T_{i}:=\sum_{j=1}^{{p}}T_{i,j}$ and
\[
\|h_{i}\|_{\Upsilon^n_{i}}:=\max_{j=1,...,{p}}\|h_{i,j}\|_{\Upsilon^n_{i,j}},
\]
it is clear that $T_{i}^{-1}{h}_{i}(\vartheta^n\varpi,{{0}})={{0}}$ and we obtain from (\ref{T-1}) that
\begin{align}
&\sup_{n\in\mathbb{Z}}\{e^{-\lambda_\kappa n-10\epsilon|n|}\|T_{i}^{-1}{h}_{i}(\vartheta^n\varpi,\cdot)\|
_{\Upsilon^n_{i}}\}
\le C_\lambda \sup_{n\in\mathbb{Z}}\{e^{-\lambda_\kappa n-10\epsilon |n|}\|{h}_{i}(\vartheta^n\varpi,\cdot)\|
_{\Upsilon^n_{i}}\}
\label{T-111}
\end{align}
for any $\{{h}_{i}(\vartheta^n\varpi, \cdot)\}_{n\in \mathbb{Z}}\subset \Upsilon^n_{i}$ such that
$
\sup_{n\in\mathbb{Z}}\{e^{-\lambda_\kappa n-10\epsilon |n|}\|{h}_{i}(\vartheta^n\varpi,\cdot)\|
_{\Upsilon^n_{i}}\}
<\infty.
$

{\bf Step 2}. {\it The contraction of ${\mathcal B}_i$}.

Next, we show that the nonlinear operator ${\mathcal B}_i$ given in (\ref{iijj-2}) is contractive. Since
$
D\bar {\bf f}(\vartheta^n\varpi, x_i)
=\partial_{x_u}\bar {\bf f}(\vartheta^n\varpi,x_i)+\partial_{x_s}\bar {\bf f}(\vartheta^n\varpi,x_i)
$
and
$
\bar{\bf \Lambda}_{\kappa}(n, \varpi)\,\vec{e}_{\kappa,\iota}\in X_\kappa\subset X_s,
$
we see from (\ref{iijj-2}) that
\begin{align}
{\mathcal B}_{i}\eta_{i}(\vartheta^n\varpi,x_i)
&=-\bar {\bf \Lambda}^{-1}(\vartheta^n\varpi)D\bar {\bf f}(\vartheta^n\varpi,\bar \varphi_i^{-1}(\vartheta^n\varpi,x_i))
\{\bar{\bf \Lambda}_{\kappa}(n, \varpi)\,\vec{e}_{\kappa,\iota}+\eta_{i}(\vartheta^n\varpi, x_i)\}
\nonumber\\
&=-\bar {\bf \Lambda}^{-1}(\vartheta^n\varpi)\partial_{x_s}\bar {\bf f}(\vartheta^n\varpi,\bar \varphi_i^{-1}(\vartheta^n\varpi,x_i))\bar{\bf \Lambda}_{\kappa}(n, \varpi)\,\vec{e}_{\kappa,\iota}
\nonumber\\
&\quad-\bar {\bf \Lambda}^{-1}(\vartheta^n\varpi)D\bar {\bf f}(\vartheta^n\varpi,\bar \varphi_i^{-1}(\vartheta^n\varpi,x_i))
\eta_{i}(\vartheta^n\varpi, x_i).
\label{def-B}
\end{align}
In order to show that
$
{\mathcal B}_{i}\eta_{i}(\vartheta^n\varpi,\cdot)\in \Upsilon^n_{i},
$
we see easily that ${\mathcal B}_{i}\eta_{i}(\vartheta^n\varpi,0)=0$. Moreover, we
notice that $D(\pi_s\bar \varphi_i^{-1})(\vartheta^n\varpi,0)={\bf 0}$ because
$
\bar \varphi_{i}^{-1}(\vartheta^n\varpi,\cdot)=\varphi_{i}^{-1}(\vartheta^n\varpi,\cdot)\circ
\varphi_{i,*}^{-1}(\vartheta^n\varpi,\cdot),
$
as known from (\ref{def-p*}),
and
$$
D(\pi_s\varphi_i^{-1})(\vartheta^n\varpi,0)
=\frac{d\{x_s-\pi_s\gamma_i(\vartheta^n\varpi,x_i)\}}{d{x_i}}\Big|_{x_i=0}
={\bf 0}
$$
by the second equality of (\ref{hhOO}) and by the fact that $x_i\in X_u$ for $i=1,...,\tau$. It follows that
\begin{align}
&\partial^2_{x_s x}\bar {\bf f}(\vartheta^n\varpi,0)
D\bar \varphi_i^{-1}(\vartheta^n\varpi,0)
\nonumber\\
&=\partial^2_{x_s x_u}\bar {\bf f}(\vartheta^n\varpi,0)
D(\pi_u\bar \varphi_i^{-1})(\vartheta^n\varpi,0)
+\partial^2_{x_s x_s}\bar {\bf f}(\vartheta^n\varpi,0)
D(\pi_s\bar \varphi_i^{-1})(\vartheta^n\varpi,0)
\nonumber\\
&=\partial^2_{x_s x_s}\bar {\bf f}(\vartheta^n\varpi,0)
D(\pi_s\bar \varphi_i^{-1})(\vartheta^n\varpi,0)={\bf 0}
\label{ff00}
\end{align}
since
$
\partial^2_{x_sx_u}\bar {\bf f}(\vartheta^n\varpi,0)={\bf 0}
$
by the third equality of (\ref{NF-1-0}).

On the other hand, notice that
\begin{align*}
\|D\bar \varphi_i^{-1}(\vartheta^n\varpi,x_i)\|
&= \|D\varphi_{i}^{-1}(\vartheta^n\varpi,\varphi_{i,*}^{-1}(\vartheta^n\varpi,x_i))
D\varphi_{i,*}^{-1}(\vartheta^n\varpi,x_i)\|\le 2(1+C^*_{4\epsilon}(\theta^n\omega)\|x_i\|^\beta)
\\
&\le 2(1+C^*_{4\epsilon}(\theta^n\omega))\quad  {\rm when}~ \|x_i\|\le 1,
\end{align*}
and that
\begin{align*}
&\|D\bar \varphi_i^{-1}(\vartheta^n\varpi,x_i)-
D\bar \varphi_i^{-1}(\vartheta^n\varpi,0)\|
\\
&=
\|D\varphi_{i}^{-1}(\vartheta^n\varpi,\varphi_{i,*}^{-1}(\vartheta^n\varpi,x_i))
D\varphi_{i,*}^{-1}(\vartheta^n\varpi,x_i)
-D\varphi_{i}^{-1}(\vartheta^n\varpi,0)
D\varphi_{i,*}^{-1}(\vartheta^n\varpi,0)
\|
\\
&\le \|D\varphi_{i}^{-1}(\vartheta^n\varpi,\varphi_{i,*}^{-1}(\vartheta^n\varpi,x_i))
-D\varphi_{i}^{-1}(\vartheta^n\varpi,0)\|\,\|
D\varphi_{i,*}^{-1}(\vartheta^n\varpi,0)
\|
\\
&\quad+\|D\varphi_{i}^{-1}(\vartheta^n\varpi,\varphi_{i,*}^{-1}(\vartheta^n\varpi,x_i))\|
\,\|D\varphi_{i,*}^{-1}(\vartheta^n\varpi,x_i)
-D\varphi_{i,*}^{-1}(\vartheta^n\varpi,0)
\|
\\
&\le L_{\epsilon,\gamma}(\theta^n\omega)\|\varphi_{i,*}^{-1}(\vartheta^n\varpi,x_i)\|^\beta
+2C^*_{4\epsilon}(\theta^n\omega)\|x_i\|^\beta
\\
&\le \{L_{\epsilon,\gamma}(\theta^n\omega)
(1+C^*_{4\epsilon}(\theta^n\omega)\|x_i\|^\beta)
+2C^*_{4\epsilon}(\theta^n\omega)\}\|x_i\|^\beta
\\
&\le L_{5\epsilon,\gamma}(\theta^n\omega)\|x_i\|^\beta\quad  {\rm when}~ \|x_i\|\le 1
\end{align*}
due to (\ref{hhOO}) and (\ref{phi-Bd}), where $L_{5\epsilon,\gamma}(\omega):=L_{\epsilon,\gamma}(\omega)(1+C^*_{4\epsilon}(\omega))
+2C^*_{4\epsilon}(\omega)>0$.
Then, we see that when $0<\|x_i\|\le 1$
\begin{align}
&\frac{\|\partial_{x_s}\bar {\bf f}(\vartheta^n\varpi,\bar \varphi_i^{-1}(\vartheta^n\varpi,x_i))\|}
{\|x_i\|^{1+\beta}}
\nonumber\\
&=\frac{\|\partial_{x_s}\bar {\bf f}(\vartheta^n\varpi,\bar \varphi_i^{-1}(\vartheta^n\varpi,x_i))-
\partial_{x_s}\bar {\bf f}(\vartheta^n\varpi,\bar \varphi_i^{-1}(\vartheta^n\varpi,0))\|}
{\|x_i\|^{1+\beta}}
\nonumber\\
&\le
\frac{\sup_{c\in (0,1)}\|\partial^2_{x_s x}\bar {\bf f}(\vartheta^n\varpi,\bar \varphi_i^{-1}(\vartheta^n\varpi,c x_i))
D\bar \varphi_i^{-1}(\vartheta^n\varpi,c x_i)\|}{\|c x_i\|^\beta}
\nonumber\\
&=
\frac{\sup_{c\in (0,1)}\|\partial^2_{x_s x}\bar {\bf f}(\vartheta^n\varpi,\bar \varphi_i^{-1}(\vartheta^n\varpi,c x_i))
D\bar \varphi_i^{-1}(\vartheta^n\varpi,c x_i)-
\partial^2_{x_s x}\bar {\bf f}(\vartheta^n\varpi,0)
D\bar \varphi_i^{-1}(\vartheta^n\varpi,0)\|}{\|c x_i\|^\beta}
\nonumber\\
&\le
\frac{\sup_{c\in (0,1)}\|\partial^2_{x_s x}\bar {\bf f}(\vartheta^n\varpi,\bar \varphi_i^{-1}(\vartheta^n\varpi,c x_i))
-
\partial^2_{x_s x}\bar {\bf f}(\vartheta^n\varpi,0)\|\,\|
D\bar \varphi_i^{-1}(\vartheta^n\varpi,0)\|}{\|c x_i\|^\beta}
\nonumber\\
&\quad+
\frac{\sup_{c\in (0,1)}\|\partial^2_{x_s x}\bar {\bf f}(\vartheta^n\varpi,\bar \varphi_i^{-1}(\vartheta^n\varpi,c x_i))\|\,\|
D\bar \varphi_i^{-1}(\vartheta^n\varpi,c x_i)-
D\bar \varphi_i^{-1}(\vartheta^n\varpi,0)\|}{\|c x_i\|^\beta}
\nonumber\\
&\le \bar L_{\epsilon,\delta}(\theta^n\omega)\frac{\sup_{c\in (0,1)}\|\bar \varphi_i^{-1}(\vartheta^n\varpi,c x_i)\|^\beta}{\|c x_i\|^\beta}+\bar M_\epsilon(\theta^n\omega)
L_{5\epsilon,\gamma}(\theta^n\omega)
\nonumber\\
&\le \bar L_{\epsilon,\delta}(\theta^n\omega)\{
2(1+C^*_{4\epsilon}(\theta^n\omega))\}^\beta+ \bar M_\epsilon(\theta^n\omega)
L_{5\epsilon,\gamma}(\theta^n\omega)
\le e^{6\epsilon|n|}\bar C_{6\epsilon}(\omega)
\label{D2ff}
\end{align}
by (\ref{con-barF}), (\ref{ff00}) and (\ref{a=b}), where
$\bar C_{6\epsilon}(\omega):= \bar L_{\epsilon,\delta}(\omega)\{
2(1+C^*_{4\epsilon}(\omega))\}^\beta+ \bar M_\epsilon(\omega)
L_{5\epsilon,\gamma}(\omega)$,
and that when $\|x_i\|>1$
$$
\frac{\|\partial_{x_s}\bar {\bf f}(\vartheta^n\varpi,\bar \varphi_i^{-1}(\vartheta^n\varpi,x_i))\|}
{\|x_i\|^{1+\beta}}\le \|\partial_{x_s}\bar {\bf f}(\vartheta^n\varpi,\bar \varphi_i^{-1}(\vartheta^n\varpi,x_i))\|\le C_\lambda\delta(\theta^n\omega)\le e^{6\epsilon|n|}\bar C_{6\epsilon}(\omega)
$$
due to the first inequality of (\ref{con-barF}). It means that \begin{align}
e^{-6\epsilon|n|}\sup_{x_i\in X_i\backslash\{{{0}}\}}\frac{\|\partial_{x_s}\bar {\bf f}(\vartheta^n\varpi,\bar \varphi_i^{-1}(\vartheta^n\varpi,x_i))\|}
{\|x_i\|^{1+\beta}}\le\bar C_{6\epsilon}(\omega),
\label{D2ff-2}
\end{align}
which together with (\ref{def-B}) and (\ref{Lp-N-2}) gives
\begin{align}
&e^{-\lambda_\kappa n-10\epsilon|n|}\frac{\|{\mathcal B}_{i}\eta_{i}(\vartheta^n\varpi,x_i)\|
_{\vartheta^n\varpi,j}}{\|x_i\|_{\vartheta^n\varpi, i}^{1+\beta}}
\nonumber\\
&\le e^{-\epsilon|n|}\hat K_\epsilon(\theta^n\omega)e^{-\epsilon|n|}\|\bar {\bf \Lambda}^{-1}(\vartheta^n\varpi)\|\,e^{-6\epsilon|n|}
\frac{\|\partial_{x_s}\bar {\bf f}(\vartheta^n\varpi,\bar \varphi_i^{-1}(\vartheta^n\varpi,x_i))\|}{\|x_i\|^{1+\beta}}
\,e^{-\lambda_\kappa n-\epsilon|n|}\|\bar{\bf \Lambda}_{\kappa}(n,\varpi)\|
\nonumber\\
&\qquad+\hat K_\epsilon(\theta^n\omega) \|\bar {\bf \Lambda}^{-1}(\vartheta^n\varpi)\|\,\|D\bar {\bf f}(\vartheta^n\varpi,\bar \varphi_i^{-1}(\vartheta^n\varpi,x_i))\|\frac{
e^{-\lambda_\kappa n-10\epsilon|n|}\|\eta_{i}(\vartheta^n\varpi, x_i)\|}{\|x_i\|^{1+\beta}_{\vartheta^n\varpi, i}}
\nonumber\\
&\le \hat K_\epsilon(\omega) K_\epsilon(\omega)^2e^{\lambda_{\max}}\bar C_{6\epsilon}(\omega)+
\hat K_\epsilon(\theta^n\omega) e^{\lambda_{\max}}K_\epsilon (\theta^{n+1}\omega) C_\lambda\delta(\theta^n\omega)\,\big\{e^{-\lambda_\kappa n-10\epsilon|n|}\|\eta_{i}(\vartheta^n\varpi,\cdot)\|_{\Upsilon^n_{i}}\big\}
\nonumber\\
&\le \hat C_{9\epsilon}(\omega)+ \frac{1}{2C_\lambda} \,\big\{e^{-\lambda_\kappa n-10\epsilon|n|}\|\eta_{i}(\vartheta^n\varpi,\cdot)\|
_{\Upsilon^n_{i}}\big\},
\label{EM-B1}
\end{align}
where $\hat C_{9\epsilon}(\omega):=\hat K_\epsilon(\omega) K_\epsilon(\omega)^2e^{\lambda_{\max}}\bar C_{6\epsilon}(\omega)>0$. In (\ref{EM-B1}), we also use the first inequality of (\ref{con-barF}) as well as the facts that
$$
\|\bar {\bf \Lambda}^{-1}(\vartheta^n\varpi)\|\le
K_\epsilon (\theta^{n+1}\omega)e^{-\lambda_\kappa+\epsilon}\le e^{\epsilon|n|}e^{\lambda_{\max}}K_\epsilon (\omega) $$
due to (\ref{ED-tmp}) with $-1$ and $\vartheta^{n+1}\varpi$ in place of $n$ and $\varpi$ respectively, and
$$
\delta(\theta^n\omega)
=1/\mathfrak{M}_\epsilon(\theta^n\omega)\le
\{2 e^{\lambda_{\max}}C_\lambda^3 K_\epsilon(\theta^n\omega) K_\epsilon (\theta^{n+1}\omega)\}^{-1}
= \{2 e^{\lambda_{\max}}C_\lambda^2
\hat K_\epsilon(\theta^n\omega) K_\epsilon (\theta^{n+1}\omega)\}^{-1}
$$
due to (\ref{MM-p}). Similarly, we have
$$
e^{-6\epsilon|n|}\sup_{x_i\in X_i\backslash\{{{0}}\}}\frac{\|\partial_{x_s}\bar {\bf f}(\vartheta^n\varpi,\bar \varphi_i^{-1}(\vartheta^n\varpi,x_i))\|}
{\|x_i\|^{\varsigma}}\le\bar C_{6\epsilon}(\omega).
$$
Therefore,
\begin{align}
e^{-\lambda_\kappa n-10\epsilon|n|}\frac{\|{\mathcal B}_{i}\eta_{i}(\vartheta^n\varpi,x_i)\|_{\vartheta^n\varpi,j}}{\|x_i\|_{\vartheta^n\varpi, i}^\varsigma}
\le \hat C_{9\epsilon}(\omega)+ \frac{1}{2C_\lambda} \,\big\{e^{-\lambda_\kappa n-10\epsilon|n|}\|\eta_{i}(\vartheta^n\varpi,\cdot)\|_{\Upsilon^n_{i}}\big\}.
\label{EM-B2}
\end{align}
Then, it is clear that
$
{\mathcal B}_{i}\eta_{i}(\vartheta^n\varpi,\cdot)\in \Upsilon^n_{i}
$
provided $\eta_{i}(\vartheta^n\varpi,\cdot)\in \Upsilon^n_{i}$. Furthermore, (\ref{EM-B1})-(\ref{EM-B2}) imply that for $\psi_0:={\bf 0}\in \Upsilon^n_{i}$ (${\bf 0}$ is the zero mapping)
\begin{align}
e^{-\lambda_\kappa n-10\epsilon|n|}\|{\mathcal B}_{i}\psi_0(\vartheta^n\varpi,\cdot)\|_{\Upsilon^n_{i}}
\le \hat C_{9\epsilon}(\omega),
\label{i000}
\end{align}
and that for $\psi,\tilde \psi\in \Upsilon^n_{i}$
\begin{align}
e^{-\lambda_\kappa n-10\epsilon|n|}\|({\mathcal B}_{i}\psi-{\mathcal B}_{i}\tilde\psi)(\vartheta^n\varpi,\cdot)\|_{\Upsilon^n_{i}}
\le\frac{1}{2C_\lambda}\,\big\{e^{-\lambda_\kappa n-10\epsilon|n|}\|(\psi-\tilde \psi)(\vartheta^n\varpi,\cdot)\|_{\Upsilon^n_{i}}\big\}.
\label{i3dd}
\end{align}

{\bf Step 3}. {\it Recursive approximation of the solution}.

Now, we construct a sequence of mappings $\psi_\ell(\vartheta^n\varpi,\cdot)\in \Upsilon^n_{i}$
by recursively defining
\begin{align}
\psi_0(\vartheta^n\varpi,\cdot):={\bf 0}\in \Upsilon^n_{i} \quad {\rm and}  \quad \psi_\ell(\vartheta^n\varpi,\cdot):=T_{i}^{-1}{\mathcal B}_{i}\psi_{\ell-1}(\vartheta^n\varpi,\cdot), \quad \forall \ell\in \mathbb{N}.
\label{SSS-def}
\end{align}
In what follows, we prove that the sequence
is well defined and convergent for each $n\in\mathbb{Z}$. In fact, it is clear that
$
\psi_0(\vartheta^n\varpi,\cdot)
$
satisfies
$
\sup_{n\in\mathbb{Z}}\{e^{-\lambda_\kappa n-10\epsilon |n|}\|\psi_0(\vartheta^n\varpi,\cdot)\|
_{\Upsilon^n_{i}}\}={\bf 0}.
$
Then, we see from (\ref{T-111}) and (\ref{i000}) that
\begin{align}
&e^{-\lambda_\kappa n-10\epsilon|n|}\|\psi_1(\vartheta^n\varpi,\cdot)\|_{\Upsilon^n_{i}}
\nonumber\\
&=e^{-\lambda_\kappa n-10\epsilon|n|}\|(T_{i}^{-1}{\mathcal B}_{i}\psi_0)(\vartheta^n\varpi,\cdot)\|_{\Upsilon^n_{i}}\le C_\lambda\sup_{n\in\mathbb{Z}}
\big\{ e^{-\lambda_\kappa n-10\epsilon |n|}\|({\mathcal B}_{i}\psi_0)(\vartheta^n\varpi,\cdot)\|_{\Upsilon^n_{i}}\big\}
\nonumber\\
&
\le \hat C_{9\epsilon}(\omega)
\label{M00}
\end{align}
for every $n\in\mathbb{Z}$.
It means that we obtain a mapping
$\psi_{\ell}(\vartheta^n\varpi,\cdot)$ with $\ell=1$ such that
\begin{align}
\sup_{{n}\in\mathbb{Z}}\{e^{-\lambda_\kappa {n}-10\epsilon |{n}|}\|\psi_{\ell}(\vartheta^n\varpi,\cdot)\|
_{\Upsilon^n_{i}}\}<\infty.
\label{e-1}
\end{align}
Suppose that (\ref{e-1}) holds for a certain $\ell\ge 1$.
Then, we see from (\ref{EM-B1})-(\ref{EM-B2}) that
\begin{align*}
&\sup_{{n}\in\mathbb{Z}}\{e^{-\lambda_\kappa {n}-10\epsilon |{n}|}\|{\mathcal B}_{i}\psi_{\ell}(\vartheta^n\varpi,\cdot)\|
_{\Upsilon^n_{i}}\}
\le \hat C_{9\epsilon}(\omega)+ \frac{1}{2C_\lambda}\,
\sup_{{n}\in\mathbb{Z}}\{e^{-\lambda_\kappa {{n}}-10\epsilon|{{n}}|}\|\psi_{\ell}(\vartheta^n\varpi,\cdot)
\|_{\Upsilon^n_{i}}\}
<\infty,
\end{align*}
implying that $\psi_{\ell+1}(\vartheta^n\varpi,\cdot):=T_{i}^{-1}{\mathcal B}_{i}\psi_{\ell}(\vartheta^n\varpi,\cdot)\in \Upsilon^n_{i}$ is well defined such that
\begin{align*}
\sup_{{n}\in\mathbb{Z}}\{e^{-\lambda_\kappa {n}-10\epsilon |{n}|}
\|\psi_{\ell+1}(\vartheta^n\varpi,\cdot)\|
_{\Upsilon^n_{i}}\}
<\infty
\end{align*}
by (\ref{T-111}), i.e., (\ref{e-1}) holds with $\ell$ being replaced by $\ell+1$. Hence, by induction on $\ell$ we conclude that $\psi_{\ell}(\vartheta^n\varpi,\cdot)\in \Upsilon^n_{i}$ are well defined for all $\ell\in\mathbb{N}$.

Next, in view of
(\ref{T-111}) and (\ref{i3dd}), we compute that
\begin{align}
&\sup_{{n}\in\mathbb{Z}}
\big\{e^{-\lambda_\kappa {n}-10\epsilon |{n}|}\|(\psi_{\ell+1}-\psi_\ell)(\vartheta^n\varpi,\cdot)\|_{\Upsilon^n_{i}}\big\}
\nonumber\\
&=\sup_{{n}\in\mathbb{Z}}
\big\{e^{-\lambda_\kappa {n}-10\epsilon |{n}|}\big\|\big\{T_{i}^{-1}({\mathcal B}_{i}\psi_\ell-{\mathcal B}_{i}\psi_{\ell-1})\big\}(\vartheta^n\varpi,\cdot)\big\|_{\Upsilon^n_{i}}\big\}
\nonumber\\
&\le C_\lambda\sup_{{n}\in\mathbb{Z}}\big\{e^{-\lambda_\kappa {n}-10\epsilon |{n}|}
\big\|({\mathcal B}_{i}\psi_\ell-{\mathcal B}_{i}\psi_{\ell-1})({n},\cdot)\big\|_{\Upsilon^n_{i}}\big\}
\nonumber\\
&\le\frac{1}{2}\sup_{{n}\in\mathbb{Z}}\big\{
e^{-\lambda_\kappa {n}-10\epsilon |{n}|}
\|(\psi_\ell-\psi_{\ell-1})({n},\cdot)\|
_{\Upsilon^n_{i}}\big\}.
\label{Mnn}
\end{align}
This enables us to define
\begin{align*}
{N}_\ell:=\sup_{n\in\mathbb{Z}}
\big\{e^{-\lambda_\kappa n-10\epsilon|n|}\|(\psi_\ell-\psi_{\ell-1})({n},\cdot)\|
_{\Upsilon^n_{i}}\big\},\qquad\forall \ell\in\mathbb{N},
\end{align*}
which satisfies that
$
{N}_{\ell+1}\le (1/2) {N}_{\ell}\le (1/2)^\ell {N}_1.
$
Then
\begin{align*}
e^{-\lambda_\kappa n-10\epsilon|n|}\|\psi_{\ell+1}(\vartheta^n\varpi,\cdot)-\psi_\ell(\vartheta^n\varpi,\cdot)\|_{\Upsilon^n_{i}}\le {N}_{\ell+1}\le (1/2)^{\ell}{N}_1,
\end{align*}
which implies that the limit
$$
\eta_{i}(\vartheta^n\varpi,\cdot):=\lim_{\ell\to +\infty}\psi_\ell(\vartheta^n\varpi,\cdot)
=\sum_{\ell=0}^{+\infty}
\{\psi_{\ell+1}(\vartheta^n\varpi,\cdot)
-\psi_\ell(\vartheta^n\varpi,\cdot)\}
$$
exists for each $n\in\mathbb{Z}$ such that $\eta_{i}(\vartheta^n\varpi,\cdot)\in \Upsilon^n_{i}$ and
\begin{align}
e^{-\lambda_\kappa n-10\epsilon|n|}\|\eta_{i}(\vartheta^n\varpi,\cdot)\|_{\Upsilon^n_{i}}
&\le  \sum_{\ell=0}^{+\infty}e^{-\lambda_\kappa n-10\epsilon|n|}\|\psi_{\ell+1}(\vartheta^n\varpi,\cdot)
-\psi_\ell(\vartheta^n\varpi,\cdot)\|_{\Upsilon^n_{i}}
\nonumber\\
&\le 2{N}_1=2 \sup_{n\in\mathbb{Z}}
\big\{e^{-\lambda_\kappa n-10\epsilon|n|}\|\psi_1(\vartheta^n\varpi,\cdot)\|
_{\Upsilon^n_{i}}\big\}
\le \check M_{v}(\omega),
\label{eta3}
\end{align}
where $\check M_{v}(\omega):=2 \hat C_{9\epsilon}(\omega)$ by (\ref{M00}).
This solves equation (\ref{eqn-jj}) such that
\begin{align}
\begin{split}
&e^{-\lambda_\kappa n-12\epsilon|n|}\frac{\|\eta_{i}(\vartheta^n\varpi,x_i)\|}{\|x_i\|^{1+\beta}}
\le \hat K_\epsilon(\omega)^2  e^{-\lambda_\kappa n-10\epsilon|n|}\|\eta_{i}(\vartheta^n\varpi,\cdot)\|_{\Upsilon^n_{i}}\le M_{v}(\omega),
\\
&e^{-\lambda_\kappa n-12\epsilon|n|}\frac{\|\eta_{i}(\vartheta^n\varpi,x_i)\|}{\|x_i\|^{\varsigma}}
\le \hat K_\epsilon(\omega)^2  e^{-\lambda_\kappa n-10\epsilon|n|}\|\eta_{i}(\vartheta^n\varpi,\cdot)\|_{\Upsilon^n_{i}}\le M_{v}(\omega),
\label{yy-yy}
\end{split}
\end{align}
where $M_{v}(\omega):=\hat K_\epsilon(\omega)^2\check M_{v}(\omega)>0$.

Hence, according to (\ref{vkvk}), we  define
\begin{align*}
{\bar {\bf v}}_{ i}^{\kappa,\iota}(\vartheta^n\varpi, x_i):=\bar{\bf \Lambda}_{\kappa}(n, \varpi)\,\vec{e}_{\kappa,\iota}+\eta_{i}(\vartheta^n\varpi, x_i),
\qquad \forall i=1,...,\tau,\quad \forall n\in\mathbb{Z},
\end{align*}
which are the solutions of equations (\ref{Distb-RRS}) for every $\varpi\in\bar {\bf \Omega}$ such that ${\bar {\bf v}_{i}^{\kappa,\iota}}(\vartheta^n\varpi,{{0}})
=\bar{\bf \Lambda}_{\kappa}(n,\varpi)\,\vec{e}_{\kappa,\iota}\in X_\kappa$ and
\begin{align}
\sup_{n\in\mathbb{Z}}\sup_{x_i\in X_i\backslash\{{{0}}\}}\bigg\{e^{-\lambda_\kappa n-12\epsilon|n|}\frac{\|{\bar {\bf v}_{i}^{\kappa,\iota}}(\vartheta^n\varpi,x_i)-{\bar {\bf v}_{i}^{\kappa,\iota}}(\vartheta^n\varpi,{{0}})\|}{\|x_i\|^{1+\beta}}
\bigg\}\le M_{v}(\omega)
\label{111}
\end{align}
and
\begin{align}
\sup_{n\in\mathbb{Z}}\sup_{x_i\in X_i\backslash\{{{0}}\}}\bigg\{e^{-\lambda_\kappa n-12\epsilon|n|}\frac{\|{\bar {\bf v}_{i}^{\kappa,\iota}}(\vartheta^n\varpi,x_i)-{\bar {\bf v}_{i}^{\kappa,\iota}}(\vartheta^n\varpi,{{0}})\|}{\|x_i\|^\varsigma}
\bigg\}\le M_{v}(\omega)
\label{222}
\end{align}
for $i=1,...,\tau$ because of (\ref{yy-yy}). Thus, we prove (\ref{v-1-sig-XX}) by (\ref{111}) (putting $n=0$) and prove (\ref{v-sig-RS}) by (\ref{222}).
 This completes the proof of {Proposition \ref{prop-ID-RS}}. \qquad$\Box$

{\bf Proof of Lemma \ref{lm-LO}}. First of all, notice that for a linear operator $T$
if there are norms $\|\cdot\|_1$ and $\|\cdot\|_2$ such that $\|{\rm id}-T\|_1<1$ and $\|({\rm id}-T)^{-1}\|_2<1$ respectively, then
\begin{align}
T^{-1}h=-\sum_{k=1}^{\infty}\Big(({\rm id}-T)^{-k}h\Big)=\sum_{k=0}^{\infty}\Big(({\rm id}-T)^k h\Big)
\label{hh-12}
\end{align}
by the uniqueness of the inverse.
We show that for any $\{h_{i,j}(\vartheta^n\varpi,\cdot)\}_{n\in\mathbb{Z}}\subset \Upsilon^n_{i,j}$ satisfying (\ref{etbd}) the linear operators $T_{i,j}$ defined by (\ref{iijj-1}) are
invertible such that $T_{i,j}^{-1}h_{i,j}(\vartheta^n\varpi,{{0}})={{0}}$ and
\begin{align}
&e^{-\lambda_\kappa n-10\epsilon|n|}\|T_{i,j}^{-1}h_{i,j}(\vartheta^n\varpi,x_i)\|_{\Upsilon^{1+\beta,n}_{i,j}}
\le C_\lambda \sup_{n\in\mathbb{Z}}\Big\{e^{-\lambda_\kappa n-10\epsilon |n|}\|{h}_{i,j}(\vartheta^n\varpi,\cdot)\|
_{\Upsilon^{1+\beta,n}_{i,j}}\Big\}.
\label{a1sg}
\end{align}
In fact, we compute from (\ref{ED-tmp-LN}) that for $k\ge 0$
\begin{align}
&e^{-\lambda_\kappa n-10\epsilon|n|}\frac{\big\|\big\{({\rm id}-T_{i,j})^{-k}h_{i,j}\big\}
(\vartheta^n\varpi,x_i)\big\|_{\vartheta^n\varpi,j}}
{\|x_i\|_{\vartheta^n\varpi, i}^{1+\beta}}
\nonumber\\
&= e^{-\lambda_\kappa n-10\epsilon|n|}
     \frac{\|\bar{\bf \Lambda}_{j}(k,\vartheta^{n-k}\varpi)
     h_{i,j}(\vartheta^{n-k}\varpi,\bar{\bf \Lambda}_{i}(-k,\vartheta^n\varpi)x_i)\|
     _{\vartheta^n\varpi,j}}{\|x_i\|_{\vartheta^n\varpi, i}^{1+\beta}}
\nonumber\\
&\le e^{-\lambda_\kappa n-10\epsilon|n|}
   e^{(\lambda_j +\epsilon) k}
     \frac{\|h_{i,j}(\vartheta^{n-k}\varpi,\bar{\bf \Lambda}_{i}(-k,\vartheta^n\varpi)x_i)\|
     _{\vartheta^{n-k}\varpi,j}}{\|x_i\|_{\vartheta^n\varpi, i}^{1+\beta}}
\nonumber\\
&\le e^{(\lambda_j-\lambda_\kappa+ 11\epsilon ) k}\,\frac{e^{-\lambda_\kappa (n-k)-10\epsilon |n-k|}\|h_{i,j}(\vartheta^{n-k}\varpi,\bar{\bf \Lambda}_{i}(-k,\vartheta^n\varpi)x_i)\|_{ \vartheta^{n-k}\varpi,j}}{\|\bar{\bf \Lambda}_{i}(-k,\vartheta^n\varpi)x_i\|_{\vartheta^{n-k}\varpi,i}
^{1+\beta}}
\nonumber\\
&\quad\cdot\frac{\|\bar{\bf \Lambda}_{i}(-k,\vartheta^n\varpi)x_i\|_{ \vartheta^{n-k}\varpi,i}^{1+\beta}}
{\|x_i\|_{\vartheta^n\varpi, i}^{1+\beta}}
\nonumber\\
&\le e^{(\lambda_j-\lambda_\kappa+11\epsilon) k}\,
e^{-(\lambda_i-\epsilon)(1+\beta)k}
\Big\{e^{-\lambda_\kappa (n-k)-10\epsilon |n-k|}\|h_{i,j}(\vartheta^{n-k}\varpi,\cdot)\|
_{\Upsilon^{1+\beta,n-k}_{i,j}}\Big\}
\nonumber\\
&\le e^{(\lambda_j-\lambda_\kappa-\lambda_i-\beta\lambda_i
+13\epsilon) k}
\sup_{m\in\mathbb{Z}}\Big\{e^{-\lambda_\kappa m-10\epsilon |m|}\|h_{i,j}(\vartheta^m\varpi,\cdot)\|_{\Upsilon^{1+\beta,m}_{i,j}}\Big\}.
\label{LN-R}
\end{align}
Similarly, we compute that for $k\ge 0$
\begin{align}
&e^{-\lambda_\kappa n-10\epsilon|n|}\frac{\big\|\big\{({\rm id}-T_{i,j})^kh_{i,j}\big\}(\vartheta^n\varpi,x_i)\big\|
_{\vartheta^n\varpi, j}}{\|x_i\|_{\vartheta^n\varpi, i}^{1+\beta}}
\nonumber\\
&\le e^{-\lambda_\kappa n-10\epsilon|n|}
     \frac{\|\bar{\bf \Lambda}_{j}(-k,\vartheta^{n+k}\varpi)h_{i,j}
     (\vartheta^{n+k}\varpi,\bar{\bf \Lambda}_{i}(k,\vartheta^n\varpi)x_i)\|
     _{\vartheta^n\varpi,j}}{\|x_i\|
     _{\vartheta^n\varpi, i}^{1+\beta}}
\nonumber\\
&\le e^{-(\lambda_j-\lambda_\kappa-\lambda_i-\beta\lambda_i
-13\epsilon) k}
\sup_{m\in\mathbb{Z}}\Big\{e^{-\lambda_\kappa m-10\epsilon |m|}\|h_{i,j}(\vartheta^m\varpi,\cdot)\|_{\Upsilon^{1+\beta,m}_{i,j}}
\Big\}.
\label{T-inv-2}
\end{align}
Notice that either
\begin{align}
e^{\lambda_j-\lambda_\kappa-\lambda_i-\beta\lambda_i+13\epsilon}\in (0,1) \quad {\rm or}\quad
e^{-(\lambda_j-\lambda_\kappa-\lambda_i-\beta\lambda_i-13\epsilon)}\in (0,1)
\label{nobunch}
\end{align}
for $i=1,...,\tau$, $\kappa={\tau}+1,...,{p}$, $j=1,...,{p}$ by (\ref{EC-RS}), (\ref{def-ep}) and (\ref{def-bbb}).
It follows that either the series
\begin{align}
T_{i,j}^{-1}h_{i,j}(\vartheta^n\varpi,x_i)
&=-\sum_{k=1}^{\infty}\Big(({\rm id}-T_{i,j})^{-k}h_{i,j}\Big)(\vartheta^n\varpi,x_i)
\label{S1}
\end{align}
or the series
\begin{align}
T_{i,j}^{-1}h_{i,j}(\vartheta^n\varpi,x_i)
&=\sum_{k=0}^{\infty}\Big(({\rm id}-T_{i,j})^kh_{i,j}\Big)(\vartheta^n\varpi,x_i)
\label{S2}
\end{align}
converges such that
\begin{align*}
e^{-\lambda_\kappa n-10\epsilon|n|}\|T_{i,j}^{-1}h_{i,j}(\vartheta^n\varpi,x_i)\|_{\Upsilon^{1+\beta,n}_{i,j}}
\le C_\lambda \sup_{n\in\mathbb{Z}}\Big\{e^{-\lambda_\kappa n-10\epsilon |n|}\|{h}_{i,j}(\vartheta^n\varpi,\cdot)\|
_{\Upsilon^{1+\beta,n}_{i,j}}\Big\},
\end{align*}
which proves (\ref{a1sg}) with a constant $C_\lambda>0$.

Moreover, we show that
\begin{align}
&e^{-\lambda_\kappa n-10\epsilon|n|}\|T_{i,j}^{-1}h_{i,j}(\vartheta^n\varpi,x_i)\|_{\Upsilon^{\varsigma,n}_{i,j}}
\le C_\lambda \sup_{n\in\mathbb{Z}}\Big\{e^{-\lambda_\kappa n-10\epsilon |n|}\|{h}_{i,j}(\vartheta^n\varpi,\cdot)\|
_{\Upsilon^{\varsigma,n}_{i,j}}\Big\}.
\label{asg}
\end{align}
In fact, since (\ref{def-ep}) and (\ref{sg12}) imply that either
$$
e^{\lambda_j-\lambda_\kappa-\varsigma\lambda_i+12\epsilon}\in (0,1) \quad {\rm or}\quad
e^{-(\lambda_j-\lambda_\kappa-\varsigma\lambda_i-12\epsilon)}\in (0,1),
$$
by the same arguments as the previous paragraph we see that either of the series (\ref{S1}) or (\ref{S2}) converges such that (\ref{asg}) holds. Remark that for all $(i,\kappa,j)$
it is possible that
$$
\lambda_j-\lambda_\kappa-\lambda_i+13\epsilon<0<
\lambda_j-\lambda_\kappa-\varsigma\lambda_i+12\epsilon
$$
since $\lambda_i>0$ and $\varsigma\in (0,1)$. This fact requires us to use different equalities (\ref{S1}) and (\ref{S2}) to estimate
$
\|T_{i,j}^{-1}h_{i,j}(\vartheta^n\varpi,x_i)\|_{\Upsilon^{1+\beta,n}_{i,j}}
$
and
$
\|T_{i,j}^{-1}h_{i,j}(\vartheta^n\varpi,x_i)\|_{\Upsilon^{\varsigma,n}_{i,j}}
$
separately, which is allowed because of (\ref{hh-12}).
Then, combining (\ref{a1sg}) with (\ref{asg}), we have  $T_{i,j}^{-1}{h}_{i,j}(\vartheta^n\varpi,\cdot)\in \Upsilon^n_{i,j}$ for every $n\in\mathbb{Z}$ and (\ref{T-1}) holds.
The proof of Lemma \ref{lm-LO} is completed. \qquad$\Box$

Remark that (\ref{eqn-jj}) is a thorough decomposition of equation (\ref{InvDis-eqn-bar}) in every subspace $X_i$. This decomposition makes it possible to use (\ref{nobunch}), coming from (\ref{EC-RS}), to guarantee the convergence of (\ref{S1}) or (\ref{S2}), which is  the key step to solve (\ref{eqn-jj}). Otherwise, without the decomposition, we would need the bunching condition (\ref{R-bunching}),
which is stronger than (\ref{EC-RS})), to
guarantee the convergence of some relevant sequences.

Moreover,  the Lyapunov norm makes the two coefficients on both sides of (\ref{T-1}) to be identical (equal to $e^{-\lambda_\kappa n-10\epsilon|n|}$). This is important to the recursive construction of the solutions of equations (\ref{eqn-jj}) (see (\ref{SSS-def})). Without the Lyapunov norm, the coefficient on the left hand side would be
$e^{-\lambda_\kappa n-c\epsilon|n|}$ with $c>10$, which would cause divergence of the sequence constructed in (\ref{SSS-def}).

\section{Solving the decomposed cohomological equation when $\alpha=0$}
\setcounter{equation}{0}

In this section, we solve equation (\ref{Distb-RRS}), i.e., proving {Proposition \ref{prop-ID-RS}}, in the case of $\alpha=0$. Since the proof for (\ref{v-sig-RS}) is exactly the same as the one given in section 7, we are mainly concerned with (\ref{v-1-sig-XX}).

The main idea to prove (\ref{v-1-sig-XX}) is almost the same as the one given in section 7, but the only difference is that we cannot use the norm $\|h_{i,j}(\vartheta^n\varpi,\cdot)\|_{\Upsilon^{1+\beta,n}_{i,j}}$, given in (\ref{sp-nor}), because (\ref{D2ff}) cannot hold when $\bar {\bf f}(\vartheta^n\varpi,\cdot)$ is only $C^2$. Thus, we need to use another
appropriate norm instead. In order to define the new norm, we first give the following lemma.

\begin{lm}\label{lm-chi}
For any constants ${\epsilon}>0$ and $\lambda>0$, there is a non-decreasing function $\chi_{\varpi,i}:(0,+\infty)\to (0,+\infty)$ such that for every $n\in\mathbb{Z}$
\begin{align}
\begin{split}
&\sup_{x_i\in X_i}\frac{e^{-6\epsilon |n|}\|\partial^2_{x_s x}\bar {\bf f}(\vartheta^n\varpi,\bar\varphi_i^{-1}(\vartheta^n\varpi,x_i))
D\bar\varphi_i^{-1}(\vartheta^n\varpi,x_i)\|}
{\chi_{\varpi,i}(\|x_i\|)}\le 1,
\\
&
\lim_{s\to 0}\chi_{\varpi,i}(s)=0,\quad\chi_{\varpi,i}(e^{\lambda k}s)\le e^{{\epsilon} k}\chi_{\varpi,i}(s),\quad \forall k\ge 0,~\forall s\in(0,+\infty).
\end{split}
\label{CHI}
\end{align}
\end{lm}
{\bf Proof}. Since we do not need the measurability of $\chi_{\varpi,i}$ with respect to $\varpi$, we fix an arbitrary $\varpi\in\bar {\bf \Omega}$ and therefore a constant means a function depending on $\varpi$ in this proof (only).

Notice that $D^2\bar {\bf f}(\varpi,\cdot)$ is $\epsilon$-slowly continuous, as mentioned just below (\ref{con-barF}), and that $D\bar \varphi_i^{-1}(\varpi,\cdot)$ is $5\epsilon$-slowly continuous at $0$, as seen in (\ref{hhOO}) and (\ref{phi-Bd}). Thus, $e^{-\epsilon |n|}D^2\bar {\bf f}(\vartheta^n\varpi,{\bar \varphi_i^{-1}(\vartheta^n\varpi,\cdot)})$ and
$e^{-5\epsilon |n|}D\bar \varphi_i^{-1}(\vartheta^n\varpi,\cdot)$ are equicontinuous at $0$ with respect to $n$ for every $\varpi\in\bar {\bf \Omega}$. It follows that
$$
\partial^2_{x_s x}\bar {\bf f}(\vartheta^n\varpi,\bar \varphi_i^{-1}(\vartheta^n\varpi,\cdot))
D\bar \varphi_i^{-1}(\vartheta^n\varpi,\cdot)
$$
is $6\epsilon$-slowly continuous at $0$.

Then, by (\ref{ff00}), we see that for any $\ell\in \mathbb{N}$ there is an $m_\ell\in \mathbb{N}$ (independent of $n$ and $x_i$) such that for all $n\in \mathbb{Z}$
$$
e^{-6\epsilon |n|}\|\partial^2_{x_s x}\bar {\bf f}(\vartheta^n\varpi,\bar \varphi_i^{-1}(\vartheta^n\varpi,x_i))
D\bar \varphi_i^{-1}(\vartheta^n\varpi,x_i)\|\le \hat C e^{-{\epsilon} \ell},\quad \forall \|x_i\|\le e^{-m_{\ell}},
$$
with a constant $\hat C>0$. In fact,
choosing sufficiently large $\hat C>0$, we  put $m_0=0$.
Inductively,
when $m_\ell$ is fixed, we choose a sufficiently lager $m_{\ell+1}$ such that
\begin{align}
m_{\ell+1}-m_\ell\ge \lambda,
\qquad \forall \ell\ge 0.
\label{mll}
\end{align}
Then, define a function $\chi_{\varpi,i}:(0,+\infty)\to (0,+\infty)$ as
\begin{align}
\chi_{\varpi,i}(s)=
\left\{
  \begin{array}{lll}
   \hat C e^{-{\epsilon} \ell},&\forall s\in I_\ell:=(e^{-m_{\ell+1}},e^{-m_{\ell}}],&
   \forall \ell\ge 0,
    \vspace{2ex}\\
   \hat C,&\forall s>1,
  \end{array}
\right.
\label{def-chi}
\end{align}
which satisfies (\ref{CHI}). Clearly, $\chi_{\varpi,i}$ is non-decreasing and the first two formulae of (\ref{CHI}) hold.

Concerning the last formulae of (\ref{CHI}), we assume that $\ell\ge 0$ and $k\ge 0$ are given. Then for any $s\in I_\ell$ (implying $s\le 1$) we have either
$e^{\lambda k}s\le 1$ or $e^{\lambda k}s> 1$. In the case of
$e^{\lambda k}s\le 1$, we see that $e^{\lambda k}s\in I_{\ell-\hat k}=(e^{-m_{\ell-\hat k+1}},e^{-m_{\ell-\hat k}}]$ with an integer $\hat k\le \ell$.
It follows that $\ell-\hat k$ is the largest integer such that
\begin{align*}
e^{\lambda k}s\le e^{-m_{\ell-\hat k}}.
\end{align*}
On the other hand, when $k\le \ell$, the fact that $s\in I_\ell$ and (\ref{mll}) imply that
\begin{align*}
e^{\lambda k}s
\le e^{\lambda k-m_\ell}\le e^{(m_\ell-m_{\ell-1})+\cdots+(m_{\ell-k+1}-m_{\ell-k})-m_\ell}
\le e^{-m_{\ell-k}}.
\end{align*}
It means that
$
\ell-\hat k\ge \ell-k
$, implying
\begin{align}
\hat k\le k,
\label{hatkk}
\end{align}
When $k>\ell$,
(\ref{hatkk}) still holds true since $k>\ell\ge \hat k$. Thus,
\begin{align*}
\chi_{\varpi,i}(e^{\lambda k}s)=\hat C e^{-{\epsilon} \ell}e^{{\epsilon} \hat k}=e^{{\epsilon} \hat k}\chi_{\varpi,i}(s)\le  e^{{\epsilon} k}\chi_{\varpi,i}(s)
\end{align*}
by (\ref{def-chi}), where $e^{\lambda k}s\in I_{\ell-\hat k}$ and $s\in I_\ell$.

In the case of $e^{\lambda k}s> 1$, we have $k>\ell$. Otherwise (i.e., $k\le \ell$),
\begin{align*}
e^{\lambda k}s\le e^{\lambda \ell}s
\le e^{\lambda\ell-m_\ell}
\le e^{-m_0}=1
\end{align*}
since $\lambda\ell\le m_\ell-m_0$ by (\ref{mll}), which contradicts the fact that $e^{\lambda k}s> 1$. Thus,
$k>\ell$ implies that
\begin{align*}
\chi_{\varpi,i}(e^{\lambda k}s)=\hat C \,\le \hat C \,e^{{\epsilon} k} e^{-{\epsilon} \ell}
=e^{{\epsilon} k}\chi_{\varpi,i}(s)
\end{align*}
by (\ref{def-chi}).
Since $\chi_{\varpi,i}(e^{\lambda k}s)\le e^{{\epsilon} k}\chi_{\varpi,i}(s)$ holds clearly for $s>1$, the last inequality of (\ref{CHI}) is proved. $\qquad\Box$

{\bf Proof of {Proposition \ref{prop-ID-RS}}}. First of all, we define a new space
\begin{align*}
\Upsilon^{\chi,n}_{i,j}
:=\bigg\{h_{i,j}(\vartheta^n\varpi,\cdot)\in C^0(X_i,X_j): h_{i,j}(\vartheta^n\varpi,0)=0,\sup_{x_i\in X_i\backslash\{0\}}\frac{\|h_{i,j}(\vartheta^n\varpi,x_i)\|_{n,j}}
{\|x_i\|_{\vartheta^n\varpi, i}\chi_{\varpi,i}(\|x_i\|_{\vartheta^n\varpi, i})}<\infty\bigg\},
\end{align*}
where $\chi_{\varpi,i}$ is given in Lemma \ref{lm-chi}, to replace
$\Upsilon^{1+\beta,n}_{i,j}$ given just above (\ref{sp-nor}), and define a new norm
$$
\|h_{i,j}(\vartheta^n\varpi,\cdot)\|_{\Upsilon^{\chi,n}_{i,j}}:=
\sup_{x_i\in X_i\backslash\{0\}}\frac{\|h_{i,j}(\vartheta^n\varpi,x_i)\|_{n,j}}
{\|x_i\|_{\vartheta^n\varpi, i}\chi_{\varpi,i}(\|x_i\|_{\vartheta^n\varpi, i})}
$$
to replace
$\|h_{i,j}(\vartheta^n\varpi,\cdot)\|_{\Upsilon^{1+\beta,n}_{i,j}}$
given in (\ref{sp-nor}).
Then,
$
\Upsilon^n_{i,j}:=\Upsilon^{\chi,n}_{i,j}\cap \Upsilon^{\varsigma,n}_{i,j}
$
is a Banach space equipped with the norm
\begin{align}
\|h_{i,j}(\vartheta^n\varpi,\cdot)\|_{\Upsilon^n_{i,j}}:=\max\big\{
\|h_{i,j}(\vartheta^n\varpi,\cdot)\|_{\Upsilon^{\chi,n}_{i,j}},
~\|h_{i,j}(\vartheta^n\varpi,\cdot)\|_{\Upsilon^{\varsigma,n}_{i,j}}\big\}.
\label{sp-nor-1}
\end{align}
This enables us to define
$
\Upsilon^n_{i}:=\bigoplus_{j=1}^{{p}}\Upsilon^n_{i,j},
$
which is also a Banach space equipped with the norm
$
\|h_{i}\|_{\Upsilon^n_{i}}:=\max_{j=1,...,{p}}\|h_{i,j}\|_{\Upsilon^n_{i,j}}.
$

In order to prove (\ref{eta3}), we observe that it suffices to replace
(\ref{EM-B1}) and (\ref{a1sg}), both of which involve the norm $\|h_{i,j}(\vartheta^n\varpi,\cdot)\|_{\Upsilon^{1+\beta,n}_{i,j}}$, with new estimates involving the norm $\|h_{i,j}(\vartheta^n\varpi,\cdot)\|_{\Upsilon^{\chi,n}_{i,j}}$.

For (\ref{EM-B1}), we first use the inequalities
\begin{align*}
&e^{-6\epsilon|n|}\sup_{x_i\in X_i\backslash\{O\}}\frac{\|\partial_{x_s}\bar {\bf f}(\vartheta^n\varpi,\bar \varphi_i^{-1}(\vartheta^n\varpi,x_i))\|}
{\|x_i\|\chi_{\varpi,i}(\|x_i\|)}
\nonumber\\
&\le e^{-6\epsilon|n|}
\frac{\sup_{c\in (0,1)}\|\partial^2_{x_s x}\bar {\bf f}(\vartheta^n\varpi,\bar \varphi_i^{-1}(\vartheta^n\varpi,c x_i))
D\bar \varphi_i^{-1}(\vartheta^n\varpi,\cdot)|_{c x_i}\|}{\chi_{\varpi,i}(\|c x_i\|)}
\le 1,
\end{align*}
obtained from (\ref{CHI}), to replace (\ref{D2ff-2}). Then the same argument as (\ref{EM-B1}) gives
\begin{align}
e^{-\lambda_\kappa n-10\epsilon|n|}\frac{\|{\mathcal B}_{i}\eta_{i}(\vartheta^n\varpi,x_i)\|_{n,j}}{\|x_i\|_{\vartheta^n\varpi, i}\chi_{\varpi,i}(\|x_i\|_{\vartheta^n\varpi, i})}
\le \hat C_{9\epsilon}(\omega)+ \frac{1}{2C_\lambda} \,\big\{e^{-\lambda_\kappa n-10\epsilon|n|}\|\eta_{i}(\vartheta^n\varpi,\cdot)\|_{\Upsilon^n_{i}}\big\}.
\label{EM-B1-1}
\end{align}

Concerning (\ref{a1sg}), we  replace it by
\begin{align}
&e^{-\lambda_\kappa n-10\epsilon|n|}\|T_{i,j}^{-1}h_{i,j}(\vartheta^n\varpi,x_i)\|_{\Upsilon^{\chi,n}_{i,j}}
\le C_\lambda \sup_{n\in\mathbb{Z}}\Big\{e^{-\lambda_\kappa n-10\epsilon |n|}\|{h}_{i,j}(\vartheta^n\varpi,\cdot)\|
_{\Upsilon^{\chi,n}_{i,j}}\Big\}.
\label{a1sg-2}
\end{align}
In fact, by (\ref{ED-tmp-LN}) and (\ref{CHI})
\begin{align*}
\chi_{\varpi,i}(\|\bar{\bf \Lambda}_{i}(\pm k,\vartheta^n\varpi)x_i\|_{\vartheta^{n\pm k}\varpi,i})
\le\chi_{\varpi,i}(e^{(\lambda_i+\epsilon) k}\|x_i\|_{\vartheta^n\varpi, i})
\le  e^{{\epsilon} k}\chi_{\varpi,i}(\|x_i\|_{\vartheta^n\varpi, i})
\end{align*}
for all $n\in\mathbb{Z}$, $k\ge 0$ and $x_i\in X_i$.
Then we compute from (\ref{ED-tmp-LN}) that for $k\ge 0$
\begin{align*}
&e^{-\lambda_\kappa n-10\epsilon|n|}\frac{\big\|\big\{({\rm id}-T_{i,j})^{-k}h_{i,j}\big\}(\vartheta^n\varpi,x_i)
\big\|_{\vartheta^n\varpi,j}}
{\|x_i\|_{\vartheta^n\varpi, i}\,\chi_{\varpi,i}(\|x_i\|_{\vartheta^n\varpi, i})}
\nonumber\\
&= e^{-\lambda_\kappa n-10\epsilon|n|}
     \frac{\|\bar{\bf \Lambda}_{j}(k,\vartheta^{n-k}\varpi)h_{i,j}
     (\vartheta^{n-k}\varpi,\bar{\bf \Lambda}_{\bar x,i}(-k,\vartheta^n\varpi)x_i)\|_{\vartheta^n\varpi,j}}{\|x_i\|
     _{\vartheta^n\varpi, i}\,\chi_{\varpi,i}(\|x_i\|_{\vartheta^n\varpi, i})}
\nonumber\\
&\le e^{(\lambda_j-\lambda_\kappa+ 11\epsilon ) k}\,\frac{e^{-\lambda_\kappa (n-k)-10\epsilon |n-k|}\|h_{i,j}(\vartheta^{n-k}\varpi,\bar{\bf \Lambda}_{i}(-k,\vartheta^{n}\varpi)x_i)\|_{ \vartheta^{n-k}\varpi,j}}{\|\bar{\bf \Lambda}_{i}(-k,\vartheta^{n}\varpi)x_i\|_{ \vartheta^{n-k}\varpi,i}\,\chi_{\varpi,i}(\|\bar{\bf \Lambda}_{i}(-k,\vartheta^{n}\varpi)x_i\|_{\vartheta^{n-k}\varpi,i})}
\nonumber\\
&\quad\cdot\frac{\|\bar{\bf \Lambda}_{i}(-k,\vartheta^{n}\varpi)x_i\|_{ \vartheta^{n-k}\varpi,i}\,\chi_{\varpi,i}(\|\bar{\bf \Lambda}_{i}(-k,\vartheta^{n}\varpi)x_i\|_{\vartheta^{n-k}\varpi,i})}
{\|x_i\|_{\vartheta^n\varpi, i}\,\chi_{\varpi,i}(\|x_i\|_{\vartheta^n\varpi, i})}
\nonumber\\
&\le e^{(\lambda_j-\lambda_\kappa+13\epsilon) k}\,
e^{-(\lambda_i-\epsilon)k}
\Big\{e^{-\lambda_\kappa (n-k)-10\epsilon |n-k|}\|h_{i,j}(\vartheta^{n-k}\varpi,\cdot)\|
_{\Upsilon^{\chi,n-k}_{i,j}}\Big\}
\nonumber\\
&\le e^{(\lambda_j-\lambda_\kappa-\lambda_i
+13\epsilon) k}
\sup_{m\in\mathbb{Z}}\Big\{e^{-\lambda_\kappa m-10\epsilon |m|}\|h_{i,j}(\vartheta^m\varpi,\cdot)\|_{\Upsilon^{\chi,m}_{i,j}}\Big\}.
\end{align*}
Similarly, we compute that for $k\ge 0$
\begin{align*}
&e^{-\lambda_\kappa n-10\epsilon|n|}\frac{\big\|\big\{({\rm id}-T_{i,j})^kh_{i,j}\big\}(\vartheta^n\varpi,x_i)\big\|_{\vartheta^n\varpi, j}}{\|x_i\|_{\vartheta^n\varpi, i}\chi_{\varpi,i}(\|x_i\|_{\vartheta^n\varpi, i})}
\nonumber\\
&\le e^{-\lambda_\kappa n-10\epsilon|n|}
     \frac{\|\bar{\bf \Lambda}_{j}(-k,\vartheta^{n+k}\varpi)
     h_{i,j}(\vartheta^{n+k},\bar{\bf \Lambda}_{i}(k,\vartheta^{n}\varpi)x_i)\|
     _{\vartheta^n\varpi,j}}{\|x_i\|_{\vartheta^n\varpi, i}\chi_{\varpi,i}(\|x_i\|_{\vartheta^n\varpi, i})}
\nonumber\\
&\le e^{-(\lambda_j-\lambda_\kappa-\lambda_i
-13\epsilon) k}
\sup_{m\in\mathbb{Z}}\Big\{e^{-\lambda_\kappa m-10\epsilon |m|}\|h_{i,j}(\vartheta^m\varpi,\cdot)\|_{\Upsilon^{\chi,m}_{i,j}}\Big\}.
\end{align*}
Notice that we have either
\begin{align}
e^{\lambda_j-\lambda_\kappa-\lambda_i+13\epsilon}\in (0,1) \quad {\rm or}\quad
e^{-(\lambda_j-\lambda_\kappa-\lambda_i-13\epsilon)}\in (0,1)
\label{nobunch-2}
\end{align}
for all $i=1,...,\tau$, $\kappa={\tau}+1,...,{p}$, $j=1,...,{p}$ by (\ref{EC-RS}) and (\ref{def-ep}). Then, using (\ref{S1}) and (\ref{S2}), we prove (\ref{a1sg-2}).

Other parts of the proof is exactly the same as the ones given in section 7, and therefore using (\ref{EM-B1-1}) and (\ref{a1sg-2}) to replace (\ref{EM-B1}) and (\ref{a1sg}) respectively, we get (\ref{eta3}), i.e.,
\begin{align*}
e^{-\lambda_\kappa n-10\epsilon|n|}\|\eta_{i}(\vartheta^n\varpi,\cdot)\|_{\Upsilon^n_{i}}
\le \check M_{v}(\omega),
\end{align*}
where we recall that $\|\cdot\|_{\Upsilon^n_{i}}:=\max\big\{
\|h_{i,j}(\vartheta^n\varpi,\cdot)\|_{\Upsilon^{\chi,n}_{i,j}},
~\|h_{i,j}(\vartheta^n\varpi,\cdot)\|_{\Upsilon^{\varsigma,n}_{i,j}}\big\}$
given in this section. Hence,
\begin{align}
\frac{\|{\bar {\bf v}_{i}^{\kappa,\iota}}(\varpi,x_i)-{\bar {\bf v}_{i}^{\kappa,\iota}}(\varpi,{{0}})\|}
{\|x_i\|\,\chi_{\varpi,i}(\hat K_\epsilon(\omega)\|x_i\|)}
&\le \hat K_\epsilon(\omega)
\frac{\|{\bar {\bf v}_{i}^{\kappa,\iota}}(\varpi,x_i)-{\bar {\bf v}_{i}^{\kappa,\iota}}(\varpi,{{0}})\|_{\varpi,i}}
{\|x_i\|_{\varpi,i}\,\chi_{\varpi,i}(\|x_i\|_{\varpi,i})}
\nonumber\\
&\le \hat K_\epsilon(\omega)\sup_{n\in\mathbb{Z}}\bigg\{e^{-\lambda_\kappa n-10\epsilon|n|}
\frac{\|{\bar {\bf v}_{i}^{\kappa,\iota}}(\vartheta^n\varpi,x_i)-{\bar {\bf v}_{i}^{\kappa,\iota}}(\vartheta^n\varpi,{{0}})\|_{\vartheta^n\varpi, i}}
{\|x_i\|_{\vartheta^n\varpi, i}\,\chi_{\varpi,i}(\|x_i\|_{\vartheta^n\varpi, i})}\bigg\}
\nonumber\\
&\le \hat K_\epsilon(\omega)\sup_{n\in\mathbb{Z}}\big\{e^{-\lambda_\kappa n-10\epsilon|n|}\|\eta_{i}(\vartheta^n\varpi,\cdot)\|
_{\Upsilon^n_{i}}\big\}\le M_{v}(\omega),
\end{align}
which proves (\ref{v-1-sig-XX}) because $\chi_{\varpi,i}(\hat K_\epsilon(\omega)\|x_i\|)\to 0$ for each $\varpi$ (or $\omega$) as $\|x_i\|\to 0$ by (\ref{CHI}).
This completes the proof. \qquad $\Box$

\section{Smoothness of the distribution}
\setcounter{equation}{0}

In this section, we prove the smoothness of distributions, namely, we will show the smoothness of $E_s(\omega,\cdot)$ given in section 3 for each fixed $\omega\in \Omega$. By the remark given below {Theorem \ref{Thm-decop}}, it suffices to find a ``canonical" basis
of $E_s(\omega, x)$ independent of on $\bar x, i$ and give the relation between the ``canonical" basis and the basis
\[
\{{v_{\bar x, i}^{\kappa,\iota}}(\omega, x)\}_{\kappa={\tau}+1,...,{{p}},~\iota=1,...,d_\kappa}
\quad{\mbox{for $x\in W_i(\omega,\bar x)$ near $\bar x$},}
\]
which is obtained in {Theorem \ref{Thm-decop}}. For this purpose, we need the following Journ\'{e}'s lemma (\cite[Theorem]{Journe-88} or \cite[Theorem 3.3.1]{KN-Book11}):

\noindent {\bf Journ\'{e}'s lemma}. {\it Let ${\mathcal{W}}^1$ and ${\mathcal{W}}^2$ be two transverse foliations of an open set $U\subset \mathbb{R}^d$ with uniformly $C^{n,\sigma}$ {\rm (}$n\ge 1$, $\sigma\in(0,1)${\rm )} leaves. Suppose that a continuous function $v:U\to \mathbb{R}$ is uniformly $C^{n,\sigma}$ along each local leaf of ${\mathcal{W}}^1$ and ${\mathcal{W}}^2$. Then $v$ is $C^{n,\sigma}$ on $U$.
}

Note that the $C^1$ version of Journ\'{e}'s lemma, i.e., $n=1$ and $\sigma=0$, also holds by \cite[Lemma 19.1.10]{KB-Book95}. Moreover,
according to \cite[p.110]{KN-Book11}), in the above lemma the notion ``uniformly $C^{n,\sigma}$ leaves"
means that those leaves vary continuously in the $C^{n,\sigma}$-topology and the continuity is uniform with respect to the base point. In view of {Theorem \ref{thm-PSU}} and Corollary \ref{cor-PSU}, we have that
all the foliations obtained in section 4 have uniformly $C^{1,\beta}$ leaves on $U_1:=\{x\in\mathbb{R}^d: \|x\|< 1/2\}$.

Moreover, the notion ``a function $v(x)$ is uniformly $C^{n,\sigma}$ along each local leaf" means that $v$ is $C^{n,\sigma}$ along each leaf and all derivatives (up to order $n$) are uniformly bounded with respect to $x$ and the leaf, and the $n$-th derivative along the leave is $\sigma$-H\"older with uniform H\"older coefficient with respect to $x$ and the leaf.

We give the main result of this section as follows.

\begin{thm}\label{thm-globC1}
Suppose that the random diffeomorphism $F:\Omega \times \mathbb{R}^d\to \mathbb{R}^d$ satisfies {\rm (\ref{c-off-RS})}
and that the Lyapunov exponents satisfy {\rm (\ref{EC-RS})}.
Then $E_s(\omega, \cdot)$ and $E_u(\omega, \cdot)$ are $C^{1,\beta_\alpha}$ in $\mathbb{R}^d$, where
\begin{align}
\beta_\alpha:=\min\{\beta,\beta_v\}=\min\{\beta,\beta_E,\beta_N\}
=\min\{\beta,\beta_E,\epsilon/\lambda_{\max},\alpha\}.
\label{def-beta}
\end{align}
\end{thm}

Remind that $\beta\in (0,1]$ is given in {Theorem \ref{thm-PSU}} and $\beta_E\in (0,1]$ is given in {Theorem \ref{thm-HD}}, both of which are determined by the Lyapunov exponents.

{\bf Proof}. In this proof,
we mainly study the smoothness of $E_s(\omega, x)$ as the proof for $E_u(\omega, x)$ is similar.

{\bf Step 1}. {\it Find the ``canonical" basis $\{\zeta_{\kappa,\iota}(\omega, x)\}_{\kappa={\tau}+1,...,{p},~\iota=1,...,d_\kappa}$ according to the solution of the L-P equation}.

First of all, we claim that
\begin{align}
E_s(\omega, x)&=\{\partial_{y_s}q_0(\omega, x,\pi_s x)\,z_s\in \mathbb{R}^d:~z_s\in X_s\},\quad \forall x\in \mathbb{R}^d,
\label{EEsu}
\end{align}
with $q_0(\omega, x,y_s)$ being the first term of the sequence $(q_n(\omega, x,y_s))_{n\ge 0}$, which is a solution of the L-P equation (\ref{eqns-foli-S-RS}) such that $q_n(\omega, x,\pi_s x)={{0}}$ for all $n\ge 0$ and $\sup_{n\ge 0}\{e^{-\varrho n}\|\partial_{y_s}q_n(\omega, x,y_s)\|\}<\infty$ for any constant  $\varrho\in (\lambda_{\tau},\lambda_{\tau+1})$.
In fact, one can further verify that
$
F(n,\omega,x+q_0(\omega, x,y_s))
=F(n,\omega,x)+q_n(\omega, x,y_s)
$
and therefore
\begin{align*}
DF(n,\omega,x)\partial_{y_s}q_0(\omega, x,
\pi_s x)
=\frac{\partial F(n,\omega,x+q_0(\omega, x,
y_s))}{\partial y_s}\Big|_{y_s=\pi_s x}
=\partial_{y_s}q_n(\omega, x,\pi_s x).
\end{align*}
Thus,
\begin{align*}
\sup_{n\ge 0}\{e^{-\varrho n}\|DF(n,\omega, x)\partial_{y_s}q_0(\omega, x,
\pi_s x)\,z_s\|\}
\le\sup_{n\ge 0}\{e^{-\varrho n}\|\partial_{y_s}q_n(\omega, x,\pi_s x)\|\,\|z_s\|\}<\infty
\end{align*}
and by the argument given just below (\ref{vbE}) we have
\begin{align}
\{\partial_{y_s}q_0(\omega, x,\pi_s x)\,z_s\in \mathbb{R}^d:~z_s\in X_s\}
\subset E_s(\omega, x).
\label{q0Es}
\end{align}
Since (\ref{D-qq12-tild}) indicates that every $\partial_{y_s}q_0(\omega, x,\pi_s x)\,\vec{e}_{\kappa,\iota}$ is close to $\vec{e}_{\kappa,\iota}\in X_{\kappa,\iota}$,
we have that the linear space $\{\partial_{y_s}q_0(\omega, x,\pi_s x)\,z_s\in \mathbb{R}^d:~z_s\in X_s\}$ has the same dimension as $E_s(\omega, x)$ and therefore (\ref{EEsu}) is proved.

Hence, for every $\kappa={\tau}+1,...,{p}$ and $\iota=1,...,d_\kappa$ we define
$$
\zeta_{{\kappa,\iota}}(\omega, x):=\partial_{y_s}q_0(\omega, x,\pi_s x)\,\vec{e}_{\kappa,\iota}\in E_s(\omega, x),\quad \forall x\in \mathbb{R}^d,
$$
and see that $\{\zeta_{{\kappa,\iota}}(\omega, x)\}_{\kappa={\tau}+1,...,{p},~\iota=1,...,d_\kappa}$
is a basis (independent of $\bar x$ and $i$) of $E_s(\omega, x)$.
Next, in order to further study $\zeta_{{\kappa,\iota}}(\omega, x)$, we write for $\kappa=\tau+1,\cdots,p$
$$
X_\kappa=X_{\kappa,1}\oplus \cdots \oplus X_{\kappa,d_\kappa}\subset X_s,\qquad
x=x_{\kappa,1}+\cdots+x_{\kappa,d_\kappa},
$$
where every $X_{\kappa,\iota}:=\{c\,\vec{e}_{\kappa,\iota}\in \mathbb{R}^d:c\in\mathbb{R}\}$ is a 1-dimensional subspace and $x_{\kappa,\iota}\in X_{\kappa,\iota}$, and let
\begin{align*}
\pi_{\kappa,\iota}x:=x_{\kappa,\iota}\in X_{\kappa,\iota},\quad
\forall \kappa={\tau}+1,...,{p},\quad\forall \iota=1,...,d_\kappa.
\end{align*}
Then, for $y_{\kappa,\iota}:=\pi_{\kappa,\iota}y_s\in X_s$ we compute that
\begin{align*}
\zeta_{{\kappa,\iota}}(\omega, x)
&=\partial_{y_s}q_0(\omega, x,
\pi_s x)\,\vec{e}_{{\kappa,\iota}}=\partial_{y_{\kappa,\iota}}
q_0(\omega, x, \pi_s x)\,\vec{e}_{{\kappa,\iota}}
\\
&=\partial_{y_{\kappa,\iota}}
(\pi_sq_0)(\omega, x,
\pi_s x)\vec{e}_{{\kappa,\iota}}+\partial_{y_{\kappa,\iota}}
(\pi_uq_0)(\omega, x,
\pi_s x)\vec{e}_{{\kappa,\iota}}
\\
&=\vec{e}_{{\kappa,\iota}}+\partial_{y_{\kappa,\iota}}
(\pi_uq_0)(\omega, x,
\pi_s x)\vec{e}_{{\kappa,\iota}}
\in  X_{{\kappa,\iota}}\oplus X_u,
\end{align*}
where we use $\pi_s q_0(\omega, x, y_s)=y_s-\pi_s x$, as known in (\ref{eqns-foli-S-RS}), and
therefore
$$
\zeta_{{\kappa,\iota}}(\omega, x)\in E_s(\omega, x)\cap (X_{{\kappa,\iota}}\oplus X_u),\quad \forall x\in \mathbb{R}^d,
$$
by (\ref{q0Es}).
Moreover, since
\begin{align*}
\|\partial_{y_{\kappa,\iota}}(\pi_uq_0)(\omega, x,
\pi_s x)\|
= \|\partial_{y_s}(q_0-\pi_sq_0)(\omega, x,
\pi_s x)\|
\le \|\partial_{y_s}q_0(\omega, x,
\pi_s x)-{\rm id}_s\|\le {\tilde \delta}(\omega)
\end{align*}
by (\ref{D-qq12-tild}), we see that $\zeta_{{\kappa,\iota}}(\omega, x)$ is close to $\vec{e}_{{\kappa,\iota}}$ and
$$
\|\zeta_{{\kappa,\iota}}(\omega, x)\|=\max\{\|\vec{e}_{{\kappa,\iota}}\|,\|\partial_{y_{{\kappa,\iota}}}
(\pi_uq_0)(\omega, x,
\pi_s x)\vec{e}_{{\kappa,\iota}}\|\}=1,\quad \forall x\in \mathbb{R}^d,
$$
for $\kappa={\tau}+1,...,{p}$ and $\iota=1,...,d_\kappa$. Notice that
$
X_s+ (X_{{\kappa,\iota}}\oplus X_u)=\mathbb{R}^d
$
and therefore
$
E_s(\omega, x)+ (X_{{\kappa,\iota}}\oplus X_u)=\mathbb{R}^d
$
since $E_s(\omega, x)$ is close to $X_s$ (see (\ref{DB-us})).
It follows that
\begin{align*}
&{\rm dim}\{E_s(\omega, x)\cap (X_{{\kappa,\iota}}\oplus X_u)\}
\\
&=
{\rm dim}\{E_s(\omega, x)\}+{\rm dim}\{X_{{\kappa,\iota}}\oplus X_u\}-{\rm dim}\{E_s(\omega, x)+ (X_{{\kappa,\iota}}\oplus X_u)\}
\\
&=(d_{{\tau}+1}+\cdots+d_{p})+(d_1+\cdots+d _{{\tau}}+1)-d=1
\end{align*}
and therefore $\zeta_{{\kappa,\iota}}(\omega, x)$ is the unique unit vector in
$E_s(\omega, x)\cap (X_{{\kappa,\iota}}\oplus X_u)$
for every $x\in\mathbb{R}^d$ that is close to $\vec{e}_{{\kappa,\iota}}$.

{\bf Step 2}. {\it Give the relation between the ``canonical" basis $\{\zeta_{\kappa,\iota}(\omega, x)\}_{\kappa={\tau}+1,...,{p},~\iota=1,...,d_\kappa}$ and the basis
$\{{v_{\bar x, i}^{\kappa,\iota}}(\omega, x)\}_{\kappa={\tau}+1,...,{p},~\iota=1,...,d_\kappa}
$}.

For the relation between the ``canonical" basis $\{\zeta_{\kappa,\iota}(\omega, x)\}_{\kappa={\tau}+1,...,{p},~\iota=1,...,d_\kappa}$ and
\[\{{v_{\bar x, i}^{\kappa,\iota}}(\omega, x)\}_{\kappa={\tau}+1,...,{p},~\iota=1,...,d_\kappa}\quad
\mbox{for $x\in W_i(\omega,\bar x)$ near $\bar x$},
\]
we notice that if a vector $\tilde \zeta_{{\kappa,\iota}}(\omega, x)\ne {{0}}$ in the 1-dimensional subspace $E_s(\omega, x)\cap (X_{{\kappa,\iota}}\oplus X_u)$ satisfies that $\tilde \zeta_{{\kappa,\iota}}(\omega, x)/\|\tilde \zeta_{{\kappa,\iota}}(\omega, x)\|$ is close to $\vec{e}_{{\kappa,\iota}}$ (the distance is $<1/2$), then $\zeta_{{\kappa,\iota}}(\omega, x)=\tilde \zeta_{{\kappa,\iota}}(\omega, x)/\|\tilde \zeta_{{\kappa,\iota}}(\omega, x)\|$. Thus, we assume that
\begin{align*}
\tilde \zeta_{{\kappa,\iota}}(\omega, x)
&:=\big(v_{\bar x, i}^{{\tau}+1,1}(\omega, x),...,v_{\bar x, i}^{{\kappa,\iota}}(\omega, x)
,...,v_{\bar x, i}^{{p},d_{p}}(\omega, x)\big) \big(-c_{{\tau}+1,1},...,1-c_{{\kappa,\iota}},...,-c_{{p},d_{p}}\big)^T
\\
&\in E_s(\omega, x)
\end{align*}
and, at the same time, assume that
$$
\tilde \zeta_{{\kappa,\iota}}(\omega, x):=\big(\vec{e}_{1,1},...,\vec{e}_{{\tau},d_{{\tau}}},
\vec{e}_{{\kappa,\iota}}\big)
\big(c_{1,1},...,c_{{\tau},d_{{\tau}}},1\big)^T\in X_{{\kappa,\iota}}\oplus X_u
$$
for $c_{i,j_i}\in\mathbb{R}$, where the special forms of coefficients in the former equality (for examples, $-c_{{\tau}+1,1}$ and $1-c_{{\kappa,\iota}}$) are just for the sake of convenient computation and
the coefficient $1$ in the latter equality is used to make sure that $\tilde \zeta_{{\kappa,\iota}}(\omega, x)\ne 0$.

In order to determine all $c_{i,j_i}\in\mathbb{R}$, we consider linear algebraic equation
\begin{align*}
&
\big(v_{\bar x, i}^{{\tau}+1,1}(\omega, x),...,v_{\bar x, i}^{{\kappa,\iota}}(\omega, x)
,...,v_{\bar x, i}^{{p},d_{p}}(\omega, x)\big) \big(-c_{{\tau}+1,1},...,1-c_{{\kappa,\iota}},...,-c_{{p},d_{p}}\big)^T
\\
&=\big(\vec{e}_{1,1},...,\vec{e}_{{\tau},d_{{\tau}}},
\vec{e}_{{\kappa,\iota}}\big)
\big(c_{1,1},...,c_{{\tau},d_{{\tau}}},1\big)^T,
\end{align*}
from which we  find $\tilde \zeta_{{\kappa,\iota}}(\omega, x)\in E_s(\omega, x)\cap (X_{{\kappa,\iota}}\oplus X_u)$.
To solve the equation, we let
\begin{align}
\delta_{\bar x, i}^{\kappa,\iota}(\omega, x):={v_{\bar x, i}^{\kappa,\iota}}(\omega, x)
-\vec{e}_{\kappa,\iota}\in \mathbb{R}^d,\quad \forall \kappa={\tau}+1,...,{p},~\forall\iota=1,...,d_\kappa,
\label{delta-v}
\end{align}
which is small since ${v_{\bar x, i}^{\kappa,\iota}}(\omega, x)$ is close to $\vec{e}_{\kappa,\iota}$ for any $x\in {\mathcal W}_i(\omega, \bar x)$ near $\bar x$. Then, putting
\begin{align}
\begin{split}
&\vec{c}:=(c_{1,1},...,c_{\tau,d_\tau},c_{\tau+1,1}
,...,c_{{p},d_{p}})^T\in \mathbb{R}^d,
\\
&\Delta(\omega, x):=({\vec 0},...,{\vec 0},
\delta_{\bar x, i}^{{\tau}+1,1}(\omega, x)
,...,\delta_{\bar x, i}^{{p},d_{p}}(\omega, x)),
\end{split}
\label{cdel-1}
\end{align}
where $\vec 0\in \mathbb{R}^d$ is the zero vector, we rewrite the above linear algebraic equation as
\begin{align*}
({\rm id}+\Delta(\omega, x))\,\vec{c}=\delta_{\bar x, i}^{\kappa,\iota}
(\omega, x).
\end{align*}
The smallness of $\Delta(\omega, x)$ implies that the equation has a unique solution
\begin{align}
\vec{c}\,(\omega, x)=({\rm id}+\Delta(\omega, x))^{-1}\delta_{\bar x, i}^{\kappa,\iota}(\omega, x),
\label{cdel-2}
\end{align}
which is also small since $\|\vec{c}\,(\omega, x)\|\le 2\|\delta_{\bar x, i}^{\kappa,\iota}(\omega, x)\|$.
Thus, we compute that
\begin{align*}
\tilde \zeta_{{\kappa,\iota}}(\omega, x)
&=\big(v_{\bar x, i}^{{\tau}+1,1}(\omega, x),...,v_{\bar x, i}^{{\kappa,\iota}}(\omega, x)
,...,v_{\bar x, i}^{{p},d_{p}}(\omega, x)\big) \big(-c_{{\tau}+1,1},...,1-c_{{\kappa,\iota}},...,-c_{{p},d_{p}}\big)^T
\\
&=\vec{e}_{{\kappa,\iota}}+\delta_{\bar x, i}^{\kappa,\iota}(\omega, x)-
\sum_{\hat\kappa={\tau}+1}^{{p}}
\sum_{\hat\iota=1}^{d_{\hat\kappa}}
c_{\hat\kappa,\hat\iota}(\omega, x)
{v_{\bar x, i}^{\hat\kappa,\hat\iota}}(\omega, x)\in E_s(\omega, x)\cap (X_{{\kappa,\iota}}\oplus X_u)
\end{align*}
and see that every $\tilde \zeta_{{\kappa,\iota}}(\omega, x)$ is
close to $\vec{e}_{{\kappa,\iota}}$ and $\|\tilde \zeta_{{\kappa,\iota}}(\omega, x)\|$ is close to $1$. Thus,
$\tilde \zeta_{{\kappa,\iota}}(\omega, x)/\|\tilde \zeta_{{\kappa,\iota}}(\omega, x)\|$ is also close to
$\vec{e}_{{\kappa,\iota}}$. According to the uniqueness of the unit vector $\zeta_{{\kappa,\iota}}( x)\in E_s(\omega, x)\cap (X_{{\kappa,\iota}}\oplus X_u)$ which is close to
$\vec{e}_{{\kappa,\iota}}$, as mentioned at the end of the previous paragraph, we have
\begin{align}
\zeta_{{\kappa,\iota}}(\omega, x)
=\frac{\tilde \zeta_{{\kappa,\iota}}(\omega, x)}{\|\tilde \zeta_{{\kappa,\iota}}(\omega, x)\|}
=\frac{\vec{e}_{{\kappa,\iota}}+\delta_{\bar x, i}^{\kappa,\iota}(\omega, x)-
\sum_{\hat\kappa={\tau}+1}^{{p}}
\sum_{\hat\iota=1}^{d_{\hat\kappa}}
c_{\hat\kappa,\hat\iota}(\omega, x)
{v_{\bar x, i}^{\hat\kappa,\hat\iota}}(\omega, x)}{1+\tilde \pi_{{\kappa,\iota}}\big\{\delta_{\bar x, i}^{\kappa,\iota}(\omega, x)
-\sum_{\hat\kappa={\tau}+1}^{{p}}
\sum_{\hat\iota=1}^{d_{\hat\kappa}}
c_{\hat\kappa,\hat\iota}(\omega, x)
{v_{\bar x, i}^{\hat\kappa,\hat\iota}}(\omega, x)\big\}},
\label{zztt}
\end{align}
where $\tilde\pi_{{\kappa,\iota}}z\in\mathbb{R}$ denotes the $(d_1+\cdots+d_{\kappa-1}+\iota)$-th component of a vector $z\in \mathbb{R}^d$ and $\|\tilde \zeta_{{\kappa,\iota}}(\omega, x)\|:=\max_{\hat\kappa=1,...,p,
\,\hat\iota=1,...,d_{\hat\kappa}}
|\tilde \pi_{\hat\kappa,\hat\iota}\tilde \zeta_{{\kappa,\iota}}(\omega, x)|$.

Now, we obtain the relation between the two bases
\[
\{\zeta_{\kappa,\iota}(\omega, x)\}_{\kappa={\tau}+1,...,{p},~\iota=1,...,d_\kappa}
\quad{\rm and}\quad\{{v_{\bar x, i}^{\kappa,\iota}}(\omega, x)\}_{\kappa={\tau}+1,...,{p},~
\iota=1,...,d_\kappa}
\]
for all $x\in {\mathcal W}_i(\omega, \bar x)$ near $\bar x$, where the former one
is independent of $\bar x$ and $i$ and is well defined in the whole space $\mathbb{R}^d$.

{\bf Step 3}. {\it Show that $\zeta_{{\kappa,\iota}}(\omega, \cdot)$ is uniformly $C^{1,\beta_v}$ along the leaf of unstable foliation by the relation, obtained in step {\rm 2}, between two bases and by using the differentiability result of $v_{\bar x, i}^{\kappa,\iota}(\omega, x)$, obtained in Theorem {\rm \ref{Thm-decop}}}.

In order to show that $\zeta_{{\kappa,\iota}}(\omega, \cdot):\mathbb{R}^d\to \mathbb{R}^d$ is differentiable at any point $\bar x\in \mathbb{R}^d$ along the leaf ${\mathcal W}_i(\omega, \bar x)$ ($i=1,...,{\tau}$), for a small neighborhood $U\in\mathbb{R}^d$,
we define $\hat U:=U\times \cdots \times U$ (the number of $U$ is $d_{\tau+1}+\cdots+d_p$) and define a $C^\infty$ mapping $G_{\kappa,\iota}:\hat U\to \mathbb{R}^d$ as
\begin{align}
G_{\kappa,\iota}
(\xi_{\tau+1,1},...\xi_{p,d_p}):
=\frac{\vec{e}_{{\kappa,\iota}}+ \xi_{\kappa,\iota}-
\sum_{\hat\kappa={\tau}+1}^{{p}}
\sum_{\hat\iota=1}^{d_{\hat\kappa}}
\tilde\pi_{{\hat\kappa,\hat\iota}}
\{({\rm id}+\Delta_\xi)^{-1} \xi_{\hat\kappa,\hat\iota}\}
(\xi_{\hat\kappa,\hat\iota}
+\vec{e}_{\hat\kappa,\hat\iota})}{1+\tilde \pi_{{\kappa,\iota}}
\big\{ \xi_{\kappa,\iota}-
\sum_{\hat\kappa={\tau}+1}^{{p}}
\sum_{\hat\iota=1}^{d_{\hat\kappa}}
\tilde\pi_{{\hat\kappa,\hat\iota}}
\{({\rm id}+\Delta_\xi)^{-1} \xi_{\hat\kappa,\hat\iota}\}
(\xi_{\hat\kappa,\hat\iota}
+\vec{e}_{\hat\kappa,\hat\iota})
\big\}},
\label{zztt-2}
\end{align}
where $\Delta_\xi:=({\vec 0},...,{\vec 0},\xi_{\tau+1,1},...,\xi_{{p},d_{p}})$.
Then, we have
\begin{align}
\zeta_{{\kappa,\iota}}(\omega, x)
=G_{\kappa,\iota}(\delta^{+}_{\bar x, i}(\omega, x)),
\label{GGG}
\end{align}
where $\delta^{+}_{\bar x, i}(\omega, x):=(\delta^{\tau+1,1}_{\bar x, i}(\omega, x),...,\delta^{p,d_p}_{\bar x, i}(\omega, x))$.
In fact, notice that, when $(\xi_{\tau+1,1},...\xi_{p,d_p})=\delta^{+}_{\bar x, i}(\omega, x)$, we have, for $\hat\kappa={\tau}+1,...,{p}$ and
$\hat\iota=1,...,d_\kappa$,
$$
\xi_{{\hat\kappa,\hat\iota}}+\vec{e}_{{\hat\kappa,\hat\iota}}
=\delta^{\hat\kappa,\hat\iota}_{\bar x, i}(\omega, x)+\vec{e}_{{\hat\kappa,\hat\iota}}={v_{\bar x, i}^{\hat\kappa,\hat\iota}}(\omega, x)
$$
by (\ref{delta-v}), and
$$
\Delta_\xi=({\vec 0},...,{\vec 0},
\delta^{\tau+1,1}_{\bar x, i}(\omega, x),...,\delta^{{p},d_{p}}_{\bar x, i}(\omega, x))
=\Delta(\omega, x)
$$
by (\ref{cdel-1}), which implies that
\begin{align*}
\tilde\pi_{{\hat\kappa,\hat\iota}}
\{({\rm id}+\Delta_\xi)^{-1} \xi_{\hat\kappa,\hat\iota}\}
=\tilde\pi_{{\hat\kappa,\hat\iota}}\{({\rm id}+\Delta(\omega, x))^{-1} \delta^{\hat\kappa, \hat\iota}_{\bar x, i}(\omega, x)\}=c_{\hat\kappa, \hat\iota}(\omega,x)
\end{align*}
by (\ref{cdel-2}). Then, it follows from (\ref{zztt-2}) and (\ref{zztt}) that
(\ref{GGG}) holds.

Next, recalling that $x_i\to x(x_i)$ is a mapping from $X_i$ into ${\mathcal W}_i(\omega,\bar x)$ with $x({{0}})=\bar x$ (given below (\ref{def-vkl-D})) and noting that $d\delta_{\bar x, i}^{\hat \kappa,\hat \iota}=dv_{\bar x, i}^{\hat \kappa,\hat \iota}$ (see (\ref{delta-v})) for all $\hat\kappa={\tau}+1,...,{p}$ and $\hat\iota=1,...,d_\kappa$,
we see from (\ref{GGG}) that
\begin{align}
\frac{d\zeta_{{\kappa,\iota}}(\omega, x(x_i))}{dx_i}\Big|_{x_i={{0}}}
=DG_{{\kappa,\iota}}(\delta^+_{\bar x, i}(\omega, \bar x))\frac{dv_{\bar x, i}^+(\omega, x(x_i))}{dx_i}\Big|_{x_i={{0}}},
\label{Com-D-1}
\end{align}
where $v^{+}_{\bar x, i}(\omega, x):=(v^{\tau+1,1}_{\bar x, i}(\omega, x),...,v^{p,d_p}_{\bar x, i}(\omega, x))$. By {Theorem \ref{Thm-decop}} and by similar arguments to the ones given below (\ref{Dv-WU}), we conclude that
$\zeta_{{\kappa,\iota}}(\omega,x)$ is differentiable at $x=\bar x$ along the leaf ${\mathcal W}_i(\omega,\bar x)$
with a $\beta_v$-H\"older continuous differential (with respect to $\bar x$), where the differential and the H\"older coefficient are uniformly bounded with respect to $\bar x$ and $i$. Since $\bar x\in \mathbb{R}^d$ can be chosen arbitrarily, we conclude that $\zeta_{\kappa,\iota}(\omega, \cdot)$
is uniformly $C^{1,\beta_v}$ along the leaf for every $\kappa={\tau}+1,...,{p}$ and
$\iota=1,...,d_\kappa$.

Now, we show that $\zeta_{\kappa,\iota}(\omega, \cdot)$ is uniformly $C^{1,\beta_\alpha}$
along every leaf ${\mathcal W}_{u}(\omega, x)$ of the unstable foliation, where $\beta_\alpha$ is given in (\ref{def-beta}). This conclusion is clear when $\tau=1$ since ${\mathcal W}_{u}(\omega, x)={\mathcal W}_1(\omega, x)$ due to (\ref{inter-foli}) and (\ref{uuss}). When $\tau\ge 2$,
by (\ref{pqpq}) we see a filtration of leaves, i.e.,
$$
{\mathcal W}_1(\omega, x)={\mathcal W}_{\le 1}^{uu}(\omega, x)
\subset \cdots \subset {\mathcal W}_{\le \tau-1}^{uu}(\omega, x) \subset {\mathcal W}_{\le \tau}^{uu}(\omega, x) :={\mathcal W}_{u}(\omega, x),
$$
and have that for any $x\in \mathbb{R}^d$,
$\{{\mathcal W}_{\le i}^{uu}(\omega, \check x)\}_{\check x\in {\mathcal W}_{\le i+1}^{uu}(\omega, x)}$ and $\{{\mathcal W}_{i+1}(\omega, \check x)\}_{\check x\in {\mathcal W}_{\le i+1}^{uu}(\omega, x)}$ ($i=1,...,\tau-1$) are transversal foliations (with uniformly $C^{1,\beta_\alpha}$ leaves by Corollary \ref{cor-PSU}) of
\[
{\mathcal W}_{\le i+1}^{uu}(\omega, x) ~~~ \mbox{(the leaf of a strong unstable foliation)},
\]
which is a $C^{2,\alpha}$ manifold and has uniformly bounded derivatives and H\"older coefficients with respect to $x$ and $i$ by {Theorem \ref{thm-PSU}}. Noting that straightening up this manifold by a $C^{2,\alpha}$ transformation (whose derivative is $C^{1,\alpha}$) will keep the uniformly $C^{1,\beta_\alpha}$ smoothness of $\zeta_{\kappa,\iota}(\omega, \cdot)$, we use the Journ\'{e}'s lemma, given at the beginning of this section, on the manifold ${\mathcal W}_{\le i+1}^{uu}(\omega, x)$. This proves inductively that $\zeta_{\kappa,\iota}(\omega, \cdot)$ is uniformly $C^{1,\beta_\alpha}$
along the leaf ${\mathcal W}_{\le i+1}^{uu}(\omega, x)$ for any $x\in \mathbb{R}^d$ and $i=1,...,\tau-1$. Hence, $\zeta_{\kappa,\iota}(\omega, \cdot)$ is uniformly $C^{1,\beta_\alpha}$
along the leaf ${\mathcal W}_{u}(\omega, x)={\mathcal W}_{\le \tau}^{uu}(\omega, x)$ of the unstable foliation.

It
is worthy mentioning that although \cite[Theorem]{Journe-88} and \cite[Theorem 3.3.1]{KN-Book11} only give the $C^{1,\beta_\alpha}$
smoothness (without uniformity, as their whole space only has two transversal foliations) of the functions, we  see that all constants in their proofs are determined by the bounds of derivatives (including their H\"older coefficients) of the leaves as well as the functions (along the leaves), which are all uniform with respect to $x$ in our paper. Hence, all the above functions have uniform smoothness.

{\bf Step 4}. {\it Show that $\zeta_{{\kappa,\iota}}(\omega, \cdot)$ is uniformly $C^{1,\beta_\alpha}$ along the leaf of stable foliation by using the classical result obtained by the L-P method}.

For the uniformly $C^{1,\beta}$ smoothness of
$
\zeta_{{\kappa,\iota}}(\omega, x)=\partial_{y_s}q_0(\omega, x,
\pi_s x)\,\vec{e}_{{\kappa,\iota}}
$
with respect to $x$ along each leaf of the stable foliation, we
fix an $\bar x\in \mathbb{R}^d$ arbitrarily. For any $x\in {\mathcal W}_{s}(\omega, \bar x)$, let $H_{x}(\omega,\cdot):X_s\to {\mathcal W}_{s}(\omega, x)\,(={\mathcal W}_{s}(\omega, \bar x))$ be defined by
$$
H_{x}(\omega, y_s):=x+q_0(\omega, x, y_s), \quad\forall y_s\in X_s,
$$
which is a $C^{2,\alpha}$ diffeomorphism by {Theorem \ref{thm-PSU}}.
Moreover, we assert that $H_{x}(\omega, y_s)$ and $q_0(\omega, x, y_s)$ satisfy the following properties:
 \begin{itemize}
\item[(P1)] $\pi_s H_{x}(\omega,y_s)=\pi_s(x+q_0(\omega, x, y_s))=y_s$;


\item[(P2)] $x+q_0(\omega, x, y_s)=\tilde x+q_0(\omega, \tilde x, y_s)$ when $\tilde x \in {\mathcal W}_{s}(\omega, x)$.
\end{itemize}
In fact, we get from (\ref{eqns-foli-S-RS}) with $n=0$ that
\begin{align*}
q_0(\omega,x,y_s)
=y_s-\pi_s x
-\sum_{k=0}^{\infty}\Lambda_u({-k-1},\theta^{k+1}\omega)
\,\Xi_{k}^{u}(q_k)(\omega,x,y_s),
\end{align*}
which implies that
$
\pi_s(x+q_0(\omega,x,y_s))=y_s.
$
This proves (P1). For (P2), we notice from (\ref{def-leaf}) that
$x+q_0(\omega,x,y_s)\in {\mathcal W}_s(\omega,x)$. On the other hand, for any $\tilde x\in {\mathcal W}_s(\omega,x)$ we similarly have
$$
\pi_s(\tilde x+q_0(\omega,\tilde x,y_s))=y_s, \qquad \tilde x+q_0(\omega,\tilde x,y_s)\in {\mathcal W}_s(\omega,\tilde x)={\mathcal W}_s(\omega,x).
$$
From the above discussion, we conclude that $x+q_0(\omega,x,y_s)$ and $\tilde x+q_0(\omega,\tilde x,y_s)$
are two points on the same leaf ${\mathcal W}_s(\omega,x)$, which have the same projection onto $X_s$, i.e., $y_s$. Hence, we see that (P2) holds since, otherwise, the leaf ${\mathcal W}_s(\omega,x)$ cannot be Lipschitzian which contradicts (\ref{D-qq12-tild}).

Then we obtain from (P2) that
$
\partial_{y_s}q_0(\omega, x, y_s)=\partial_{y_s}q_0(\omega, \tilde x, y_s)
$
for all $\tilde x \in {\mathcal W}_{s}(\omega, x)$,
which together with (P1) implies that
\begin{align}
\partial_{y_s}q_0(\omega,H_{x}(\omega,y_s),\pi_sH_{x}(\omega,y_s))
=\partial_{y_s}q_0(\omega,H_{x}(\omega,y_s),y_s)
=\partial_{y_s}q_0(\omega,x,y_s).
\label{HHxx}
\end{align}
Remind that the above $\partial_{y_s}q_0(\omega,H_{x}(\omega,y_s),\pi_sH_{x}(\omega,y_s))$
is the value of $\partial_{y_s}q_0(\omega,x,y_s)$
at the point $(x,y_s)=(H_{x}(\omega,y_s),\pi_sH_{x}(\omega,y_s))$.
Thus, using (\ref{HHxx}), the derivative of $\zeta_{{\kappa,\iota}}(\omega,x)=\partial_{y_s}q_0(\omega, x,
\pi_s x)\,\vec{e}_{{\kappa,\iota}}$
at $x=H_{x}(\omega,\pi_s x)$ along ${\mathcal W}_{s}(\omega,\bar x)$ ($={\mathcal W}_{s}(\omega, x)$) can be computed as follows:
\begin{align}
\frac{d \zeta_{{\kappa,\iota}}(\omega,H_{x}(\omega,y_s))}{d y_s}\bigg|_{y_s=\pi_s x}
&=
\frac{d\big\{\partial_{y_s}q_0(\omega,H_{x}(\omega,y_s),
\pi_sH_{x}(\omega,y_s))\,\vec{e}_{{\kappa,\iota}} \big\}}{d y_s}\bigg|_{y_s=\pi_s x}
\nonumber\\
&=\frac{d\{\partial_{y_s}q_0(\omega, x,y_s)\,\vec{e}_{{\kappa,\iota}} \}}{d y_s}\bigg|_{y_s=\pi_s x}
\nonumber\\
&=\partial^2_{y_sy_s}q_0(\omega, x,\pi_s x)\,\vec{e}_{{\kappa,\iota}}.
\label{ze-q}
\end{align}
Since $x\mapsto \partial^2_{y_sy_s}q_{0,j}(\omega,x,y_s)$
is H\"older continuous with uniformly bounded H\"older coefficients with respect to $x,y_s$ by {Theorem \ref{thm-PSU}}, employing a similar argument to (\ref{Hd-qq-33}),
we conclude that $\zeta_{\kappa,\iota}(\omega, \cdot)$ is uniformly $C^{1,\beta}$ (and therefore $C^{1,\beta_\alpha}$) along all leaves of the stable foliation because $\bar x$
is chosen arbitrarily on $U_1$.

Now, we obtain that $\zeta_{\kappa,\iota}(\omega, \cdot)$ is $C^{1,\beta_\alpha}$ using Journ\'{e}'s lemma again because
the leaves of stable and unstable foliations are transverse in $\mathbb{R}^d$ and uniformly $C^{1,\beta_\alpha}$ on $U_1$ by {Theorem \ref{thm-PSU}}. Since $\{\zeta_{\kappa,\iota}(\omega, x)\}_{\kappa={\tau}+1,...,{p}, ~\iota=1,...,d_\kappa}$ is a basis of
$E_s (\omega, x)$, we see that
$E_s (\omega, \cdot)$ is $C^{1,\beta_\alpha}$ on $U_1$.

{Finally, we extend the
smoothness of $E_s (\omega, \cdot)$ from $U_1:=\{x\in\mathbb{R}^d: \|x\|< 1/2\}$ into $\mathbb{R}^d$.
In fact, the following facts hold:
\begin{itemize}
\item[$\bullet$]
$V(\theta^n\omega)\subset U_1$ for all $n\in\mathbb{Z}$ since
the diameter $\rho(\theta^n\omega)$ of $V(\omega)$ is $<1/4$ (by (\ref{MM-p})).

\item[$\bullet$] Both $W_1:=\bigcup_{n\in\mathbb{Z}}F(n,\omega,U_1)$ and $W_2:=\mathbb{R}^d\backslash W_1$ are invariant under $F(n,\omega,\cdot)$.
\end{itemize}
Hence, we conclude that $F(n,\omega,\cdot)=\Lambda(n,\omega)$ for all $n\in\mathbb{N}$ on the invariant set $W_2$ since $W_2\cap V(\theta^n\omega)=\emptyset$, implying that $F(\theta^n\omega,\cdot)=\Lambda(\theta^n\omega)$ for all $n\in\mathbb{N}$ on $W_2$, as mentioned below (\ref{c-off-RS}). Then, by the uniqueness of the stable distribution (see the discussion for (\ref{vbE})),
\begin{align}\label{Es1}
E_s(\omega,x)=X_s\qquad \forall x\in W_2.
\end{align}
On the other hand, for any $x\in W_1$, we define
\begin{align}\label{Es2}
E_s(\omega, x):=\{DF(n,\omega,x)\}^{-1}E_s(\theta^n\omega, F(n,\omega,x))
\end{align}
with an integer $n$ such that $F(n,\omega,x)\in U_1$. Notice that (\ref{Es1}) still holds near the boundary of $W_1$. Thus,
combining (\ref{Es1}) with (\ref{Es2}), we obtain that $E_s(\omega, \cdot)$ is $C^{1,\beta}$ in $\mathbb{R}^d$, and,
similarly, $E_u (\omega, \cdot)$ is $C^{1,\beta_\alpha}$ in $\mathbb{R}^d$.}
The proof of {Theorem \ref{thm-globC1}} is completed. \qquad $\Box$

\section{Proof of the main result}
\setcounter{equation}{0}

Having those preparations in sections 3-8, we are now ready to prove our main result.

{\bf Proof of Theorem \ref{Main-ThmC1-RS}}.
We see from {Theorem \ref{thm-globC1}} that $E_s(\omega):=\{E_s(\omega, x)\subset \mathbb{R}^d: x\in \mathbb{R}^d\}$ and $E_u(\omega):=\{E_u(\omega, x)\subset \mathbb{R}^d: x\in \mathbb{R}^d\}$ are $C^{1,\beta_\alpha}$ distributions. Remind that we have already known that $E_s(\omega)$ and $E_u(\omega)$ are integrable by the existence of stable and unstable foliations and by (\ref{EEsu}).
Then, using the well-known Frobenius' theorem (see e.g. \cite[p.541]{PSW}), for each $\omega\in\Omega$ we obtain the $C^{1,\beta_\alpha}$ smoothness of random stable foliation $\{{\mathcal W}_s(\omega, x)\}_{x\in\mathbb{R}^d}$ and random  unstable foliation $\{{\mathcal W}_u(\omega, x)\}_{x\in\mathbb{R}^d}$.
Noting that $\pi_s \{x+q_0(\omega,x,y_s)\}=y_s$ and $\pi_u \{x+p_0(\omega,x,y_u)\}=y_u$ by (\ref{eqns-foli-S-RS}) and (\ref{eqns-foli-U-RS}), we get from
(\ref{def-leaf}) that $\mathcal W_s(\omega,x)$ and $\mathcal W_u(\omega,x)$ are the graphs of mappings
\begin{align}\label{varp-pq}
y_s\mapsto \pi_u \{x+q_0(\omega,x,y_s)\} \quad {\rm and} \quad y_u\mapsto \pi_s \{x+p_0(\omega,x,y_u)\}
\end{align}
respectively, which implies that $q_0(\omega,\cdot,\cdot)$
and $p_0(\omega,\cdot,\cdot)$ are $C^{1,\beta_\alpha}$ for each $\omega\in\Omega$.

In this proof, the key step is to show the measurability of the first order derivatives of $q_0(\omega,\cdot,\cdot)$
and $p_0(\omega,\cdot,\cdot)$ with respect to $\omega$.
For the purpose,
one way is to note that the differential of the distribution along the leaf, obtained in Theorem \ref{Thm-decop}, is measurable. Then, we can modify the Journ\'{e}'s Lemma and the
Frobenius' theorem in the random case to show the measurability of the derivatives of foliations. However, we adopt another more straightforward method in the following.

First of all, straightening up the stable and unstable manifolds obtained in Theorem \ref{thm-PSU} by a tempered $C^{2,\alpha}$ random differmorphism, we can assume that
\begin{align}\label{SWSU}
{\mathcal W}_s(\omega, 0)=X_s\quad{\rm and}\quad {\mathcal W}_u(\omega,0)=X_u
\end{align}
at the very beginning of the proof, i.e., at the beginning of section 3.
{Then, we define a random continuous mapping $\phi(\omega,\cdot):\mathbb{R}^d\to \mathbb{R}^d$ as
\begin{align}\label{def-phi++}
\phi(\omega,x):=x+\pi_s p_0(\omega,x,0)+\pi_u q_0(\omega,x,0),
\end{align}
and give the following lemma.

\begin{lm}\label{lm-decp}
The inverse
of $\phi(\omega,\cdot)$ exists such that
$
\phi^{-1}(\omega,z)\in\mathbb{R}^d
$
for each $\omega\in \Omega$ and $z\in\mathbb{R}^d$ is the unique fixed point of the mapping
\begin{align}\label{def-phi--}
x\mapsto z+\pi_s p_0(\omega,\pi_s z,\pi_u x)+\pi_u q_0(\omega,\pi_u z,\pi_s x),
\end{align}
and that $\phi^{-1}(\omega,\cdot)$ is also a random continuous mapping. Moreover,
\begin{align}\label{decp}
\phi(\theta\omega, \cdot)\circ F(\omega, \cdot) \circ \phi^{-1}(\omega, \cdot)=F_s(\omega,\cdot)
+F_u(\omega,\cdot),
\end{align}
where
$
F_s(\omega,\cdot):=\pi_s F(\omega,\cdot)|_{X_s}
$
and
$
F_u(\omega,\cdot):=\pi_u F(\omega,\cdot)|_{X_u}.
$
\end{lm}
}

This lemma means that the random homeomorphism
$
\phi(\omega,\cdot):\mathbb{R}^d\to \mathbb{R}^d,
$
defined in (\ref{def-phi++}), decouples $F$ into a contraction and an expansion. A similar result was given in \cite[Proposition 7.4.3]{Arn98}, but here we construct $\phi^{-1}(\omega,\cdot)$ via $p_0$ and $q_0$.

{\bf Proof of Lemma \ref{lm-decp}}.
In order to prove that $\phi(\omega,\cdot)$ is invertible, we first
observe that
the partial derivatives of $\pi_s p_0(\omega,\pi_s z,\pi_u x)$ and
$\pi_u q_0(\omega,\pi_u z,\pi_s x)$ with respect to $\pi_u x$ and $\pi_s x$ respectively are bounded by $1/2$ due to (\ref{D-qq12-tild}), where
\begin{align*}
\pi_s p_0(\omega,\pi_s z,\pi_u x)
&=\pi_s z+p_0(\omega,\pi_s z,\pi_u x)-\pi_u \{\pi_s z+p_0(\omega,\pi_s z,\pi_u x)\}-\pi_s z
\\
&=p_0(\omega,\pi_s z,\pi_u x)-\pi_u x
\end{align*}
as known from the equality given below (\ref{varp-pq}), and similarly
$\pi_u q_0(\omega,\pi_u z,\pi_s x)
=q_0(\omega,\pi_u z,\pi_s x)-\pi_s x$.
Thus, the unique fixed point can be obtained by the Banach Contraction Principle, and the fixed point is measurable with respect to $\omega$ and continuous with respect to $z$.

Then, we recall from the proof of \cite[Theorem 3.8]{LiLu2} that $(p_k)_{k\in\mathbb{Z}}$ or $(q_k)_{k\in\mathbb{Z}}$ is the difference between trajectories with two initial points on the same unstable or stable leaf. Namely,
\begin{align}\label{intersection}
p_0(\omega,x,0)=\check x-x ~{\rm with}~
\check x\in X_s\cap {\mathcal W}_u(\omega, x),\qquad
q_0(\omega,x,0)=\hat x-x ~{\rm with}~
\hat x\in X_u\cap {\mathcal W}_s(\omega, x)
\end{align}
for all $x\in \mathbb{R}^d$, and
\begin{align}\label{diff-k}
p_k(\omega,x,0)=F(k,\omega,\check x)-F(k,\omega,x),
\qquad
q_k(\omega,x,0)=F(k,\omega,\hat x)-F(k,\omega,x)
\end{align}
for all $k\in\mathbb{Z}$. We also notice from (\ref{def-leaf}) that
\begin{align*}
\check x+p_0(\omega,\check x,\pi_u x)\in {\mathcal W}_u(\omega, \check x)={\mathcal W}_u(\omega, x),\qquad
\hat x+q_0(\omega,\hat x,\pi_s x)\in {\mathcal W}_s(\omega, \hat x)={\mathcal W}_s(\omega, x),
\end{align*}
and therefore $\check x+p_0(\omega,\check x,\pi_u x)=\hat x+q_0(\omega,\hat x,\pi_s x)=x$ by (\ref{varp-pq}).

Hence, (\ref{def-phi++}) and (\ref{intersection}) mean that $\phi(\omega,x)=\check x+\hat x$, which is totally determined by the pair of leaves ${\mathcal W}_u(\omega, x)$ and ${\mathcal W}_s(\omega, x)$,
and therefore $\phi(\omega,\cdot)$ is one-to-one because different $x$ corresponds to different pair of leaves, i.e., $\phi(\omega,\cdot)$ is invertible.  Moreover, in view of (\ref{def-phi++})-(\ref{def-phi--}) and (\ref{intersection}), we obtain that
\begin{align*}
\phi^{-1}(\omega,\phi(\omega,x))
&=\phi(\omega,x)+\pi_s p_0(\omega,\pi_s \phi(\omega,x),\pi_u x)+\pi_u q_0(\omega,\pi_u \phi(\omega,x),\pi_s x)
\\
&=x+\pi_s p_0(\omega,x,0)+\pi_u q_0(\omega,x,0)+\pi_s p_0(\omega,\pi_s (x+p_0(\omega,x,0)),\pi_u x)
\\
&\quad+\pi_u q_0(\omega,\pi_u (x+q_0(\omega,x,0)),\pi_s x)
\\
&=\pi_s (\check x+p_0(\omega,\check x,\pi_u x))
+\pi_u (\hat x+q_0(\omega,\hat x,\pi_s x))
=x,
\end{align*}
which shows that
$
\phi^{-1}(\omega,z)\in\mathbb{R}^d
$
for each $\omega\in \Omega$ and $z\in\mathbb{R}^d$ is the unique fixed point of the mapping (\ref{def-phi--}).

Next, in order to prove (\ref{decp}),
we notice from (\ref{diff-k}) that
\begin{align}\label{pkpk1}
p_k(\theta\omega, F(\omega,x),0)
&=F(k,\theta\omega,F(\omega,\check x))-F(k,\theta\omega,F(\omega,x))
\nonumber\\
&=F(k+1,\omega,\check x)-F(k+1,\omega,x)=p_{k+1}(\omega, x,0),
\end{align}
and therefore the definition of $\Xi_{k}^s(p_k)$ given below (\ref{eqns-foli-U-RS}) yields
\begin{align*}
&\Xi_{k}^s(p_k)(\theta\omega,F(\omega,x),0)
\\
&=\pi_sf\big(\theta^{k+1} \omega,p_k(\theta\omega, F(\omega,x),0)+F(k,\theta\omega,F(\omega,x))\big)
-\pi_sf\big(\theta^{k+1}\omega,F(k,\theta\omega,F(\omega,x))\big)
\\
&=\pi_sf\big(\theta^{k+1} \omega,p_{k+1}(\omega, x,0)+F(k+1,\omega,x)\big)
-\pi_sf\big(\theta^{k+1}\omega,F(k+1,\omega,x)\big)
\\
&=\Xi_{k+1}^s(p_{k+1})(\omega,x,0).
\end{align*}
Similarly, $\Xi_{k}^u(q_k)(\theta\omega,F(\omega,x),0)=
\Xi_{k+1}^u(q_{k+1})(\omega,x,0)$. Then, it follows from (\ref{def-phi++}) and (\ref{eqns-foli-U-RS})-(\ref{eqns-foli-S-RS}) that
\begin{align*}
\phi(\theta\omega, \cdot)&\circ F(\omega,x)
=
F(\omega,x)+\pi_s p_0(\theta\omega,F(\omega,x),0)+\pi_u q_0(\theta\omega,F(\omega,x),0)
\\
&=F(\omega,x)+\Lambda_s(\omega)\sum_{k=-\infty}^{-1}\Lambda_s({-k-2},\theta^{k+2}\omega)
\,\Xi_{k}^{s}(p_k)(\theta\omega,F(\omega,x),0)
\\
&\quad-\Lambda_u(\omega)\sum_{k=0}^{\infty}\Lambda_u({-k-2},\theta^{k+2}\omega)
\,\Xi_{k}^{u}(q_k)(\theta\omega,F(\omega,x),0)
\\
&=F(\omega,x)+\Lambda_s(\omega)\sum_{k=-\infty}^{0}\Lambda_s({-k-1},\theta^{k+1}\omega)
\,\Xi_{k-1}^{s}(p_{k-1})(\theta\omega,F(\omega,x),0)
\\
&\quad-\Lambda_u(\omega)\sum_{k=1}^{\infty}\Lambda_u({-k-1},\theta^{k+1}\omega)
\,\Xi_{k-1}^{u}(q_{k-1})(\theta\omega,F(\omega,x),0)
\\
&=\Lambda(\omega)x+\Lambda_s(\omega)\sum_{k=-\infty}^{-1}\Lambda_s({-k-1},\theta^{k+1}\omega)
\,\Xi_{k}^{s}(p_{k})(\omega,x,0)+ \pi_sf(\omega,x)+\Xi_{0}^{s}(p_{0})(\omega,x,0)
\\
&\quad-\Lambda_u(\omega)\sum_{k=0}^{\infty}\Lambda_u({-k-1},\theta^{k+1}\omega)
\,\Xi_{k}^{u}(q_{k})(\omega,x,0)
+ \pi_uf(\omega,x)+\Xi_{0}^{u}(q_{0})(\omega,x,0)
\\
&=\Lambda_s(\omega)\pi_sx+\Lambda_s(\omega)\pi_s p_0(\omega,x,0)
+\pi_sf(\omega,x+p_0(\omega,x,0))
\\
&\quad+\Lambda_u(\omega)\pi_ux+\Lambda_u(\omega)\pi_u q_0(\omega,x,0)+\pi_uf(\omega,x+q_0(\omega,x,0))
\\
&=\Lambda_s(\omega)\pi_s(x+p_0(\omega,x,0))
+\pi_sf(\omega,\pi_s(x+p_0(\omega,x,0)))
\\
&\quad+\Lambda_u(\omega)\pi_u(x+ q_0(\omega,x,0))+\pi_uf(\omega,\pi_u(x+q_0(\omega,x,0)))
\\
&=\pi_s F(\omega,\cdot)\circ \pi_s\phi(\omega,x)+\pi_u F(\omega,\cdot)\circ \pi_u\phi(\omega,x)
=\{F_s(\omega,\cdot)+F_u(\omega,\cdot)\}\circ \phi(\omega,x),
\end{align*}
where
$
x+p_0(\omega,x,0)=\pi_s(x+p_0(\omega,x,0))
$
and
$
x+q_0(\omega,x,0)=\pi_u(x+q_0(\omega,x,0))
$
due to (\ref{intersection}). This proves (\ref{decp})
and the proof of Lemma \ref{lm-decp} is completed.
\qquad $\Box$

We continue to prove Theorem \ref{Main-ThmC1-RS}. By Lemma \ref{lm-decp} and the classical Implicit Function Theorem,  for each $\omega$
the smoothness of $\phi^{\pm 1}(\omega,\cdot)$
is equivalent to the smoothness of $q_0(\omega,\cdot,\cdot)$ and $p_0(\omega,\cdot,\cdot)$. Furthermore, the derivatives of $\phi^{\pm 1}(\omega,\cdot)$ has the same measurability as the derivatives of $q_0(\omega,\cdot,\cdot)$ and $p_0(\omega,\cdot,\cdot)$ because the proof of the Implicit Function Theorem is essentially based on the Contraction Principle.

Notice that $q_0(\omega,x,y_s)$ is tempered $\sigma^2$-H\"older continuous with respect to $x$ on $U_1:=\{x\in\mathbb{R}^d:\|x\|\le 1/2\}$ (see (\ref{Hd-qq-22}) in the case of $j=\tau$) and is tempered $C^{2,\alpha}$ with respect to $y_s$ (see Theorem \ref{thm-PSU}). Moreover, the domain on which $q_0(\omega,\cdot,y_s)$ is tempered $\sigma^2$-H\"older continuous can be extended from $U_1$ to
$$
U_\phi(\omega):=\{z=\phi(\omega,x):x\in U_1\}~ {\rm with}~ \|\pi_u z\|<1,
$$
where $U_\phi(\omega)$ is bounded by the constant $2$ in view of (\ref{def-phi++}) and (\ref{D-qq12-tild}) (noting that $p_0(\omega,x,\pi_u x)=q_0(\omega,x,\pi_s x)=0$), i.e.,
\begin{align*}
\big\|\partial_{y_s}q_{0}(\omega, z, y_s)-
\partial_{y_s}q_{0}(\omega, \pi_s z, y_s)\big\|
\le \tilde L_\epsilon(\omega)\|\pi_u z\|^{\sigma^2}, \quad \forall z,y_s\in U_\phi(\omega).
\end{align*}
Actually, in the proof of (\ref{Hd-qq-22}), $U_1$ can be replaced with any bounded domain (see (\ref{HH-L})) and the condition $\|x-\tilde x\|<1$ is satisfied here, i.e., $\|z-\pi_s z\|=\|\pi_u z\|<1$. A similar conclusion also holds for $p_0(\omega,x,y_u)$. Hence, we have the following facts.
\begin{itemize}
\item[(B1)] The mapping $\phi(\omega,\cdot)$ is a random homeomorphism such that
\begin{align}\label{ph-ph}
\|\phi(\omega,x)-\phi(\omega,\tilde x)\|\le L_1(\omega)\|x-\tilde x\|^{\hat\sigma},\quad \forall x,\tilde x\in U_1,
\end{align}
    with a random
    variable $L_1(\omega)>0$ tempered from above and a constant $\hat\sigma:=\sigma^2\in (0,1)$.
\item[(B2)]
For any $z\in U_\phi(\omega)$ with $\|\pi_u z\|<1$, we have
\begin{align}\label{ph-su}
\|\phi^{-1}(\omega,z)-\phi^{-1}(\omega,\pi_s z)\|\le L_2(\omega)\|\pi_u z\|^{\hat\sigma}
\end{align}
with a random variable $L_2(\omega)>0$ tempered from above.

\item[(B3)]
$\phi(\omega,\cdot)$ is a $C^{1,\beta_\alpha}$ diffeomophism since $q_0(\omega,\cdot,\cdot)$
and $p_0(\omega,\cdot,\cdot)$ are $C^{1,\beta_\alpha}$ for every $\omega\in\Omega$, as mentioned below (\ref{varp-pq}), but we do not know the measurability of $D\phi^{\pm 1}(\omega,x)$ with respect to $\omega$.
\end{itemize}
In order to show that $D\phi^{\pm 1}(\omega,x)$ are measurable with respect to $\omega$, observing (\ref{def-phi++}), we need to prove that $\pi_u q_0(\omega,\cdot,0)$ and $\pi_s p_0(\omega,\cdot,0)$ have measurable derivatives. In what follows, we only deal with $q_0$ and the proof for $p_0$ is similar.

Regarding the measurability of the derivative of $q_0(\omega,\cdot,0)$, we consider the leaf ${\mathcal W}_s(\omega, x)$ of the stable foliation of $F(\omega,\cdot)$, and claim that ${\mathcal W}_s(\omega, x)$ is flat when $\|\pi_u x\|>1/2$, i.e.,
\begin{align}\label{Ws-out}
{\mathcal W}_s(\omega, x)=\{x+y_s\in \mathbb{R}^d:y_s\in X_s\}=\{z\in\mathbb{R}^d:\pi_u z=\pi_u x\}\quad {\rm if}~ \|\pi_u x\|>1/2.
\end{align}
In fact, using (\ref{ED-ge}) in the case of $j=\tau$ and $x=0$, we have
\begin{align}\label{min-u}
\|\Lambda_u(n,\omega)\pi_u x\|
\ge (1/K_\epsilon(\omega))e^{(\lambda_\tau-2\epsilon) n}\|\pi_u x\|,\quad \forall x\in \mathbb{R}^d,
\end{align}
which implies that
$
\|\Lambda_u(n,\omega)\pi_u x\|
\ge 1/(2K_\epsilon(\omega))
$
for all $n\in \mathbb{N}$ when $\|\pi_u x\|>1/2$ since $\lambda_\tau-2\epsilon>0$. Then, (\ref{Ws-out}) is true by observing two facts that $F(\omega,x)=\Lambda(\omega)x$ for all $x\notin V(\omega)$, whose diameter is $\rho(\omega)<1/(2K_\epsilon(\omega))$ (by (\ref{MM-p})),
and that the stable foliation is unique.

Since ${\mathcal W}_s(\omega, x)$ is the graph of the
mapping
$
y_s\mapsto \pi_u \{x+q_0(\omega,x,y_s)\},
$
as mentioned in (\ref{varp-pq}), we obtain from (\ref{SWSU}) and (\ref{Ws-out}) that
\begin{align}\label{uq0=0}
\quad\pi_u q_0(\omega,x,y_s)=0
\quad {\rm if}~ \pi_u x=0~{\rm or}~\|\pi_u x\|>1/2,\quad \forall y_s\in X_s.
\end{align}
Then, for all $z\in\mathbb{R}^d$ such that $\pi_u z=0$ or $\|\pi_u z\|>1/2$, we know from Lemma \ref{lm-decp} that
\begin{align}\label{uph-0out}
\pi_u\phi^{-1}(\omega,z)
&=\pi_u\{z+\pi_s p_0(\omega,\pi_s z,\pi_u \phi^{-1}(\omega,z))+\pi_u q_0(\omega,\pi_u z,\pi_s \phi^{-1}(\omega,z))\}
\nonumber\\
&=\pi_u\{z+\pi_s p_0(\omega,\pi_s z,\pi_u \phi^{-1}(\omega,z))\}=\pi_u z.
\end{align}
{In what follows, for any $x\in U_1\backslash X_s$, we want to find an integer $N\in\mathbb{N}$ (depending on $\omega$ and $x$) such that $\|\pi_uF(N,\omega,x)\|>1/2$, and then establish a relation between $\pi_u q_0(\omega,x,0)$
and $\pi_u q_0(\theta^N\omega,F(N,\omega,x),0)$, which is equal to $0$ by (\ref{uq0=0}), to show the measurability of $\partial_x (\pi_u q_0)(\omega,x,0)$ (see (\ref{qyy}) below). Basically, such an integer $N$ exists since
the norm of the unstable projection $\pi_uF(n,\omega,x)$ tends to infinity as $n\to \infty$ for any $x\in U_1\backslash X_s$, but the integer $N$ depends on $\omega$ and we do not know its measurability. In order to solve this problem, after the following preparations, we will give our Lemma \ref{lm-N}.
}

Noting that $F_s(\omega,\cdot)$ and
$F_u(\omega,\cdot)$ are tempered $C^{2,\alpha}$ random diffeomorphisms such that $DF_s(\omega,{{0}})=\Lambda_s (\omega)$ and $D F_u(\omega,{{0}})=\Lambda_u (\omega)$,
we obtain from \cite[Main Theorem]{LZZ} that there are locally tempered $C^{1,\beta_\alpha}$ (without loss of generality, we still use $\beta_\alpha$ to denote the H\"older exponent here) random diffeomorphisms $\psi_s(\omega,\cdot):X_s \to X_s$ and $\psi_u(\omega,\cdot):X_u \to X_u$
such that $\psi_s(\omega,0)=\psi_u(\omega,0)=0$ and
\begin{align}\label{lin-su}
\psi_s\circ F_s \circ \psi_s^{-1}(\omega,x)=\Lambda_s(\omega) ,\qquad
\psi_u\circ F_u \circ \psi_u^{-1}(\omega,x)=\Lambda_u(\omega).
\end{align}
Remark that although $\psi_s^{\pm 1}(\omega,\cdot)$ and $\psi_u^{\pm 1}(\omega,\cdot)$ are local mappings, a similar manner as (\ref{Es2}) helps us extend them to global ones, which are tempered $C^{1,\beta_\alpha}$ on any domain bounded by a tempered random variable.
Combining this fact with (\ref{ph-ph}) and putting
\begin{align}\label{def-Phi}
\Phi(\omega,\cdot)
:=\psi(\omega,\cdot)\circ\phi(\omega,\cdot)
\quad \mbox{for $\omega\in\Omega$}
\end{align}
with $\psi(\omega,x):=\psi_s(\omega,\pi_s x)
+\psi_u(\omega,\pi_u x)$, we obtain that
\begin{align}\label{Ph-Hold}
\|\Phi(\omega,x)-\Phi(\omega,\tilde x)\|\le L_3(\omega)\|x-\tilde x\|^{\hat\sigma},\quad \forall x,\tilde x\in U_1,
\end{align}
where $L_3(\omega)>0$ is a random variable tempered from above.

For the estimate of $\Phi^{-1}(\omega,\cdot)$, we
set $y:=\Phi(\omega,x)$ for $x\in U_1$ and notice that
\begin{align}\label{psiyy}
\|\psi_u^{-1}(\omega,\pi_u y)\|\le \hat L_\epsilon(\omega)\|\pi_u y\|
\end{align}
with an $\epsilon$-slowly varying tempered random variable
$\hat L_\epsilon(\omega)>0$
because $\psi_u^{-1}(\omega,\cdot)$ is tempered $C^{1,\beta_\alpha}$ on a bounded domain $\{y= \Phi(\omega,x):x\in U_1\}$ (the boundedness is indicated by (\ref{Ph-Hold})). Then, we conclude that $\psi^{-1}(y)\in U_\phi(\omega)$ by (\ref{def-Phi}) and therefore (\ref{ph-su}) yields
\begin{align}\label{Ph-Lip}
\|\Phi^{-1}(\omega,y)-\Phi^{-1}(\omega,\pi_s y)\|
&=\|\phi^{-1}(\omega,\cdot)\circ\psi^{-1}(\omega,y)
-\phi^{-1}(\omega,\cdot)\circ\psi_s^{-1}(\omega,\pi_s y)\|
\nonumber\\
&\le L_2(\omega)\|\psi_u^{-1}(\omega,\pi_u y)\|^{\hat \sigma}\le L_3(\omega)\|\pi_u y\|^{\hat \sigma}
\end{align}
provided that $\|\pi_u y\|<1/\hat L_\epsilon(\omega)$ (implying $\|\psi_u^{-1}(\omega,\pi_u y)\|<1$ by (\ref{psiyy})).

Moreover, for any $y=\Phi(\omega,x)$ with $x\in \mathbb{R}^d$, we also have that
\begin{align}\label{Phi-u0}
\pi_u\Phi^{-1}(\omega,\pi_s y)=\pi_u\phi^{-1}(\omega,\cdot)
\circ\psi^{-1}_s(\omega,\pi_s y)=0
\end{align}
by (\ref{uph-0out}) with $z:=\psi^{-1}_s(\omega,\pi_s y)\in X_s$,
\begin{align}\label{Phi-uout}
\pi_u\Phi^{-1}(\omega,y)=
\psi_u^{-1}(\omega,\pi_u y)\quad {\rm if}~ \|\psi_u^{-1}(\omega,\pi_u y)\|>1/2
\end{align}
by (\ref{uph-0out}) with $z:=\psi^{-1}(\omega, y)$, and that
\begin{align}\label{FA-conj}
\Phi(\theta\omega,\cdot)\circ F(\omega,\cdot)\circ\Phi^{-1}(\omega,\cdot)=\Lambda(\omega)
\end{align}
by (\ref{decp}) and (\ref{lin-su}).

Having the above preparations, we give the
following lemma.
\begin{lm}\label{lm-N}
For $U_1:=\{x\in\mathbb{R}^d:\|x\|\le 1/2\}$ and for any point $x\in U_1\backslash X_s$, there is a real $N(\omega,x)>0$ such that
\begin{align}\label{barx-out}
\|\pi_u F(n,\omega,x)\|>1/2,\quad
\forall n\ge N(\omega,x),~\forall\omega\in\Omega,
\end{align}
where $N(\omega,x)$ is measurable with respect to $\omega$.
\end{lm}

{\bf Proof of Lemma \ref{lm-N}}.
Since $\psi_u(\omega,\cdot)$ is tempered $C^{1,\beta_\alpha}$, we have
$$
\|\psi_u(\omega,x_u)\|\le \hat L_\epsilon(\omega)\|x_u\|\le (1/2)\hat L_\epsilon(\omega),\quad {\rm if}~ \|x_u\|\le 1/2,
$$
with an $\epsilon$-slowly varying tempered random variable $\hat L_\epsilon(\omega)>0$.
It means that, for any $y\in \mathbb{R}^d$,
$
\|\psi_u^{-1}(\omega,\pi_u y)\|>1/2
$
if
$
\|\pi_u y\|>\hat L_\epsilon(\omega)
$
and therefore (\ref{Phi-uout}) gives
\begin{align}\label{Ph-L}
\|\pi_u\Phi^{-1}(\omega,y)\|=
\|\psi_u^{-1}(\omega,\pi_u y)\|>1/2 \quad {\rm if}~ \|\pi_u y\|>\hat L_\epsilon(\omega).
\end{align}
By (\ref{Ph-Lip}) and (\ref{Phi-u0}), we also get
$$
\|\pi_u \Phi^{-1}(\omega,y)\|=\|\pi_u\Phi^{-1}(\omega,y)
-\pi_u\Phi^{-1}(\omega,\pi_s y)\|
\le L_3(\omega)\|\pi_u y\|^{\hat \sigma}
$$
and therefore
\begin{align}\label{uyux}
\|\pi_u y\|\ge (1/L_3(\omega))^{1/\hat \sigma}\|\pi_u \Phi^{-1}(\omega,y)\|^{1/\hat \sigma} \quad {\rm if}~ \|\pi_u y\|<1/\hat L_\epsilon(\omega).
\end{align}
Moreover, (\ref{min-u}) implies that
\begin{align}\label{est-n}
\|\Lambda_u(n,\omega)\pi_u y\|>e^{\epsilon n}\hat L_\epsilon(\omega)\ge \hat L_\epsilon(\theta^n\omega),\quad
\forall n\ge \frac{1}{\lambda_\tau-3\epsilon}\ln\frac{K_\epsilon(\omega)\hat L_\epsilon(\omega)}{\|\pi_u y\|}>0.
\end{align}

Next, for any $x\in U_1\backslash X_s$, set $y:=\Phi(\omega,x)$ and consider the case that  $\|\pi_u y\|<1/\hat L_\epsilon(\omega)$,
where
$
\|\pi_u y\|\ge (1/L_3(\omega))^{1/\hat \sigma}\|\pi_u x\|^{1/\hat \sigma}>0
$
by (\ref{uyux}). This enables us to choose
$$
N(\omega,x):=\frac{1}{\lambda_\tau-3\epsilon}\ln\frac{K_\epsilon(\omega)\hat L_\epsilon(\omega)L_3(\omega)^{1/\hat \sigma}}{\|\pi_u x\|^{1/\hat \sigma}}\ge \frac{1}{\lambda_\tau-3\epsilon}\ln\frac{K_\epsilon(\omega)\hat L_\epsilon(\omega)}{\|\pi_u y\|}>0,
$$
which is clearly measurable with respect to $\omega$.
Then, we
obtain from (\ref{FA-conj}), (\ref{Ph-L}) and (\ref{est-n}) that
\begin{align*}
\|\pi_u F(n,\omega,x)\|
=\|\pi_u F(n,\omega,\cdot) \circ \Phi^{-1}(\omega,y)\|
=
\|\pi_u \Phi^{-1}(\theta^{n}\omega,\cdot)
\circ\Lambda(n,\omega)y\|>1/2
\end{align*}
for all $n\ge N(\omega,x)$,
which proves (\ref{barx-out}) in the case of $\|\pi_u y\|<1/\hat L_\epsilon(\omega)$. In another case that $\|\pi_u y\|\ge 1/\hat L_\epsilon(\omega)$, we choose
$$
N(\omega,x):=\frac{1}{\lambda_\tau-3\epsilon}\ln\{K_\epsilon(\omega)\hat L_\epsilon(\omega)^2\}\ge \frac{1}{\lambda_\tau-3\epsilon}\ln\frac{K_\epsilon(\omega)\hat L_\epsilon(\omega)}{\|\pi_u y\|}>0,
$$
which is clearly measurable with respect to $\omega$, we still prove (\ref{barx-out}). This completes the proof of Lemma \ref{lm-N}. \qquad $\Box$

We continue to prove Theorem \ref{Main-ThmC1-RS}.
Fix an arbitrary point $\bar x\in U_1\backslash {X_s}$ and let $N(\omega,\bar x)>0$ be given in Lemma \ref{lm-N}. Then, we make a decomposition $\Omega=\cup_{m=1}^{\infty}\Omega_m$ with
\begin{align}\label{omg-m}
\Omega_m:=\{\omega\in\Omega: N(\omega,\bar x)\in (m-1,m]\},\quad \forall m\in \mathbb{N}.
\end{align}
Clearly, every $\Omega_m$ is measurable since $N(\omega,\bar x)$ is measurable with respect to $\omega$ and $(m-1,m]$ is a Borel set in $\mathbb{R}$.
By Lemma \ref{lm-N} and (\ref{omg-m}),
there is a small neighborhood $U(\omega_m, \bar x)\subset\mathbb{R}^d$ of $\bar x$ (not necessarily a tempered ball here) such that
\begin{align*}
\|\pi_u F(m,\omega_m,x)\|>1/2,\quad \forall x\in U(\omega_m, \bar x),~\forall \omega_m\in\Omega_m.
\end{align*}
This fact together with (\ref{uq0=0}) implies that
$
\pi_u q_0(\theta^m\omega_m,F(m,\omega_m,x),y_s)=0
$
for all $y_s\in X_s$,
and therefore (\ref{varp-pq}) tells that the leaf
$\mathcal W_s(\theta^{m}\omega_m,F(m,\omega_m,x))$
is the graph of the mapping
$$
y_s\mapsto \pi_u F(m,\omega_m,x).
$$
Combining this fact with (\ref{inv-fol}), we conclude that
\begin{align}\label{oth-ws}
\mathcal W_s(\omega_m,x)
&=F(-m,\theta^{m}\omega_m, \mathcal W_s(\theta^{m}\omega_m,F(m,\omega_m,x)))
\nonumber\\
&=\{F(-m,\theta^{m}\omega_m, y_s+\pi_u F(m,\omega_m,x)):y_s\in X_s\}, \quad \forall x\in U(\omega_m, \bar x).
\end{align}
Then, comparing (\ref{oth-ws}) with (\ref{varp-pq}) (replacing $y_s$ by $\pi_s F(-m,\theta^{m}\omega_m, y_s+\pi_u F(m,\omega_m,x))$), we have that
$$
\pi_u F(-m,\theta^{m}\omega_m, y_s+\pi_u F(m,\omega_m,x))=\pi_u\{x+q_0(\omega_m,x,\cdot)\circ \pi_s F(-m,\theta^{m}\omega_m, y_s+\pi_u F(m,\omega_m,x))\}.
$$
Setting $y_s=0$ in the above equality, we get by (\ref{SWSU}) that
\begin{align}\label{qyy}
\pi_u q_0(\omega_m,x,0)= \pi_u F(-m,\theta^{m}\omega_m, \pi_u F(m,\omega_m,x))-\pi_u x, \quad \forall x\in U(\omega_m, \bar x),
\end{align}
whose right-hand side clearly has a derivative measurable with respect to $\omega_m$ by our assumption that $F(\omega,\cdot)$ is a $C^{2,\alpha}$ random diffeomorphism, i.e., $F^{\pm 1}(\omega,\cdot)$ together with their derivatives are measurable with respect to $\omega$. Hence, $\partial_x(\pi_u q_0)(\omega_m,x,0)$ is measurable with respect to $\omega_m\in \Omega_m$. Since $\Omega$ is a union of countably many measurable sets $\Omega_m$, we further conclude that $\partial_x(\pi_u q_0)(\omega,x,0)$ is measurable with respect to $\omega\in\Omega$ for $x\in U(\omega,\bar x)$. Namely, in the case of $\bar x\in U_1\backslash{X_s}$, we have that
\begin{itemize}
\item[(C)]
\mbox{
$\partial_x(\pi_u q_0)(\omega,\bar x,0)$ is measurable with respect to $\omega\in \Omega$.
}
\end{itemize}

Next, in the case of $\bar x\in {X_s}\cap U_1$, we find a sequence $\{\bar x_n\}_{n\in\mathbb{N}}$ in $U_1\backslash {X_s}$ such that
$$
\lim_{n\to \infty}\bar x_n=\bar x.
$$
Since $q_0(\omega,\cdot,0)$ is continuously differentiable (actually $C^{1,\beta_\alpha}$) in $\mathbb{R}^d$, as mentioned at the beginning of this proof, we get
$$
\lim_{n\to \infty}\partial_x(\pi_u q_0)(\omega,\bar x_n,0)=\partial_x(\pi_u q_0)(\omega,\bar x,0),\quad \forall \omega\in \Omega,
$$
which together with (C) implies that $\partial_x(\pi_u q_0)(\omega,x,0)$ is measurable with respect to $\omega$ for $x\in U_1$. Similarly, we also prove that
$\partial_x(\pi_s p_0)(\omega,x,0)$ is measurable with respect to $\omega$ for $x\in U_1$.

Hence, $D\phi(\omega,x)$ is measurable with respect to $\omega\in \Omega$ by (\ref{def-phi++}). Moreover, noting that $\phi^{-1}(\omega,x)$ is a random mapping (see Lemma \ref{lm-decp}) and that
$$
D\phi^{-1}(\omega,x)=\{D\phi(\omega,\phi^{-1}(\omega,x))\}^{-1},
$$
we obtain that $D\phi^{-1}(\omega,x)$ is measurable with respect to $\omega$.
It follows that the mapping $\Phi(\omega,\cdot)$ defined in (\ref{def-Phi}) is a $C^{1,\beta_\alpha}$ random diffeomorphism in $U_1\supset U(\omega)$, which proves (ii) of Theorem \ref{Main-ThmC1-RS}. Since
$\Phi(\omega,\cdot)$ straightens up the stable and unstable foliations by (\ref{def-Phi}), we also prove
(i) of Theorem \ref{Main-ThmC1-RS} and
the proof is completed. \qquad $\Box$

\section*{Appendix}
\setcounter{equation}{0}
\renewcommand\thesection{A}

In this appendix,
we prove (\ref{cutoff-w}).
Putting $\hat \rho(\omega):=\rho(\omega)^2>0$, we define
\begin{align*}
g(\omega,s):=
\left\{
  \begin{array}{lll}
   e^{\,\hat \rho(\omega)^2(s-\hat \rho(\omega)/2)^{-1}(s-\hat \rho(\omega))^{-1}}, & \hat \rho(\omega)/2< s< \hat \rho(\omega),
   \\
   0, & \mbox{other cases},
  \end{array}
\right.
\end{align*}
and
\begin{align*}
{h}(\omega,s):=\frac{\int_{s}^{+\infty}g(\omega,\varsigma)d\varsigma}
{\int_{-\infty}
^{+\infty}g(\omega,\varsigma)d\varsigma}.
\end{align*}
It is clear that $h(\omega,\cdot):\mathbb{R}\to (0,1)$ is
smooth such that $h(\omega,s)=1$ for all $s\le \hat \rho(\omega)/2$ and $h(\omega,s)=0$ for all $s\ge \hat \rho(\omega)$. For the estimate of its derivative,
we have that the maximum value of $g(\omega,s)$ is $e^{-16}$ for all $\omega$ and $s$. Computing the points at which the values of $g(\omega,\cdot)$ are $e^{-16}/2$, i.e.,
$$
\frac{3}{4}\hat \rho(\omega)\pm \bigg(\frac{1}{16}-\frac{1}{\ln2+16}\bigg)^{1/2}
\hat \rho(\omega) \in (\hat \rho(\omega)/2, \hat \rho(\omega)),
$$
we have the  estimates
\begin{align*}
\frac{1}{\hat \rho(\omega)}\int_{0}
^{+\infty}g(\omega,\varsigma)d\varsigma\ge \bigg(\frac{1}{16}-\frac{1}{\ln2+16}\bigg)^{1/2}
e^{-16},
\end{align*}
and therefore
\begin{align}
\hat \rho(\omega)|Dh(\omega,s)|\le \frac{\hat \rho(\omega) |g(\omega,s)|}{\int_{0}
^{+\infty}g(\omega,\varsigma)d\varsigma}\le \bigg(\frac{1}{16}-\frac{1}{\ln2+16}\bigg)^{-1/2}.
\label{h-1}
\end{align}
For the second order derivative of $h(\omega,\cdot)$, noting that
$
\max_{a\ge 0}\{e^{-a}a^r\}\le e^{-r}r^r
$
for $r\ge 1$, we have
\begin{align*}
&\hat \rho(\omega)|Dg(\omega,s)|
\\
&\le e^{\,\hat \rho(\omega)^2(s-\hat \rho(\omega)/2)^{-1}(s-\hat \rho(\omega))^{-1}}
\frac{\hat \rho(\omega)^3|2s-3\hat \rho(\omega)/2|}{(s-\hat \rho(\omega)/2)^2(s-\hat \rho(\omega))^2}
\\
&\le \bigg\{e^{\,\hat \rho(\omega)^2(s-\hat \rho(\omega)/2)^{-1}(s-\hat \rho(\omega))^{-1}}
\frac{\hat \rho(\omega)^4}{(s-\hat \rho(\omega)/2)^2(s-\hat \rho(\omega))^2}
\bigg\}\frac{|2s-3\hat \rho(\omega)/2|}{\hat \rho(\omega)}
\\
&\le 2e^{-2}
\end{align*}
for $s\in(\hat \rho(\omega)/2, \,\hat \rho(\omega))$, and
therefore $\hat \rho(\omega)|Dg(\omega,s)|\le 2e^{-2}$ for all
$s\in\mathbb{R}$ since $g(\omega,s)=0$ outside $(\hat \rho(\omega)/2, \hat \rho(\omega))$. Thus,
\begin{align}
\hat \rho(\omega)^2|D^2{h}(\omega,s)|\le \frac{\hat \rho(\omega)^2 |Dg(\omega,s)|}{\int_{0}
^{+\infty}g(\omega,\varsigma)d\varsigma}\le 2\bigg(\frac{1}{16}-\frac{1}{\ln2+16}\bigg)^{-1/2}e^{14}.
\label{h-2}
\end{align}
Similarly, for the third order derivative we have the estimate
\begin{align}
\hat \rho(\omega)^3|D^3{h}(\omega,s)|\le (4e^{-4}+(9/2)e^{-3}+8e^{-2})
\bigg(\frac{1}{16}-\frac{1}{\ln2+16}\bigg)^{-1/2}e^{16}
\label{h-3}
\end{align}
for all $s\in \mathbb{R}$.

Next, let $U(\omega)$ and $V(\omega)$ be neighborhoods of the origin with diameters $\rho(\omega)/2$ and $\rho(\omega)$ respectively, and define $u(\omega,\cdot):\mathbb{R}^d\to (0,1)$ as
\begin{align*}
u(\omega,x):=h\Big(\omega, \sum_{j=1,...,{p}}\,\sum_{\iota=1,...,d_j}(\pi_{j,\iota}x)^2\Big).
\end{align*}
Then we verify that
\begin{align*}
u(\omega,x)=
\left\{
\begin{array}{lll}
1,& \forall x\in U(\omega),
\vspace{1ex}\\
0,& \forall x\in \mathbb{R}^d\backslash V(\omega).
\end{array}
\right.
\end{align*}
Moreover, (\ref{h-1})-(\ref{h-3}) give
\begin{align*}
\rho^{r}(\omega)\|D^ru(\omega,x)\|
\le 2r\rho^{2r}(\omega)\Big|D^rh\Big(\sum_{j=1,...,{p}}\,\sum_{\iota=1,...,d_j}(\pi_{j,\iota}x)^2\Big)\Big|\le C_u,~~~r=1,2,3,
\end{align*}
where
$$
C_u:=6(4e^{-4}+(9/2)e^{-3}+8e^{-2})
\bigg(\frac{1}{16}-\frac{1}{\ln2+16}\bigg)^{-1/2}e^{16}>0.
$$
This proves (\ref{cutoff-w}).


\end{document}